\def\secl{1}
\def\seg{2}
\def\ser{3}
\def\sewz{4}
\def\seex{5}
\def\seoe{6}
\def\tU{\tilde{U}}
\def\tX{\tilde{X}}
\def\tx{\tilde{x}}
\def\tP{\tilde{P}}

\def\Se{Sasakian-Einstein }

\expandafter\chardef\csname pre amssym.def at\endcsname=\the\catcode`\@
\catcode`\@=11

\def\undefine#1{\let#1\undefined}
\def\newsymbol#1#2#3#4#5{\let\next@\relax
 \ifnum#2=\@ne\let\next@\msafam@\else
 \ifnum#2=\tw@\let\next@\msbfam@\fi\fi
 \mathchardef#1="#3\next@#4#5}
\def\mathhexbox@#1#2#3{\relax
 \ifmmode\mathpalette{}{\m@th\mathchar"#1#2#3}%
 \else\leavevmode\hbox{$\m@th\mathchar"#1#2#3$}\fi}
\def\hexnumber@#1{\ifcase#1 0\or 1\or 2\or 3\or 4\or 5\or 6\or 7\or 8\or
 9\or A\or B\or C\or D\or E\or F\fi}

\font\tenmsa=msam10
\font\sevenmsa=msam7
\font\fivemsa=msam5
\newfam\msafam
\textfont\msafam=\tenmsa
\scriptfont\msafam=\sevenmsa
\scriptscriptfont\msafam=\fivemsa
\edef\msafam@{\hexnumber@\msafam}
\mathchardef\dabar@"0\msafam@39
\def\dashrightarrow{\mathrel{\dabar@\dabar@\mathchar"0\msafam@4B}}
\def\dashleftarrow{\mathrel{\mathchar"0\msafam@4C\dabar@\dabar@}}

\def\ulcorner{\delimiter"4\msafam@70\msafam@70 }
\def\urcorner{\delimiter"5\msafam@71\msafam@71 }
\def\llcorner{\delimiter"4\msafam@78\msafam@78 }
\def\lrcorner{\delimiter"5\msafam@79\msafam@79 }
\def\yen{{\mathhexbox@\msafam@55 }}
\def\checkmark{{\mathhexbox@\msafam@58 }}
\def\circledR{{\mathhexbox@\msafam@72 }}
\def\maltese{{\mathhexbox@\msafam@7A }}

\font\tenmsb=msbm10
\font\sevenmsb=msbm7
\font\fivemsb=msbm5
\newfam\msbfam
\textfont\msbfam=\tenmsb
\scriptfont\msbfam=\sevenmsb
\scriptscriptfont\msbfam=\fivemsb
\edef\msbfam@{\hexnumber@\msbfam}

\catcode`\@=\csname pre amssym.def at\endcsname

\expandafter\ifx\csname pre amssym.tex at\endcsname\relax \else \endinput\fi
\expandafter\chardef\csname pre amssym.tex at\endcsname=\the\catcode`\@
\catcode`\@=11
\newsymbol\boxdot 1200
\newsymbol\boxplus 1201
\newsymbol\boxtimes 1202
\newsymbol\square 1003
\newsymbol\blacksquare 1004
\newsymbol\centerdot 1205
\newsymbol\lozenge 1006
\newsymbol\blacklozenge 1007
\newsymbol\circlearrowright 1308
\newsymbol\circlearrowleft 1309
\undefine\rightleftharpoons
\newsymbol\rightleftharpoons 130A
\newsymbol\leftrightharpoons 130B
\newsymbol\boxminus 120C
\newsymbol\Vdash 130D
\newsymbol\Vvdash 130E
\newsymbol\vDash 130F
\newsymbol\twoheadrightarrow 1310
\newsymbol\twoheadleftarrow 1311
\newsymbol\leftleftarrows 1312
\newsymbol\rightrightarrows 1313
\newsymbol\upuparrows 1314
\newsymbol\downdownarrows 1315
\newsymbol\upharpoonright 1316
 
\newsymbol\downharpoonright 1317
\newsymbol\upharpoonleft 1318
\newsymbol\downharpoonleft 1319
\newsymbol\rightarrowtail 131A
\newsymbol\leftarrowtail 131B
\newsymbol\leftrightarrows 131C
\newsymbol\rightleftarrows 131D
\newsymbol\Lsh 131E
\newsymbol\Rsh 131F
\newsymbol\rightsquigarrow 1320
\newsymbol\leftrightsquigarrow 1321
\newsymbol\looparrowleft 1322
\newsymbol\looparrowright 1323
\newsymbol\circeq 1324
\newsymbol\succsim 1325
\newsymbol\gtrsim 1326
\newsymbol\gtrapprox 1327
\newsymbol\multimap 1328
\newsymbol\therefore 1329
\newsymbol\because 132A
\newsymbol\doteqdot 132B
 
\newsymbol\triangleq 132C
\newsymbol\precsim 132D
\newsymbol\lesssim 132E
\newsymbol\lessapprox 132F
\newsymbol\eqslantless 1330
\newsymbol\eqslantgtr 1331
\newsymbol\curlyeqprec 1332
\newsymbol\curlyeqsucc 1333
\newsymbol\preccurlyeq 1334
\newsymbol\leqq 1335
\newsymbol\leqslant 1336
\newsymbol\lessgtr 1337
\newsymbol\backprime 1038
\newsymbol\risingdotseq 133A
\newsymbol\fallingdotseq 133B
\newsymbol\succcurlyeq 133C
\newsymbol\geqq 133D
\newsymbol\geqslant 133E
\newsymbol\gtrless 133F
\newsymbol\sqsubset 1340
\newsymbol\sqsupset 1341
\newsymbol\vartriangleright 1342
\newsymbol\vartriangleleft 1343
\newsymbol\trianglerighteq 1344
\newsymbol\trianglelefteq 1345
\newsymbol\bigstar 1046
\newsymbol\between 1347
\newsymbol\blacktriangledown 1048
\newsymbol\blacktriangleright 1349
\newsymbol\blacktriangleleft 134A
\newsymbol\vartriangle 134D
\newsymbol\blacktriangle 104E
\newsymbol\triangledown 104F
\newsymbol\eqcirc 1350
\newsymbol\lesseqgtr 1351
\newsymbol\gtreqless 1352
\newsymbol\lesseqqgtr 1353
\newsymbol\gtreqqless 1354
\newsymbol\Rrightarrow 1356
\newsymbol\Lleftarrow 1357
\newsymbol\veebar 1259
\newsymbol\barwedge 125A
\newsymbol\doublebarwedge 125B
\undefine\angle
\newsymbol\angle 105C
\newsymbol\measuredangle 105D
\newsymbol\sphericalangle 105E
\newsymbol\varpropto 135F
\newsymbol\smallsmile 1360
\newsymbol\smallfrown 1361
\newsymbol\Subset 1362
\newsymbol\Supset 1363
\newsymbol\Cup 1264
 
\newsymbol\Cap 1265
 
\newsymbol\curlywedge 1266
\newsymbol\curlyvee 1267
\newsymbol\leftthreetimes 1268
\newsymbol\rightthreetimes 1269
\newsymbol\subseteqq 136A
\newsymbol\supseteqq 136B
\newsymbol\bumpeq 136C
\newsymbol\Bumpeq 136D
\newsymbol\lll 136E
 
\newsymbol\ggg 136F
 
\newsymbol\circledS 1073
\newsymbol\pitchfork 1374
\newsymbol\dotplus 1275
\newsymbol\backsim 1376
\newsymbol\backsimeq 1377
\newsymbol\complement 107B
\newsymbol\intercal 127C
\newsymbol\circledcirc 127D
\newsymbol\circledast 127E
\newsymbol\circleddash 127F
\newsymbol\lvertneqq 2300
\newsymbol\gvertneqq 2301
\newsymbol\nleq 2302
\newsymbol\ngeq 2303
\newsymbol\nless 2304
\newsymbol\ngtr 2305
\newsymbol\nprec 2306
\newsymbol\nsucc 2307
\newsymbol\lneqq 2308
\newsymbol\gneqq 2309
\newsymbol\nleqslant 230A
\newsymbol\ngeqslant 230B
\newsymbol\lneq 230C
\newsymbol\gneq 230D
\newsymbol\npreceq 230E
\newsymbol\nsucceq 230F
\newsymbol\precnsim 2310
\newsymbol\succnsim 2311
\newsymbol\lnsim 2312
\newsymbol\gnsim 2313
\newsymbol\nleqq 2314
\newsymbol\ngeqq 2315
\newsymbol\precneqq 2316
\newsymbol\succneqq 2317
\newsymbol\precnapprox 2318
\newsymbol\succnapprox 2319
\newsymbol\lnapprox 231A
\newsymbol\gnapprox 231B
\newsymbol\nsim 231C
\newsymbol\ncong 231D
\newsymbol\diagup 231E
\newsymbol\diagdown 231F
\newsymbol\varsubsetneq 2320
\newsymbol\varsupsetneq 2321
\newsymbol\nsubseteqq 2322
\newsymbol\nsupseteqq 2323
\newsymbol\subsetneqq 2324
\newsymbol\supsetneqq 2325
\newsymbol\varsubsetneqq 2326
\newsymbol\varsupsetneqq 2327
\newsymbol\subsetneq 2328
\newsymbol\supsetneq 2329
\newsymbol\nsubseteq 232A
\newsymbol\nsupseteq 232B
\newsymbol\nparallel 232C
\newsymbol\nmid 232D
\newsymbol\nshortmid 232E
\newsymbol\nshortparallel 232F
\newsymbol\nvdash 2330
\newsymbol\nVdash 2331
\newsymbol\nvDash 2332
\newsymbol\nVDash 2333
\newsymbol\ntrianglerighteq 2334
\newsymbol\ntrianglelefteq 2335
\newsymbol\ntriangleleft 2336
\newsymbol\ntriangleright 2337
\newsymbol\nleftarrow 2338
\newsymbol\nrightarrow 2339
\newsymbol\nLeftarrow 233A
\newsymbol\nRightarrow 233B
\newsymbol\nLeftrightarrow 233C
\newsymbol\nleftrightarrow 233D
\newsymbol\divideontimes 223E
\newsymbol\varnothing 203F
\newsymbol\nexists 2040
\newsymbol\Finv 2060
\newsymbol\Game 2061
\newsymbol\mho 2066
\newsymbol\eth 2067
\newsymbol\eqsim 2368
\newsymbol\beth 2069
\newsymbol\gimel 206A
\newsymbol\daleth 206B
\newsymbol\lessdot 236C
\newsymbol\gtrdot 236D
\newsymbol\ltimes 226E
\newsymbol\rtimes 226F
\newsymbol\shortmid 2370
\newsymbol\shortparallel 2371
\newsymbol\smallsetminus 2272
\newsymbol\thicksim 2373
\newsymbol\thickapprox 2374
\newsymbol\approxeq 2375
\newsymbol\succapprox 2376
\newsymbol\precapprox 2377
\newsymbol\curvearrowleft 2378
\newsymbol\curvearrowright 2379
\newsymbol\digamma 207A
\newsymbol\varkappa 207B
\newsymbol\Bbbk 207C
\newsymbol\hslash 207D
\undefine\hbar
\newsymbol\hbar 207E
\newsymbol\backepsilon 237F
\catcode`\@=\csname pre amssym.tex at\endcsname

\magnification=1200
\hsize=468truept
\vsize=646truept
\voffset=-10pt
\parskip=4pt
\baselineskip=14truept
\count0=1

\dimen100=\hsize

\def\leftill#1#2#3#4{
\medskip
\line{$
\vcenter{
\hsize = #1truept \hrule\hbox{\vrule\hbox to  \hsize{\hss \vbox{\vskip#2truept
\hbox{{\copy100 \the\count105}: #3}\vskip2truept}\hss }
\vrule}\hrule}
\dimen110=\dimen100
\advance\dimen110 by -36truept
\advance\dimen110 by -#1truept
\hss \vcenter{\hsize = \dimen110
\medskip
\noindent { #4\par\medskip}}$}
\advance\count105 by 1
}
\def\rightill#1#2#3#4{
\medskip
\line{
\dimen110=\dimen100
\advance\dimen110 by -36truept
\advance\dimen110 by -#1truept
$\vcenter{\hsize = \dimen110
\medskip
\noindent { #4\par\medskip}}
\hss \vcenter{
\hsize = #1truept \hrule\hbox{\vrule\hbox to  \hsize{\hss \vbox{\vskip#2truept
\hbox{{\copy100 \the\count105}: #3}\vskip2truept}\hss }
\vrule}\hrule}
$}
\advance\count105 by 1
}
\def\midill#1#2#3{\medskip
\line{$\hss
\vcenter{
\hsize = #1truept \hrule\hbox{\vrule\hbox to  \hsize{\hss \vbox{\vskip#2truept
\hbox{{\copy100 \the\count105}: #3}\vskip2truept}\hss }
\vrule}\hrule}
\dimen110=\dimen100
\advance\dimen110 by -36truept
\advance\dimen110 by -#1truept
\hss $}
\advance\count105 by 1
}
\def\insectnum{\copy110\the\count120
\advance\count120 by 1
}

\font\ninerm=cmr9

\font\tenrm=cmr10 at 10pt

\font\sc=cmcsc10

\def\msb{\fam\msbfam\tenmsb}

\def\bba{{\msb A}}

\def\bbc{{\msb C}}

\def\bbp{{\msb P}}
\def\bbq{{\msb Q}}
\def\bbr{{\msb R}}

\def\bbz{{\msb Z}}

\def\grG{\Gamma}

\def\grO{\Omega}

\def\gra{\alpha}
\def\grb{\beta}

\def\grd{\delta}

\def\grg{\gamma}

\def\grk{\kappa}

\def\gro{\omega}

\def\grr{\rho}
\def\grs{\sigma}
\def\grt{\tau}

\font\svtnrm=cmr17

\font\aa=eufm10

\def\gs{{\Got s}}

\def\gu{{\Got u}}

\def\gt{{\Got t}}
\def\Got#1{\hbox{\aa#1}}

\def\gsp1{{\Got s}{\Got p}(1)}

\def\mh-1{\hat{\mu}^{-1}(0)}
\def\n-1c{\nu^{-1}(c)}
\def\m-1{\mu^{-1}(0)}
\def\p-1{\pi^{-1}}
\def\p'-1{\prime{\pi}^{-1}}

\def\cala{{\cal A}}

\def\cale{{\cal E}}
\def\calf{{\cal F}}
\def\calg{{\cal G}}
\def\calh{{\cal H}}

\def\call{{\cal L}}
\def\calm{{\cal M}}

\def\calr{{\cal R}}
\def\cals{{\cal S}}

\def\calz{{\cal Z}}

\def\la#1{\hbox to #1pc{\leftarrowfill}}
\def\ra#1{\hbox to #1pc{\rightarrowfill}}

\def\fract#1#2{\raise4pt\hbox{$ #1 \atop #2 $}}
\def\decdnar#1{\phantom{\hbox{$\scriptstyle{#1}$}}
\left\downarrow\vbox{\vskip15pt\hbox{$\scriptstyle{#1}$}}\right.}

\def\bowtie{\hbox to 1pt{\hss}\raise.66pt\hbox{$\scriptstyle{>}$}
\kern-4.9pt\triangleleft}
\def\hsmash{\triangleright\kern-4.4pt\raise.66pt\hbox{$\scriptstyle{<}$}}
\def\boxit#1{\vbox{\hrule\hbox{\vrule\kern3pt
\vbox{\kern3pt#1\kern3pt}\kern3pt\vrule}\hrule}}

\def\za{\vrule height6pt width4pt depth1pt}

\font\aa=eufm10

\def\Got#1{\hbox{\aa#1}}

\def\cala{{\cal A}}

\def\cale{{\cal E}}
\def\calf{{\cal F}}
\def\calg{{\cal G}}
\def\calh{{\cal H}}

\def\call{{\cal L}}
\def\calm{{\cal M}}

\def\calr{{\cal R}}
\def\cals{{\cal S}}

\def\calz{{\cal Z}}

\def\ga{{\Got a}}

\def\gs{{\Got s}}
\def\gt{{\Got t}}
\def\gu{{\Got u}}

\catcode`!=11 
 
  

\def\PiC{P\kern-.12em\lower.5ex\hbox{I}\kern-.075emC}
\def\PiCTeX{\PiC\kern-.11em\TeX}

\def\!ifnextchar#1#2#3{%
  \let\!testchar=#1%
  \def\!first{#2}%
  \def\!second{#3}%
  \futurelet\!nextchar\!testnext}
\def\!testnext{%
  \ifx \!nextchar \!spacetoken 
    \let\!next=\!skipspacetestagain
  \else
    \ifx \!nextchar \!testchar
      \let\!next=\!first
    \else 
      \let\!next=\!second 
    \fi 
  \fi
  \!next}
\def\\{\!skipspacetestagain} 
  \expandafter\def\\ {\futurelet\!nextchar\!testnext} 
\def\\{\let\!spacetoken= } \\  

\def\!tfor#1:=#2\do#3{%
  \edef\!fortemp{#2}%
  \ifx\!fortemp\!empty 
    \else
    \!tforloop#2\!nil\!nil\!!#1{#3}%
  \fi}
\def\!tforloop#1#2\!!#3#4{%
  \def#3{#1}%
  \ifx #3\!nnil
    \let\!nextwhile=\!fornoop
  \else
    #4\relax
    \let\!nextwhile=\!tforloop
  \fi 
  \!nextwhile#2\!!#3{#4}}

\def\!etfor#1:=#2\do#3{%
  \def\!!tfor{\!tfor#1:=}%
  \edef\!!!tfor{#2}%
  \expandafter\!!tfor\!!!tfor\do{#3}}

\def\!cfor#1:=#2\do#3{%
  \edef\!fortemp{#2}%
  \ifx\!fortemp\!empty 
  \else
    \!cforloop#2,\!nil,\!nil\!!#1{#3}%
  \fi}
\def\!cforloop#1,#2\!!#3#4{%
  \def#3{#1}%
  \ifx #3\!nnil
    \let\!nextwhile=\!fornoop 
  \else
    #4\relax
    \let\!nextwhile=\!cforloop
  \fi
  \!nextwhile#2\!!#3{#4}}

\def\!ecfor#1:=#2\do#3{%
  \def\!!cfor{\!cfor#1:=}%
  \edef\!!!cfor{#2}%
  \expandafter\!!cfor\!!!cfor\do{#3}}

\def\!empty{}
\def\!nnil{\!nil}
\def\!fornoop#1\!!#2#3{}

\def\!ifempty#1#2#3{%
  \edef\!emptyarg{#1}%
  \ifx\!emptyarg\!empty
    #2%
  \else
    #3%
  \fi}
 
\def\!getnext#1\from#2{%
  \expandafter\!gnext#2\!#1#2}%
\def\!gnext\\#1#2\!#3#4{%
  \def#3{#1}%
  \def#4{#2\\{#1}}%
  \ignorespaces}

%
\def\!getnextvalueof#1\from#2{%
  \expandafter\!gnextv#2\!#1#2}%
\def\!gnextv\\#1#2\!#3#4{%
  #3=#1%
  \def#4{#2\\{#1}}%
  \ignorespaces}

\def\!copylist#1\to#2{%
  \expandafter\!!copylist#1\!#2}
\def\!!copylist#1\!#2{%
  \def#2{#1}\ignorespaces}

\def\!wlet#1=#2{%
  \let#1=#2 
  \wlog{\string#1=\string#2}}
 
\def\!listaddon#1#2{%
  \expandafter\!!listaddon#2\!{#1}#2}
\def\!!listaddon#1\!#2#3{%
  \def#3{#1\\#2}}
 

\def\!rightappend#1\withCS#2\to#3{\expandafter\!!rightappend#3\!#2{#1}#3}
\def\!!rightappend#1\!#2#3#4{\def#4{#1#2{#3}}}

\def\!leftappend#1\withCS#2\to#3{\expandafter\!!leftappend#3\!#2{#1}#3}
\def\!!leftappend#1\!#2#3#4{\def#4{#2{#3}#1}}

\def\!lop#1\to#2{\expandafter\!!lop#1\!#1#2}
\def\!!lop\\#1#2\!#3#4{\def#4{#1}\def#3{#2}}



\def\!loop#1\repeat{\def\!body{#1}\!iterate}
\def\!iterate{\!body\let\!next=\!iterate\else\let\!next=\relax\fi\!next}
 
\def\!!loop#1\repeat{\def\!!body{#1}\!!iterate}
\def\!!iterate{\!!body\let\!!next=\!!iterate\else\let\!!next=\relax\fi\!!next}
 
\def\!removept#1#2{\edef#2{\expandafter\!!removePT\the#1}}
{\catcode`p=12 \catcode`t=12 \gdef\!!removePT#1pt{#1}}

\def\placevalueinpts of <#1> in #2 {%
  \!removept{#1}{#2}}
 
\def\!mlap#1{\hbox to 0pt{\hss#1\hss}}
\def\!vmlap#1{\vbox to 0pt{\vss#1\vss}}
 
\def\!not#1{%
  #1\relax
    \!switchfalse
  \else
    \!switchtrue
  \fi
  \if!switch
  \ignorespaces}


 

\let\!!!wlog=\wlog              
\def\wlog#1{}    

\newdimen\headingtoplotskip     
\newdimen\linethickness         
\newdimen\longticklength        
\newdimen\plotsymbolspacing     
\newdimen\shortticklength       
\newdimen\stackleading          
\newdimen\tickstovaluesleading  
\newdimen\totalarclength        
\newdimen\valuestolabelleading  

\newbox\!boxA                   
\newbox\!boxB                   
\newbox\!picbox                 
\newbox\!plotsymbol             
\newbox\!putobject              
\newbox\!shadesymbol            

\newcount\!countA               
\newcount\!countB               
\newcount\!countC               
\newcount\!countD               
\newcount\!countE               
\newcount\!countF               
\newcount\!countG               
\newcount\!fiftypt              
\newcount\!intervalno           
\newcount\!npoints              
\newcount\!nsegments            
\newcount\!ntemp                
\newcount\!parity               
\newcount\!scalefactor          
\newcount\!tfs                  
\newcount\!tickcase             

\newdimen\!Xleft                
\newdimen\!Xright               
\newdimen\!Xsave                
\newdimen\!Ybot                 
\newdimen\!Ysave                
\newdimen\!Ytop                 
\newdimen\!angle                
\newdimen\!arclength            
\newdimen\!areabloc             
\newdimen\!arealloc             
\newdimen\!arearloc             
\newdimen\!areatloc             
\newdimen\!bshrinkage           
\newdimen\!checkbot             
\newdimen\!checkleft            
\newdimen\!checkright           
\newdimen\!checktop             
\newdimen\!dimenA               
\newdimen\!dimenB               
\newdimen\!dimenC               
\newdimen\!dimenD               
\newdimen\!dimenE               
\newdimen\!dimenF               
\newdimen\!dimenG               
\newdimen\!dimenH               
\newdimen\!dimenI               
\newdimen\!distacross           
\newdimen\!downlength           
\newdimen\!dp                   
\newdimen\!dshade               
\newdimen\!dxpos                
\newdimen\!dxprime              
\newdimen\!dypos                
\newdimen\!dyprime              
\newdimen\!ht                   
\newdimen\!leaderlength         
\newdimen\!lshrinkage           
\newdimen\!midarclength         
\newdimen\!offset               
\newdimen\!plotheadingoffset    
\newdimen\!plotsymbolxshift     
\newdimen\!plotsymbolyshift     
\newdimen\!plotxorigin          
\newdimen\!plotyorigin          
\newdimen\!rootten              
\newdimen\!rshrinkage           
\newdimen\!shadesymbolxshift    
\newdimen\!shadesymbolyshift    
\newdimen\!tenAa                
\newdimen\!tenAc                
\newdimen\!tenAe                
\newdimen\!tshrinkage           
\newdimen\!uplength             
\newdimen\!wd                   
\newdimen\!wmax                 
\newdimen\!wmin                 
\newdimen\!xB                   
\newdimen\!xC                   
\newdimen\!xE                   
\newdimen\!xM                   
\newdimen\!xS                   
\newdimen\!xaxislength          
\newdimen\!xdiff                
\newdimen\!xleft                
\newdimen\!xloc                 
\newdimen\!xorigin              
\newdimen\!xpivot               
\newdimen\!xpos                 
\newdimen\!xprime               
\newdimen\!xright               
\newdimen\!xshade               
\newdimen\!xshift               
\newdimen\!xtemp                
\newdimen\!xunit                
\newdimen\!xxE                  
\newdimen\!xxM                  
\newdimen\!xxS                  
\newdimen\!xxloc                
\newdimen\!yB                   
\newdimen\!yC                   
\newdimen\!yE                   
\newdimen\!yM                   
\newdimen\!yS                   
\newdimen\!yaxislength          
\newdimen\!ybot                 
\newdimen\!ydiff                
\newdimen\!yloc                 
\newdimen\!yorigin              
\newdimen\!ypivot               
\newdimen\!ypos                 
\newdimen\!yprime               
\newdimen\!yshade               
\newdimen\!yshift               
\newdimen\!ytemp                
\newdimen\!ytop                 
\newdimen\!yunit                
\newdimen\!yyE                  
\newdimen\!yyM                  
\newdimen\!yyS                  
\newdimen\!yyloc                
\newdimen\!zpt                  

\newif\if!axisvisible           
\newif\if!gridlinestoo          
\newif\if!keepPO                
\newif\if!placeaxislabel        
\newif\if!switch                
\newif\if!xswitch               

\newtoks\!axisLaBeL             
\newtoks\!keywordtoks           

\newwrite\!replotfile           

\newhelp\!keywordhelp{The keyword mentioned in the error message in unknown. 
Replace NEW KEYWORD in the indicated response by the keyword that 
should have been specified.}    

\!wlet\!!origin=\!xM                   
\!wlet\!!unit=\!uplength               
\!wlet\!Lresiduallength=\!dimenG       
\!wlet\!Rresiduallength=\!dimenF       
\!wlet\!axisLength=\!distacross        
\!wlet\!axisend=\!ydiff                
\!wlet\!axisstart=\!xdiff              
\!wlet\!axisxlevel=\!arclength         
\!wlet\!axisylevel=\!downlength        
\!wlet\!beta=\!dimenE                  
\!wlet\!gamma=\!dimenF                 
\!wlet\!shadexorigin=\!plotxorigin     
\!wlet\!shadeyorigin=\!plotyorigin     
\!wlet\!ticklength=\!xS                
\!wlet\!ticklocation=\!xE              
\!wlet\!ticklocationincr=\!yE          
\!wlet\!tickwidth=\!yS                 
\!wlet\!totalleaderlength=\!dimenE     
\!wlet\!xone=\!xprime                  
\!wlet\!xtwo=\!dxprime                 
\!wlet\!ySsave=\!yM                    
\!wlet\!ybB=\!yB                       
\!wlet\!ybC=\!yC                       
\!wlet\!ybE=\!yE                       
\!wlet\!ybM=\!yM                       
\!wlet\!ybS=\!yS                       
\!wlet\!ybpos=\!yyloc                  
\!wlet\!yone=\!yprime                  
\!wlet\!ytB=\!xB                       
\!wlet\!ytC=\!xC                       
\!wlet\!ytE=\!downlength               
\!wlet\!ytM=\!arclength                
\!wlet\!ytS=\!distacross               
\!wlet\!ytpos=\!xxloc                  
\!wlet\!ytwo=\!dyprime                 

\!zpt=0pt                              
\!xunit=1pt
\!yunit=1pt
\!arearloc=\!xunit
\!areatloc=\!yunit
\!dshade=5pt
\!leaderlength=24in
\!tfs=256                              
\!wmax=5.3pt                           
\!wmin=2.7pt                           
\!xaxislength=\!xunit
\!xpivot=\!zpt
\!yaxislength=\!yunit 
\!ypivot=\!zpt
\plotsymbolspacing=.4pt
  \!dimenA=50pt \!fiftypt=\!dimenA     

\!rootten=3.162278pt                   
\!tenAa=8.690286pt                     
\!tenAc=2.773839pt                     
\!tenAe=2.543275pt                     

\def\!cosrotationangle{1}      
\def\!sinrotationangle{0}      
\def\!xpivotcoord{0}           
\def\!xref{0}                  
\def\!xshadesave{0}            
\def\!ypivotcoord{0}           
\def\!yref{0}                  
\def\!yshadesave{0}            
\def\!zero{0}                  

\let\wlog=\!!!wlog
%
  
\def\normalgraphs{%
  \longticklength=.4\baselineskip
  \shortticklength=.25\baselineskip
  \tickstovaluesleading=.25\baselineskip
  \valuestolabelleading=.8\baselineskip
  \linethickness=.4pt
  \stackleading=.17\baselineskip
  \headingtoplotskip=1.5\baselineskip
  \visibleaxes
  \ticksout
  \nogridlines
  \unloggedticks}
%
\def\setplotarea x from #1 to #2, y from #3 to #4 {%
  \!arealloc=\!M{#1}\!xunit \advance \!arealloc -\!xorigin
  \!areabloc=\!M{#3}\!yunit \advance \!areabloc -\!yorigin
  \!arearloc=\!M{#2}\!xunit \advance \!arearloc -\!xorigin
  \!areatloc=\!M{#4}\!yunit \advance \!areatloc -\!yorigin
  \!initinboundscheck
  \!xaxislength=\!arearloc  \advance\!xaxislength -\!arealloc
  \!yaxislength=\!areatloc  \advance\!yaxislength -\!areabloc
  \!plotheadingoffset=\!zpt
  \!dimenput {{\setbox0=\hbox{}\wd0=\!xaxislength\ht0=\!yaxislength\box0}}
     [bl] (\!arealloc,\!areabloc)}
%
\def\visibleaxes{%
  \def\!axisvisibility{\!axisvisibletrue}}

%

\def\!fixkeyword#1{%
  \errhelp=\!keywordhelp
  \errmessage{Unrecognized keyword `#1': \the\!keywordtoks{NEW KEYWORD}'}}

\!keywordtoks={enter `i\fixkeyword}

\def\fixkeyword#1{%
  \!nextkeyword#1 }


\def\axis {%
  \def\!nextkeyword##1 {%
    \expandafter\ifx\csname !axis##1\endcsname \relax
      \def\!next{\!fixkeyword{##1}}%
    \else
      \def\!next{\csname !axis##1\endcsname}%
    \fi
    \!next}%
  \!offset=\!zpt
  \!axisvisibility
  \!placeaxislabelfalse
  \!nextkeyword}

\def\!axisbottom{%
  \!axisylevel=\!areabloc
  \def\!tickxsign{0}%
  \def\!tickysign{-}%
  \def\!axissetup{\!axisxsetup}%
  \def\!axislabeltbrl{t}%
  \!nextkeyword}

\def\!axistop{%
  \!axisylevel=\!areatloc
  \def\!tickxsign{0}%
  \def\!tickysign{+}%
  \def\!axissetup{\!axisxsetup}%
  \def\!axislabeltbrl{b}%
  \!nextkeyword}

\def\!axisleft{%
  \!axisxlevel=\!arealloc
  \def\!tickxsign{-}%
  \def\!tickysign{0}%
  \def\!axissetup{\!axisysetup}%
  \def\!axislabeltbrl{r}%
  \!nextkeyword}

\def\!axisright{%
  \!axisxlevel=\!arearloc
  \def\!tickxsign{+}%
  \def\!tickysign{0}%
  \def\!axissetup{\!axisysetup}%
  \def\!axislabeltbrl{l}%
  \!nextkeyword}

\def\!axisshiftedto#1=#2 {%
  \if 0\!tickxsign
    \!axisylevel=\!M{#2}\!yunit
    \advance\!axisylevel -\!yorigin
  \else
    \!axisxlevel=\!M{#2}\!xunit
    \advance\!axisxlevel -\!xorigin
  \fi
  \!nextkeyword}

\def\!axisvisible{%
  \!axisvisibletrue  
  \!nextkeyword}

\def\!axisinvisible{%
  \!axisvisiblefalse
  \!nextkeyword}

\def\!axislabel#1 {%
  \!axisLaBeL={#1}%
  \!placeaxislabeltrue
  \!nextkeyword}

\expandafter\def\csname !axis/\endcsname{%
  \!axissetup 
  \if!placeaxislabel
    \!placeaxislabel
  \fi
  \if +\!tickysign 
    \!dimenA=\!axisylevel
    \advance\!dimenA \!offset 
    \advance\!dimenA -\!areatloc 
    \ifdim \!dimenA>\!plotheadingoffset
      \!plotheadingoffset=\!dimenA 
    \fi
  \fi}

\def\grid #1 #2 {%
  \!countA=#1\advance\!countA 1
  \axis bottom invisible ticks length <\!zpt> andacross quantity {\!countA} /
  \!countA=#2\advance\!countA 1
  \axis left   invisible ticks length <\!zpt> andacross quantity {\!countA} / }

\def\plotheading#1 {%
  \advance\!plotheadingoffset \headingtoplotskip
  \!dimenput {#1} [B] <.5\!xaxislength,\!plotheadingoffset>
    (\!arealloc,\!areatloc)}

\def\!axisxsetup{%
  \!axisxlevel=\!arealloc
  \!axisstart=\!arealloc
  \!axisend=\!arearloc
  \!axisLength=\!xaxislength
  \!!origin=\!xorigin
  \!!unit=\!xunit
  \!xswitchtrue
  \if!axisvisible 
    \!makeaxis
  \fi}

\def\!axisysetup{%
  \!axisylevel=\!areabloc
  \!axisstart=\!areabloc
  \!axisend=\!areatloc
  \!axisLength=\!yaxislength
  \!!origin=\!yorigin
  \!!unit=\!yunit
  \!xswitchfalse
  \if!axisvisible
    \!makeaxis
  \fi}

\def\!makeaxis{%
  \setbox\!boxA=\hbox{
    \beginpicture
      \!setdimenmode
      \setcoordinatesystem point at {\!zpt} {\!zpt}   
      \putrule from {\!zpt} {\!zpt} to
        {\!tickysign\!tickysign\!axisLength} 
        {\!tickxsign\!tickxsign\!axisLength}
    \endpicturesave <\!Xsave,\!Ysave>}%
    \wd\!boxA=\!zpt
    \!placetick\!axisstart}

\def\!placeaxislabel{%
  \advance\!offset \valuestolabelleading
  \if!xswitch
    \!dimenput {\the\!axisLaBeL} [\!axislabeltbrl]
      <.5\!axisLength,\!tickysign\!offset> (\!axisxlevel,\!axisylevel)
    \advance\!offset \!dp  
    \advance\!offset \!ht  
  \else
    \!dimenput {\the\!axisLaBeL} [\!axislabeltbrl]
      <\!tickxsign\!offset,.5\!axisLength> (\!axisxlevel,\!axisylevel)
  \fi
  \!axisLaBeL={}}

%


\def\arrow <#1> [#2,#3]{%
  \!ifnextchar<{\!arrow{#1}{#2}{#3}}{\!arrow{#1}{#2}{#3}<\!zpt,\!zpt> }}

\def\!arrow#1#2#3<#4,#5> from #6 #7 to #8 #9 {%
%
  \!xloc=\!M{#8}\!xunit   
  \!yloc=\!M{#9}\!yunit
  \!dxpos=\!xloc  \!dimenA=\!M{#6}\!xunit  \advance \!dxpos -\!dimenA
  \!dypos=\!yloc  \!dimenA=\!M{#7}\!yunit  \advance \!dypos -\!dimenA
  \let\!MAH=\!M
  \!setdimenmode
  \!xshift=#4\relax  \!yshift=#5\relax
  \!reverserotateonly\!xshift\!yshift
  \advance\!xshift\!xloc  \advance\!yshift\!yloc
%
  \!xS=-\!dxpos  \advance\!xS\!xshift
  \!yS=-\!dypos  \advance\!yS\!yshift
  \!start (\!xS,\!yS)
  \!ljoin (\!xshift,\!yshift)
%
  \!Pythag\!dxpos\!dypos\!arclength
  \!divide\!dxpos\!arclength\!dxpos  
  \!dxpos=32\!dxpos  \!removept\!dxpos\!!cos
  \!divide\!dypos\!arclength\!dypos  
  \!dypos=32\!dypos  \!removept\!dypos\!!sin
%
  \!halfhead{#1}{#2}{#3}
  \!halfhead{#1}{-#2}{-#3}
  \let\!M=\!MAH
  \ignorespaces}
%
  \def\!halfhead#1#2#3{%
    \!dimenC=-#1%
    \divide \!dimenC 2 
    \!dimenD=#2\!dimenC
    \!rotate(\!dimenC,\!dimenD)by(\!!cos,\!!sin)to(\!xM,\!yM)
    \!dimenC=-#1
    \!dimenD=#3\!dimenC
    \!dimenD=.5\!dimenD
    \!rotate(\!dimenC,\!dimenD)by(\!!cos,\!!sin)to(\!xE,\!yE)
    \!start (\!xshift,\!yshift)
    \advance\!xM\!xshift  \advance\!yM\!yshift
    \advance\!xE\!xshift  \advance\!yE\!yshift
    \!qjoin (\!xM,\!yM) (\!xE,\!yE) 
    \ignorespaces}

\def\betweenarrows #1#2 from #3 #4 to #5 #6 {%
  \!xloc=\!M{#3}\!xunit  \!xxloc=\!M{#5}\!xunit%
  \!yloc=\!M{#4}\!yunit  \!yyloc=\!M{#6}\!yunit%
  \!dxpos=\!xxloc  \advance\!dxpos by -\!xloc
  \!dypos=\!yyloc  \advance\!dypos by -\!yloc
  \advance\!xloc .5\!dxpos
  \advance\!yloc .5\!dypos
  \let\!MBA=\!M
  \!setdimenmode
  \ifdim\!dypos=\!zpt
    \ifdim\!dxpos<\!zpt \!dxpos=-\!dxpos \fi
    \put {\!lrarrows{\!dxpos}{#1}}#2{} at {\!xloc} {\!yloc}
  \else
    \ifdim\!dxpos=\!zpt
      \ifdim\!dypos<\!zpt \!dypos=-\!zpt \fi
      \put {\!udarrows{\!dypos}{#1}}#2{} at {\!xloc} {\!yloc}
    \fi
  \fi
  \let\!M=\!MBA
  \ignorespaces}

\def\!lrarrows#1#2{
  {\setbox\!boxA=\hbox{$\mkern-2mu\mathord-\mkern-2mu$}%
   \setbox\!boxB=\hbox{$\leftarrow$}\!dimenE=\ht\!boxB
   \setbox\!boxB=\hbox{}\ht\!boxB=2\!dimenE
   \hbox to #1{$\mathord\leftarrow\mkern-6mu
     \cleaders\copy\!boxA\hfil
     \mkern-6mu\mathord-$%
     \kern.4em $\vcenter{\box\!boxB}$$\vcenter{\hbox{#2}}$\kern.4em
     $\mathord-\mkern-6mu
     \cleaders\copy\!boxA\hfil
     \mkern-6mu\mathord\rightarrow$}}}

\def\!udarrows#1#2{
  {\setbox\!boxB=\hbox{#2}%
   \setbox\!boxA=\hbox to \wd\!boxB{\hss$\vert$\hss}%
   \!dimenE=\ht\!boxA \advance\!dimenE \dp\!boxA \divide\!dimenE 2
   \vbox to #1{\offinterlineskip
      \vskip .05556\!dimenE
      \hbox to \wd\!boxB{\hss$\mkern.4mu\uparrow$\hss}\vskip-\!dimenE
      \cleaders\copy\!boxA\vfil
      \vskip-\!dimenE\copy\!boxA
      \vskip\!dimenE\copy\!boxB\vskip.4em
      \copy\!boxA\vskip-\!dimenE
      \cleaders\copy\!boxA\vfil
      \vskip-\!dimenE \hbox to \wd\!boxB{\hss$\mkern.4mu\downarrow$\hss}
      \vskip .05556\!dimenE}}}

%

\def\putbar#1breadth <#2> from #3 #4 to #5 #6 {%
  \!xloc=\!M{#3}\!xunit  \!xxloc=\!M{#5}\!xunit%
  \!yloc=\!M{#4}\!yunit  \!yyloc=\!M{#6}\!yunit%
  \!dypos=\!yyloc  \advance\!dypos by -\!yloc
  \!dimenI=#2  
  \ifdim \!dimenI=\!zpt 
    \putrule#1from {#3} {#4} to {#5} {#6} 
  \else 
    \let\!MBar=\!M
    \!setdimenmode 
    \divide\!dimenI 2
    \ifdim \!dypos=\!zpt             
      \advance \!yloc -\!dimenI 
      \advance \!yyloc \!dimenI
    \else
      \advance \!xloc -\!dimenI 
      \advance \!xxloc \!dimenI
    \fi
    \putrectangle#1corners at {\!xloc} {\!yloc} and {\!xxloc} {\!yyloc}
    \let\!M=\!MBar 
  \fi
  \ignorespaces}

\def\setbars#1breadth <#2> baseline at #3 = #4 {%
  \edef\!barshift{#1}%
  \edef\!barbreadth{#2}%
  \edef\!barorientation{#3}%
  \edef\!barbaseline{#4}%
  \def\!bardobaselabel{\!bardoendlabel}%
  \def\!bardoendlabel{\!barfinish}%
  \let\!drawcurve=\!barcurve
  \!setbars}
\def\!setbars{%
  \futurelet\!nextchar\!!setbars}
\def\!!setbars{%
  \if b\!nextchar
    \def\!!!setbars{\!setbarsbget}%
  \else 
    \if e\!nextchar
      \def\!!!setbars{\!setbarseget}%
    \else
      \def\!!!setbars{\relax}%
    \fi
  \fi
  \!!!setbars}
\def\!setbarsbget baselabels (#1) {%
  \def\!barbaselabelorientation{#1}%
  \def\!bardobaselabel{\!!bardobaselabel}%
  \!setbars}
\def\!setbarseget endlabels (#1) {%
  \edef\!barendlabelorientation{#1}%
  \def\!bardoendlabel{\!!bardoendlabel}%
  \!setbars}

\def\!barcurve #1 #2 {%
  \if y\!barorientation
    \def\!basexarg{#1}%
    \def\!baseyarg{\!barbaseline}%
  \else
    \def\!basexarg{\!barbaseline}%
    \def\!baseyarg{#2}%
  \fi
  \expandafter\putbar\!barshift breadth <\!barbreadth> from {\!basexarg}
    {\!baseyarg} to {#1} {#2}
  \def\!endxarg{#1}%
  \def\!endyarg{#2}%
  \!bardobaselabel}

\def\!!bardobaselabel "#1" {%
  \put {#1}\!barbaselabelorientation{} at {\!basexarg} {\!baseyarg}
  \!bardoendlabel}
 
\def\!!bardoendlabel "#1" {%
  \put {#1}\!barendlabelorientation{} at {\!endxarg} {\!endyarg}
  \!barfinish}

\def\!barfinish{%
  \!ifnextchar/{\!finish}{\!barcurve}}

%
%
%
\def\putrectangle{%
  \!ifnextchar<{\!putrectangle}{\!putrectangle<\!zpt,\!zpt> }}
\def\!putrectangle<#1,#2> corners at #3 #4 and #5 #6 {%
%
  \!xone=\!M{#3}\!xunit  \!xtwo=\!M{#5}\!xunit%
  \!yone=\!M{#4}\!yunit  \!ytwo=\!M{#6}\!yunit%
  \ifdim \!xtwo<\!xone
    \!dimenI=\!xone  \!xone=\!xtwo  \!xtwo=\!dimenI
  \fi
  \ifdim \!ytwo<\!yone
    \!dimenI=\!yone  \!yone=\!ytwo  \!ytwo=\!dimenI
  \fi
  \!dimenI=#1\relax  \advance\!xone\!dimenI  \advance\!xtwo\!dimenI
  \!dimenI=#2\relax  \advance\!yone\!dimenI  \advance\!ytwo\!dimenI
  \let\!MRect=\!M
  \!setdimenmode
%
  \!shaderectangle
%
  \!dimenI=.5\linethickness
  \advance \!xone  -\!dimenI
  \advance \!xtwo   \!dimenI
  \putrule from {\!xone} {\!yone} to {\!xtwo} {\!yone} 
  \putrule from {\!xone} {\!ytwo} to {\!xtwo} {\!ytwo} 
%
  \advance \!xone   \!dimenI
  \advance \!xtwo  -\!dimenI%
  \advance \!yone  -\!dimenI
  \advance \!ytwo   \!dimenI
  \putrule from {\!xone} {\!yone} to {\!xone} {\!ytwo} 
  \putrule from {\!xtwo} {\!yone} to {\!xtwo} {\!ytwo} 
  \let\!M=\!MRect
  \ignorespaces}
 

\def\shaderectanglesoff{%
  \def\!shaderectangle{}%
  \ignorespaces}

\shaderectanglesoff
 
\def\!!shaderectangle{%
  \!dimenA=\!xtwo  \advance \!dimenA -\!xone
  \!dimenB=\!ytwo  \advance \!dimenB -\!yone
  \ifdim \!dimenA<\!dimenB
    \!startvshade (\!xone,\!yone,\!ytwo)
    \!lshade      (\!xtwo,\!yone,\!ytwo)
  \else
    \!starthshade (\!yone,\!xone,\!xtwo)
    \!lshade      (\!ytwo,\!xone,\!xtwo)
  \fi
  \ignorespaces}
  
\def\frame{%
  \!ifnextchar<{\!frame}{\!frame<\!zpt> }}
\long\def\!frame<#1> #2{%
  \beginpicture
    \setcoordinatesystem units <1pt,1pt> point at 0 0 
    \put {#2} [Bl] at 0 0 
    \!dimenA=#1\relax
    \!dimenB=\!wd \advance \!dimenB \!dimenA
    \!dimenC=\!ht \advance \!dimenC \!dimenA
    \!dimenD=\!dp \advance \!dimenD \!dimenA
    \let\!MFr=\!M
    \!setdimenmode
    \putrectangle corners at {-\!dimenA} {-\!dimenD} and {\!dimenB} {\!dimenC}
    \!setcoordmode
    \let\!M=\!MFr
  \endpicture
  \ignorespaces}
 
\def\rectangle <#1> <#2> {%
  \setbox0=\hbox{}\wd0=#1\ht0=#2\frame {\box0}}

%

\def\plot{%
  \!ifnextchar"{\!plotfromfile}{\!drawcurve}}
\def\!plotfromfile"#1"{%
  \expandafter\!drawcurve \input #1 /}

\def\setquadratic{%
  \let\!drawcurve=\!qcurve
  \let\!!Shade=\!!qShade
  \let\!!!Shade=\!!!qShade}

\def\setlinear{%
  \let\!drawcurve=\!lcurve
  \let\!!Shade=\!!lShade
  \let\!!!Shade=\!!!lShade}

\def\sethistograms{%
  \let\!drawcurve=\!hcurve}

\def\!qcurve #1 #2 {%
  \!start (#1,#2)
  \!Qjoin}
\def\!Qjoin#1 #2 #3 #4 {%
  \!qjoin (#1,#2) (#3,#4)             
  \!ifnextchar/{\!finish}{\!Qjoin}}

\def\!lcurve #1 #2 {%
  \!start (#1,#2)
  \!Ljoin}
\def\!Ljoin#1 #2 {%
  \!ljoin (#1,#2)                    
  \!ifnextchar/{\!finish}{\!Ljoin}}

\def\!finish/{\ignorespaces}

\def\!hcurve #1 #2 {%
  \edef\!hxS{#1}%
  \edef\!hyS{#2}%
  \!hjoin}
\def\!hjoin#1 #2 {%
  \putrectangle corners at {\!hxS} {\!hyS} and {#1} {#2}
  \edef\!hxS{#1}%
  \!ifnextchar/{\!finish}{\!hjoin}}

\def\vshade #1 #2 #3 {%
  \!startvshade (#1,#2,#3)
  \!Shadewhat}

\def\hshade #1 #2 #3 {%
  \!starthshade (#1,#2,#3)
  \!Shadewhat}

\def\!Shadewhat{%
  \futurelet\!nextchar\!Shade}
\def\!Shade{%
  \if <\!nextchar
    \def\!nextShade{\!!Shade}%
  \else
    \if /\!nextchar
      \def\!nextShade{\!finish}%
    \else
      \def\!nextShade{\!!!Shade}%
    \fi
  \fi
  \!nextShade}
\def\!!lShade<#1> #2 #3 #4 {%
  \!lshade <#1> (#2,#3,#4)                 
  \!Shadewhat}
\def\!!!lShade#1 #2 #3 {%
  \!lshade (#1,#2,#3)
  \!Shadewhat} 
\def\!!qShade<#1> #2 #3 #4 #5 #6 #7 {%
  \!qshade <#1> (#2,#3,#4) (#5,#6,#7)      
  \!Shadewhat}
\def\!!!qShade#1 #2 #3 #4 #5 #6 {%
  \!qshade (#1,#2,#3) (#4,#5,#6)
  \!Shadewhat} 

\setlinear

\def\setdashpattern <#1>{%
  \def\!Flist{}\def\!Blist{}\def\!UDlist{}%
  \!countA=0
  \!ecfor\!item:=#1\do{%
    \!dimenA=\!item\relax
    \expandafter\!rightappend\the\!dimenA\withCS{\\}\to\!UDlist%
    \advance\!countA  1
    \ifodd\!countA
      \expandafter\!rightappend\the\!dimenA\withCS{\!Rule}\to\!Flist%
      \expandafter\!leftappend\the\!dimenA\withCS{\!Rule}\to\!Blist%
    \else 
      \expandafter\!rightappend\the\!dimenA\withCS{\!Skip}\to\!Flist%
      \expandafter\!leftappend\the\!dimenA\withCS{\!Skip}\to\!Blist%
    \fi}%
  \!leaderlength=\!zpt
  \def\!Rule##1{\advance\!leaderlength  ##1}%
  \def\!Skip##1{\advance\!leaderlength  ##1}%
  \!Flist%
  \ifdim\!leaderlength>\!zpt 
  \else
    \def\!Flist{\!Skip{24in}}\def\!Blist{\!Skip{24in}}\ignorespaces
    \def\!UDlist{\\{\!zpt}\\{24in}}\ignorespaces
    \!leaderlength=24in
  \fi
  \!dashingon}

\def\!dashingon{%
  \def\!advancedashing{\!!advancedashing}%
  \def\!drawlinearsegment{\!lineardashed}%
  \def\!puthline{\!putdashedhline}%
  \def\!putvline{\!putdashedvline}%
  \ignorespaces}%
\def\!dashingoff{%
  \def\!advancedashing{\relax}%
  \def\!drawlinearsegment{\!linearsolid}%
  \def\!puthline{\!putsolidhline}%
  \def\!putvline{\!putsolidvline}%
  \ignorespaces}

\def\setdots{%
  \!ifnextchar<{\!setdots}{\!setdots<5pt>}}
\def\!setdots<#1>{%
  \!dimenB=#1\advance\!dimenB -\plotsymbolspacing
  \ifdim\!dimenB<\!zpt
    \!dimenB=\!zpt
  \fi
\setdashpattern <\plotsymbolspacing,\!dimenB>}
 
\def\setdotsnear <#1> for <#2>{%
  \!dimenB=#2\relax  \advance\!dimenB -.05pt  
  \!dimenC=#1\relax  \!countA=\!dimenC 
  \!dimenD=\!dimenB  \advance\!dimenD .5\!dimenC  \!countB=\!dimenD
  \divide \!countB  \!countA
  \ifnum 1>\!countB 
    \!countB=1
  \fi
  \divide\!dimenB  \!countB
  \setdots <\!dimenB>}
 
\def\setdashes{%
  \!ifnextchar<{\!setdashes}{\!setdashes<5pt>}}
\def\!setdashes<#1>{\setdashpattern <#1,#1>}
 
\def\setdashesnear <#1> for <#2>{%
  \!dimenB=#2\relax  
  \!dimenC=#1\relax  \!countA=\!dimenC 
  \!dimenD=\!dimenB  \advance\!dimenD .5\!dimenC  \!countB=\!dimenD
  \divide \!countB  \!countA
  \ifodd \!countB 
  \else 
    \advance \!countB  1
  \fi
  \divide\!dimenB  \!countB
  \setdashes <\!dimenB>}
 
\def\setsolid{%
  \def\!Flist{\!Rule{24in}}\def\!Blist{\!Rule{24in}}%
  \def\!UDlist{\\{24in}\\{\!zpt}}%
  \!dashingoff}  
\setsolid


 
  
 
\def\!divide#1#2#3{%
  \!dimenB=#1
  \!dimenC=#2
  \!dimenD=\!dimenB
  \divide \!dimenD \!dimenC
  \!dimenA=\!dimenD
  \multiply\!dimenD \!dimenC
  \advance\!dimenB -\!dimenD
  \!dimenD=\!dimenC
    \ifdim\!dimenD<\!zpt \!dimenD=-\!dimenD 
  \fi
  \ifdim\!dimenD<64pt
    \!divstep[\!tfs]\!divstep[\!tfs]%
  \else 
    \!!divide
  \fi
  #3=\!dimenA\ignorespaces}

\def\!!divide{%
  \ifdim\!dimenD<256pt
    \!divstep[64]\!divstep[32]\!divstep[32]%
  \else 
    \!divstep[8]\!divstep[8]\!divstep[8]\!divstep[8]\!divstep[8]%
    \!dimenA=2\!dimenA
  \fi}

\def\!divstep[#1]{
  \!dimenB=#1\!dimenB
  \!dimenD=\!dimenB
    \divide \!dimenD by \!dimenC
  \!dimenA=#1\!dimenA
    \advance\!dimenA by \!dimenD%
  \multiply\!dimenD by \!dimenC
    \advance\!dimenB by -\!dimenD}
 
\def\Divide <#1> by <#2> forming <#3> {%
  \!divide{#1}{#2}{#3}}

 
 

 

\def\ellipticalarc axes ratio #1:#2 #3 degrees from #4 #5 center at #6 #7 {%
  \!angle=#3pt\relax
  \ifdim\!angle>\!zpt 
    \def\!sign{}
  \else 
    \def\!sign{-}\!angle=-\!angle
  \fi
  \!xxloc=\!M{#6}\!xunit
  \!yyloc=\!M{#7}\!yunit     
  \!xxS=\!M{#4}\!xunit
  \!yyS=\!M{#5}\!yunit
  \advance\!xxS -\!xxloc
  \advance\!yyS -\!yyloc
  \!divide\!xxS{#1pt}\!xxS 
  \!divide\!yyS{#2pt}\!yyS 
  \let\!MC=\!M
  \!setdimenmode
  \!xS=#1\!xxS  \advance\!xS\!xxloc
  \!yS=#2\!yyS  \advance\!yS\!yyloc
  \!start (\!xS,\!yS)%
  \!loop\ifdim\!angle>14.9999pt
    \!rotate(\!xxS,\!yyS)by(\!cos,\!sign\!sin)to(\!xxM,\!yyM) 
    \!rotate(\!xxM,\!yyM)by(\!cos,\!sign\!sin)to(\!xxE,\!yyE)
    \!xM=#1\!xxM  \advance\!xM\!xxloc  \!yM=#2\!yyM  \advance\!yM\!yyloc
    \!xE=#1\!xxE  \advance\!xE\!xxloc  \!yE=#2\!yyE  \advance\!yE\!yyloc
    \!qjoin (\!xM,\!yM) (\!xE,\!yE)
    \!xxS=\!xxE  \!yyS=\!yyE 
    \advance \!angle -15pt
  \repeat
  \ifdim\!angle>\!zpt
    \!angle=100.53096\!angle
    \divide \!angle 360 
    \!sinandcos\!angle\!!sin\!!cos
    \!rotate(\!xxS,\!yyS)by(\!!cos,\!sign\!!sin)to(\!xxM,\!yyM) 
    \!rotate(\!xxM,\!yyM)by(\!!cos,\!sign\!!sin)to(\!xxE,\!yyE)
    \!xM=#1\!xxM  \advance\!xM\!xxloc  \!yM=#2\!yyM  \advance\!yM\!yyloc
    \!xE=#1\!xxE  \advance\!xE\!xxloc  \!yE=#2\!yyE  \advance\!yE\!yyloc
    \!qjoin (\!xM,\!yM) (\!xE,\!yE)
  \fi
  \let\!M=\!MC
  \ignorespaces}

\def\!rotate(#1,#2)by(#3,#4)to(#5,#6){%
  \!dimenA=#3#1\advance \!dimenA -#4#2
  \!dimenB=#3#2\advance \!dimenB  #4#1
  \divide \!dimenA 32  \divide \!dimenB 32 
  #5=\!dimenA  #6=\!dimenB
  \ignorespaces}
\def\!sin{4.17684}
\def\!cos{31.72624}

\def\!sinandcos#1#2#3{%
 \!dimenD=#1
 \!dimenA=\!dimenD
 \!dimenB=32pt
 \!removept\!dimenD\!value
 \!dimenC=\!dimenD
 \!dimenC=\!value\!dimenC \divide\!dimenC by 64 
 \advance\!dimenB by -\!dimenC
 \!dimenC=\!value\!dimenC \divide\!dimenC by 96 
 \advance\!dimenA by -\!dimenC
 \!dimenC=\!value\!dimenC \divide\!dimenC by 128 
 \advance\!dimenB by \!dimenC%
 \!removept\!dimenA#2
 \!removept\!dimenB#3
 \ignorespaces}




\def\putrule#1from #2 #3 to #4 #5 {%
  \!xloc=\!M{#2}\!xunit  \!xxloc=\!M{#4}\!xunit%
  \!yloc=\!M{#3}\!yunit  \!yyloc=\!M{#5}\!yunit%
  \!dxpos=\!xxloc  \advance\!dxpos by -\!xloc
  \!dypos=\!yyloc  \advance\!dypos by -\!yloc
  \ifdim\!dypos=\!zpt
    \def\!!Line{\!puthline{#1}}\ignorespaces
  \else
    \ifdim\!dxpos=\!zpt
      \def\!!Line{\!putvline{#1}}\ignorespaces
    \else 
       \def\!!Line{}
    \fi
  \fi
  \let\!ML=\!M
  \!setdimenmode
  \!!Line%
  \let\!M=\!ML
  \ignorespaces}

\def\!putsolidhline#1{%
  \ifdim\!dxpos>\!zpt 
    \put{\!hline\!dxpos}#1[l] at {\!xloc} {\!yloc}
  \else 
    \put{\!hline{-\!dxpos}}#1[l] at {\!xxloc} {\!yyloc}
  \fi
  \ignorespaces}
 
\def\!putsolidvline#1{%
  \ifdim\!dypos>\!zpt 
    \put{\!vline\!dypos}#1[b] at {\!xloc} {\!yloc}
  \else 
    \put{\!vline{-\!dypos}}#1[b] at {\!xxloc} {\!yyloc}
  \fi
  \ignorespaces}
 
\def\!hline#1{\hbox to #1{\leaders \hrule height\linethickness\hfill}}
\def\!vline#1{\vbox to #1{\leaders \vrule width\linethickness\vfill}}

\def\!putdashedhline#1{%
  \ifdim\!dxpos>\!zpt 
    \!DLsetup\!Flist\!dxpos
    \put{\hbox to \!totalleaderlength{\!hleaders}\!hpartialpattern\!Rtrunc}
      #1[l] at {\!xloc} {\!yloc} 
  \else 
    \!DLsetup\!Blist{-\!dxpos}
    \put{\!hpartialpattern\!Ltrunc\hbox to \!totalleaderlength{\!hleaders}}
      #1[r] at {\!xloc} {\!yloc} 
  \fi
  \ignorespaces}
 
\def\!putdashedvline#1{%
  \!dypos=-\!dypos
  \ifdim\!dypos>\!zpt 
    \!DLsetup\!Flist\!dypos 
    \put{\vbox{\vbox to \!totalleaderlength{\!vleaders}
      \!vpartialpattern\!Rtrunc}}#1[t] at {\!xloc} {\!yloc} 
  \else 
    \!DLsetup\!Blist{-\!dypos}
    \put{\vbox{\!vpartialpattern\!Ltrunc
      \vbox to \!totalleaderlength{\!vleaders}}}#1[b] at {\!xloc} {\!yloc} 
  \fi
  \ignorespaces}

\def\!DLsetup#1#2{
  \let\!RSlist=#1
  \!countB=#2
  \!countA=\!leaderlength
  \divide\!countB by \!countA
  \!totalleaderlength=\!countB\!leaderlength
  \!Rresiduallength=#2%
  \advance \!Rresiduallength by -\!totalleaderlength
  \!Lresiduallength=\!leaderlength
  \advance \!Lresiduallength by -\!Rresiduallength
  \ignorespaces}
 
\def\!hleaders{%
  \def\!Rule##1{\vrule height\linethickness width##1}%
  \def\!Skip##1{\hskip##1}%
  \leaders\hbox{\!RSlist}\hfill}
 
\def\!hpartialpattern#1{%
  \!dimenA=\!zpt \!dimenB=\!zpt 
  \def\!Rule##1{#1{##1}\vrule height\linethickness width\!dimenD}%
  \def\!Skip##1{#1{##1}\hskip\!dimenD}%
  \!RSlist}
 
\def\!vleaders{%
  \def\!Rule##1{\hrule width\linethickness height##1}%
  \def\!Skip##1{\vskip##1}%
  \leaders\vbox{\!RSlist}\vfill}
 
\def\!vpartialpattern#1{%
  \!dimenA=\!zpt \!dimenB=\!zpt 
  \def\!Rule##1{#1{##1}\hrule width\linethickness height\!dimenD}%
  \def\!Skip##1{#1{##1}\vskip\!dimenD}%
  \!RSlist}
 
\def\!Rtrunc#1{\!trunc{#1}>\!Rresiduallength}
\def\!Ltrunc#1{\!trunc{#1}<\!Lresiduallength}
 
\def\!trunc#1#2#3{%
  \!dimenA=\!dimenB         
  \advance\!dimenB by #1%
  \!dimenD=\!dimenB  \ifdim\!dimenD#2#3\!dimenD=#3\fi
  \!dimenC=\!dimenA  \ifdim\!dimenC#2#3\!dimenC=#3\fi
  \advance \!dimenD by -\!dimenC}

\def\!start (#1,#2){%
  \!plotxorigin=\!xorigin  \advance \!plotxorigin by \!plotsymbolxshift
  \!plotyorigin=\!yorigin  \advance \!plotyorigin by \!plotsymbolyshift
  \!xS=\!M{#1}\!xunit \!yS=\!M{#2}\!yunit
  \!rotateaboutpivot\!xS\!yS
  \!copylist\!UDlist\to\!!UDlist
  \!getnextvalueof\!downlength\from\!!UDlist
  \!distacross=\!zpt
  \!intervalno=0 
  \global\totalarclength=\!zpt
  \ignorespaces}

\def\!ljoin (#1,#2){%
  \advance\!intervalno by 1
  \!xE=\!M{#1}\!xunit \!yE=\!M{#2}\!yunit
  \!rotateaboutpivot\!xE\!yE
  \!xdiff=\!xE \advance \!xdiff by -\!xS
  \!ydiff=\!yE \advance \!ydiff by -\!yS
  \!Pythag\!xdiff\!ydiff\!arclength
  \global\advance \totalarclength by \!arclength%
  \!drawlinearsegment
  \!xS=\!xE \!yS=\!yE
  \ignorespaces}

\def\!linearsolid{%
  \!npoints=\!arclength
  \!countA=\plotsymbolspacing
  \divide\!npoints by \!countA
  \ifnum \!npoints<1 
    \!npoints=1 
  \fi
  \divide\!xdiff by \!npoints
  \divide\!ydiff by \!npoints
  \!xpos=\!xS \!ypos=\!yS
  \loop\ifnum\!npoints>-1
    \!plotifinbounds
    \advance \!xpos by \!xdiff
    \advance \!ypos by \!ydiff
    \advance \!npoints by -1
  \repeat
  \ignorespaces}

\def\!lineardashed{%
  \ifdim\!distacross>\!arclength
    \advance \!distacross by -\!arclength  
  \else
    \loop\ifdim\!distacross<\!arclength
      \!divide\!distacross\!arclength\!dimenA
      \!removept\!dimenA\!t
      \!xpos=\!t\!xdiff \advance \!xpos by \!xS
      \!ypos=\!t\!ydiff \advance \!ypos by \!yS
      \!plotifinbounds
      \advance\!distacross by \plotsymbolspacing
      \!advancedashing
    \repeat  
    \advance \!distacross by -\!arclength
  \fi
  \ignorespaces}

\def\!!advancedashing{%
  \advance\!downlength by -\plotsymbolspacing
  \ifdim \!downlength>\!zpt
  \else
    \advance\!distacross by \!downlength
    \!getnextvalueof\!uplength\from\!!UDlist
    \advance\!distacross by \!uplength
    \!getnextvalueof\!downlength\from\!!UDlist
  \fi}

\def\inboundscheckoff{%
  \def\!plotifinbounds{\!plot(\!xpos,\!ypos)}%
  \def\!initinboundscheck{\relax}\ignorespaces}
 
\inboundscheckoff
 
\def\!!plotifinbounds{%
  \ifdim \!xpos<\!checkleft
  \else
    \ifdim \!xpos>\!checkright
    \else
      \ifdim \!ypos<\!checkbot
      \else
         \ifdim \!ypos>\!checktop
         \else
           \!plot(\!xpos,\!ypos)
         \fi 
      \fi
    \fi
  \fi}

\def\!!initinboundscheck{%
  \!checkleft=\!arealloc     \advance\!checkleft by \!xorigin
  \!checkright=\!arearloc    \advance\!checkright by \!xorigin
  \!checkbot=\!areabloc      \advance\!checkbot by \!yorigin
  \!checktop=\!areatloc      \advance\!checktop by \!yorigin}

%


\def\!logten#1#2{%
  \expandafter\!!logten#1\!nil
  \!removept\!dimenF#2%
  \ignorespaces}

\def\!!logten#1#2\!nil{%
  \if -#1%
    \!dimenF=\!zpt
    \def\!next{\ignorespaces}%
  \else
    \if +#1%
      \def\!next{\!!logten#2\!nil}%
    \else
      \if .#1%
        \def\!next{\!!logten0.#2\!nil}%
      \else
        \def\!next{\!!!logten#1#2..\!nil}%
      \fi
    \fi
  \fi
  \!next}

\def\!!!logten#1#2.#3.#4\!nil{%
  \!dimenF=1pt 
  \if 0#1%
    \!!logshift#3pt 
  \else 
    \!logshift#2/
    \!dimenE=#1.#2#3pt 
  \fi 
  \ifdim \!dimenE<\!rootten
    \multiply \!dimenE 10 
    \advance  \!dimenF -1pt
  \fi
  \!dimenG=\!dimenE
    \advance\!dimenG 10pt
  \advance\!dimenE -10pt 
  \multiply\!dimenE 10 
  \!divide\!dimenE\!dimenG\!dimenE
  \!removept\!dimenE\!t
  \!dimenG=\!t\!dimenE
  \!removept\!dimenG\!tt
  \!dimenH=\!tt\!tenAe
    \divide\!dimenH 100
  \advance\!dimenH \!tenAc
  \!dimenH=\!tt\!dimenH
    \divide\!dimenH 100   
  \advance\!dimenH \!tenAa
  \!dimenH=\!t\!dimenH
    \divide\!dimenH 100 
  \advance\!dimenF \!dimenH}

\def\!logshift#1{%
  \if #1/%
    \def\!next{\ignorespaces}%
  \else
    \advance\!dimenF 1pt 
    \def\!next{\!logshift}%
  \fi 
  \!next}
 
 \def\!!logshift#1{%
   \advance\!dimenF -1pt
   \if 0#1%
     \def\!next{\!!logshift}%
   \else
     \if p#1%
       \!dimenF=1pt
       \def\!next{\!dimenE=1p}%
     \else
       \def\!next{\!dimenE=#1.}%
     \fi
   \fi
   \!next}

\def\beginpicture{%
  \setbox\!picbox=\hbox\bgroup%
  \!xleft=\maxdimen  
  \!xright=-\maxdimen
  \!ybot=\maxdimen
  \!ytop=-\maxdimen}
 
\def\endpicture{%
  \ifdim\!xleft=\maxdimen
    \!xleft=\!zpt \!xright=\!zpt \!ybot=\!zpt \!ytop=\!zpt 
  \fi
  \global\!Xleft=\!xleft \global\!Xright=\!xright
  \global\!Ybot=\!ybot \global\!Ytop=\!ytop
  \egroup%
  \ht\!picbox=\!Ytop  \dp\!picbox=-\!Ybot
  \ifdim\!Ybot>\!zpt
  \else 
    \ifdim\!Ytop<\!zpt
      \!Ybot=\!Ytop
    \else
      \!Ybot=\!zpt
    \fi
  \fi
  \hbox{\kern-\!Xleft\lower\!Ybot\box\!picbox\kern\!Xright}}
 
\def\endpicturesave <#1,#2>{%
  \endpicture \global #1=\!Xleft \global #2=\!Ybot \ignorespaces}

\def\setcoordinatesystem{%
  \!ifnextchar{u}{\!getlengths }
    {\!getlengths units <\!xunit,\!yunit>}}
\def\!getlengths units <#1,#2>{%
  \!xunit=#1\relax
  \!yunit=#2\relax
  \!ifcoordmode 
    \let\!SCnext=\!SCccheckforRP
  \else
    \let\!SCnext=\!SCdcheckforRP
  \fi
  \!SCnext}
\def\!SCccheckforRP{%
  \!ifnextchar{p}{\!cgetreference }
    {\!cgetreference point at {\!xref} {\!yref} }}
\def\!cgetreference point at #1 #2 {%
  \edef\!xref{#1}\edef\!yref{#2}%
  \!xorigin=\!xref\!xunit  \!yorigin=\!yref\!yunit  
  \!initinboundscheck 
  \ignorespaces}
\def\!SCdcheckforRP{%
  \!ifnextchar{p}{\!dgetreference}%
    {\ignorespaces}}
\def\!dgetreference point at #1 #2 {%
  \!xorigin=#1\relax  \!yorigin=#2\relax
  \ignorespaces}

\long\def\put#1#2 at #3 #4 {%
  \!setputobject{#1}{#2}%
  \!xpos=\!M{#3}\!xunit  \!ypos=\!M{#4}\!yunit  
  \!rotateaboutpivot\!xpos\!ypos%
  \advance\!xpos -\!xorigin  \advance\!xpos -\!xshift
  \advance\!ypos -\!yorigin  \advance\!ypos -\!yshift
  \kern\!xpos\raise\!ypos\box\!putobject\kern-\!xpos%
  \!doaccounting\ignorespaces}
 
\long\def\multiput #1#2 at {%
  \!setputobject{#1}{#2}%
  \!ifnextchar"{\!putfromfile}{\!multiput}}
\def\!putfromfile"#1"{%
  \expandafter\!multiput \input #1 /}
\def\!multiput{%
  \futurelet\!nextchar\!!multiput}
\def\!!multiput{%
  \if *\!nextchar
    \def\!nextput{\!alsoby}%
  \else
    \if /\!nextchar
      \def\!nextput{\!finishmultiput}%
    \else
      \def\!nextput{\!alsoat}%
    \fi
  \fi
  \!nextput}
\def\!finishmultiput/{%
  \setbox\!putobject=\hbox{}%
  \ignorespaces}
 
\def\!alsoat#1 #2 {%
  \!xpos=\!M{#1}\!xunit  \!ypos=\!M{#2}\!yunit  
  \!rotateaboutpivot\!xpos\!ypos%
  \advance\!xpos -\!xorigin  \advance\!xpos -\!xshift
  \advance\!ypos -\!yorigin  \advance\!ypos -\!yshift
  \kern\!xpos\raise\!ypos\copy\!putobject\kern-\!xpos%
  \!doaccounting
  \!multiput}
 
\def\!alsoby*#1 #2 #3 {%
  \!dxpos=\!M{#2}\!xunit \!dypos=\!M{#3}\!yunit 
  \!rotateonly\!dxpos\!dypos
  \!ntemp=#1%
  \!!loop\ifnum\!ntemp>0
    \advance\!xpos by \!dxpos  \advance\!ypos by \!dypos
    \kern\!xpos\raise\!ypos\copy\!putobject\kern-\!xpos%
    \advance\!ntemp by -1
  \repeat
  \!doaccounting 
  \!multiput}
 
\def\accountingon{\def\!doaccounting{\!!doaccounting}\ignorespaces}

\accountingon
\def\!!doaccounting{%
  \!xtemp=\!xpos  
  \!ytemp=\!ypos
  \ifdim\!xtemp<\!xleft 
     \!xleft=\!xtemp 
  \fi
  \advance\!xtemp by  \!wd 
  \ifdim\!xright<\!xtemp 
    \!xright=\!xtemp
  \fi
  \advance\!ytemp by -\!dp
  \ifdim\!ytemp<\!ybot  
    \!ybot=\!ytemp
  \fi
  \advance\!ytemp by  \!dp
  \advance\!ytemp by  \!ht 
  \ifdim\!ytemp>\!ytop  
    \!ytop=\!ytemp  
  \fi}
 
\long\def\!setputobject#1#2{%
  \setbox\!putobject=\hbox{#1}%
  \!ht=\ht\!putobject  \!dp=\dp\!putobject  \!wd=\wd\!putobject
  \wd\!putobject=\!zpt
  \!xshift=.5\!wd   \!yshift=.5\!ht   \advance\!yshift by -.5\!dp
  \edef\!putorientation{#2}%
  \expandafter\!SPOreadA\!putorientation[]\!nil%
  \expandafter\!SPOreadB\!putorientation<\!zpt,\!zpt>\!nil\ignorespaces}
 
\def\!SPOreadA#1[#2]#3\!nil{\!etfor\!orientation:=#2\do\!SPOreviseshift}
 
\def\!SPOreadB#1<#2,#3>#4\!nil{\advance\!xshift by -#2\advance\!yshift by -#3}
 
\def\!SPOreviseshift{%
  \if l\!orientation 
    \!xshift=\!zpt
  \else 
    \if r\!orientation 
      \!xshift=\!wd
    \else 
      \if b\!orientation
        \!yshift=-\!dp
      \else 
        \if B\!orientation 
          \!yshift=\!zpt
        \else 
          \if t\!orientation 
            \!yshift=\!ht
          \fi 
        \fi
      \fi
    \fi
  \fi}

\long\def\!dimenput#1#2(#3,#4){%
  \!setputobject{#1}{#2}%
  \!xpos=#3\advance\!xpos by -\!xshift
  \!ypos=#4\advance\!ypos by -\!yshift
  \kern\!xpos\raise\!ypos\box\!putobject\kern-\!xpos%
  \!doaccounting\ignorespaces}

\def\!setdimenmode{%
  \let\!M=\!M!!\ignorespaces}
\def\!setcoordmode{%
  \let\!M=\!M!\ignorespaces}
\def\!ifcoordmode{%
  \ifx \!M \!M!}
\def\!ifdimenmode{%
  \ifx \!M \!M!!}
\def\!M!#1#2{#1#2} 
\def\!M!!#1#2{#1}
\!setcoordmode
\let\setdimensionmode=\!setdimenmode
\let\setcoordinatemode=\!setcoordmode




\def\!stack[#1]{%
  \let\!lglue=\hfill \let\!rglue=\hfill
  \expandafter\let\csname !#1glue\endcsname=\relax
  \!ifnextchar<{\!!stack}{\!!stack<\stackleading>}}
\def\!!stack<#1>#2{%
  \vbox{\def\!valueslist{}\!ecfor\!value:=#2\do{%
    \expandafter\!rightappend\!value\withCS{\\}\to\!valueslist}%
    \!lop\!valueslist\to\!value
    \let\\=\cr\lineskiplimit=\maxdimen\lineskip=#1%
    \baselineskip=-1000pt\halign{\!lglue##\!rglue\cr \!value\!valueslist\cr}}%
  \ignorespaces}


\def\!lines[#1]#2{%
  \let\!lglue=\hfill \let\!rglue=\hfill
  \expandafter\let\csname !#1glue\endcsname=\relax
  \vbox{\halign{\!lglue##\!rglue\cr #2\crcr}}%
  \ignorespaces}


\def\!Lines[#1]#2{%
  \let\!lglue=\hfill \let\!rglue=\hfill
  \expandafter\let\csname !#1glue\endcsname=\relax
  \vtop{\halign{\!lglue##\!rglue\cr #2\crcr}}%
  \ignorespaces}

 
 
 
\def\setplotsymbol(#1#2){%
  \!setputobject{#1}{#2}
  \setbox\!plotsymbol=\box\!putobject%
  \!plotsymbolxshift=\!xshift 
  \!plotsymbolyshift=\!yshift 
  \ignorespaces}
 
\setplotsymbol({\fiverm .})

 
\def\!!plot(#1,#2){%
  \!dimenA=-\!plotxorigin \advance \!dimenA by #1
  \!dimenB=-\!plotyorigin \advance \!dimenB by #2
  \kern\!dimenA\raise\!dimenB\copy\!plotsymbol\kern-\!dimenA%
  \ignorespaces}
 
\def\!!!plot(#1,#2){%
  \!dimenA=-\!plotxorigin \advance \!dimenA by #1
  \!dimenB=-\!plotyorigin \advance \!dimenB by #2
  \kern\!dimenA\raise\!dimenB\copy\!plotsymbol\kern-\!dimenA%
  \!countE=\!dimenA
  \!countF=\!dimenB
  \immediate\write\!replotfile{\the\!countE,\the\!countF.}%
  \ignorespaces}

\def\savelinesandcurves on "#1" {%
  \immediate\closeout\!replotfile
  \immediate\openout\!replotfile=#1%
  \let\!plot=\!!!plot}

\def\dontsavelinesandcurves {%
  \let\!plot=\!!plot}
\dontsavelinesandcurves

{\catcode`\%=11\xdef\!Commentsignal{
\def\writesavefile#1 {%
  \immediate\write\!replotfile{\!Commentsignal #1}%
  \ignorespaces}

\def\replot"#1" {%
  \expandafter\!replot\input #1 /}
\def\!replot#1,#2. {%
  \!dimenA=#1sp
  \kern\!dimenA\raise#2sp\copy\!plotsymbol\kern-\!dimenA
  \futurelet\!nextchar\!!replot}
\def\!!replot{%
  \if /\!nextchar 
    \def\!next{\!finish}%
  \else
    \def\!next{\!replot}%
  \fi
  \!next}


 
 
\def\!Pythag#1#2#3{%
  \!dimenE=#1\relax                                     
  \ifdim\!dimenE<\!zpt 
    \!dimenE=-\!dimenE 
  \fi
  \!dimenF=#2\relax
  \ifdim\!dimenF<\!zpt 
    \!dimenF=-\!dimenF 
  \fi
  \advance \!dimenF by \!dimenE
  \ifdim\!dimenF=\!zpt 
    \!dimenG=\!zpt
  \else 
    \!divide{8\!dimenE}\!dimenF\!dimenE
    \advance\!dimenE by -4pt
      \!dimenE=2\!dimenE
    \!removept\!dimenE\!!t
    \!dimenE=\!!t\!dimenE
    \advance\!dimenE by 64pt
    \divide \!dimenE by 2
    \!dimenH=7pt
    \!!Pythag\!!Pythag\!!Pythag
    \!removept\!dimenH\!!t
    \!dimenG=\!!t\!dimenF
    \divide\!dimenG by 8
  \fi
  #3=\!dimenG
  \ignorespaces}

\def\!!Pythag{
  \!divide\!dimenE\!dimenH\!dimenI
  \advance\!dimenH by \!dimenI
    \divide\!dimenH by 2}

\def\placehypotenuse for <#1> and <#2> in <#3> {%
  \!Pythag{#1}{#2}{#3}}

 
 
 
\def\!qjoin (#1,#2) (#3,#4){%
  \advance\!intervalno by 1
  \!ifcoordmode
    \edef\!xmidpt{#1}\edef\!ymidpt{#2}%
  \else
    \!dimenA=#1\relax \edef\!xmidpt{\the\!dimenA}%
    \!dimenA=#2\relax \edef\!ymidpt{\the\!dimenA}%
  \fi
  \!xM=\!M{#1}\!xunit  \!yM=\!M{#2}\!yunit   \!rotateaboutpivot\!xM\!yM
  \!xE=\!M{#3}\!xunit  \!yE=\!M{#4}\!yunit   \!rotateaboutpivot\!xE\!yE
%
  \!dimenA=\!xM  \advance \!dimenA by -\!xS
  \!dimenB=\!xE  \advance \!dimenB by -\!xM
  \!xB=3\!dimenA \advance \!xB by -\!dimenB
  \!xC=2\!dimenB \advance \!xC by -2\!dimenA
%
  \!dimenA=\!yM  \advance \!dimenA by -\!yS%
  \!dimenB=\!yE  \advance \!dimenB by -\!yM%
  \!yB=3\!dimenA \advance \!yB by -\!dimenB%
  \!yC=2\!dimenB \advance \!yC by -2\!dimenA%
%
  \!xprime=\!xB  \!yprime=\!yB
  \!dxprime=.5\!xC  \!dyprime=.5\!yC
  \!getf \!midarclength=\!dimenA
  \!getf \advance \!midarclength by 4\!dimenA
  \!getf \advance \!midarclength by \!dimenA
  \divide \!midarclength by 12
%
  \!arclength=\!dimenA
  \!getf \advance \!arclength by 4\!dimenA
  \!getf \advance \!arclength by \!dimenA
  \divide \!arclength by 12
  \advance \!arclength by \!midarclength
  \global\advance \totalarclength by \!arclength
%
%
  \ifdim\!distacross>\!arclength 
    \advance \!distacross by -\!arclength
  \else
    \!initinverseinterp
    \loop\ifdim\!distacross<\!arclength
      \!inverseinterp
      \!xpos=\!t\!xC \advance\!xpos by \!xB
        \!xpos=\!t\!xpos \advance \!xpos by \!xS
      \!ypos=\!t\!yC \advance\!ypos by \!yB
        \!ypos=\!t\!ypos \advance \!ypos by \!yS
      \!plotifinbounds
      \advance\!distacross \plotsymbolspacing
      \!advancedashing
    \repeat  
    \advance \!distacross by -\!arclength
  \fi
  \!xS=\!xE
  \!yS=\!yE
  \ignorespaces}

\def\!getf{\!Pythag\!xprime\!yprime\!dimenA%
  \advance\!xprime by \!dxprime
  \advance\!yprime by \!dyprime}

\def\!initinverseinterp{%
  \ifdim\!arclength>\!zpt
    \!divide{8\!midarclength}\!arclength\!dimenE
    \ifdim\!dimenE<\!wmin \!setinverselinear
    \else 
      \ifdim\!dimenE>\!wmax \!setinverselinear
      \else
        \def\!inverseinterp{\!inversequad}\ignorespaces
%
%
         \!removept\!dimenE\!Ew
         \!dimenF=-\!Ew\!dimenE
         \advance\!dimenF by 32pt
         \!dimenG=8pt 
         \advance\!dimenG by -\!dimenE
         \!dimenG=\!Ew\!dimenG
         \!divide\!dimenF\!dimenG\!beta
         \!gamma=1pt
         \advance \!gamma by -\!beta
      \fi
    \fi
  \fi
  \ignorespaces}

\def\!inversequad{%
  \!divide\!distacross\!arclength\!dimenG
  \!removept\!dimenG\!v
  \!dimenG=\!v\!gamma
  \advance\!dimenG by \!beta
  \!dimenG=\!v\!dimenG
  \!removept\!dimenG\!t}

\def\!setinverselinear{%
  \def\!inverseinterp{\!inverselinear}%
  \divide\!dimenE by 8 \!removept\!dimenE\!t
  \!countC=\!intervalno \multiply \!countC 2
  \!countB=\!countC     \advance \!countB -1
  \!countA=\!countB     \advance \!countA -1
  \wlog{\the\!countB th point (\!xmidpt,\!ymidpt) being plotted 
    doesn't lie in the}%
  \wlog{ middle third of the arc between the \the\!countA th 
    and \the\!countC th points:}%
  \wlog{ [arc length \the\!countA\space to \the\!countB]/[arc length 
    \the \!countA\space to \the\!countC]=\!t.}%
  \ignorespaces}
 
\def\!inverselinear{%
  \!divide\!distacross\!arclength\!dimenG
  \!removept\!dimenG\!t}

 

\def\startrotation{%
  \let\!rotateaboutpivot=\!!rotateaboutpivot
  \let\!rotateonly=\!!rotateonly
  \!ifnextchar{b}{\!getsincos }%
    {\!getsincos by {\!cosrotationangle} {\!sinrotationangle} }}
\def\!getsincos by #1 #2 {%
  \edef\!cosrotationangle{#1}%
  \edef\!sinrotationangle{#2}%
  \!ifcoordmode 
    \let\!ROnext=\!ccheckforpivot
  \else
    \let\!ROnext=\!dcheckforpivot
  \fi
  \!ROnext}
\def\!ccheckforpivot{%
  \!ifnextchar{a}{\!cgetpivot}%
    {\!cgetpivot about {\!xpivotcoord} {\!ypivotcoord} }}
\def\!cgetpivot about #1 #2 {%
  \edef\!xpivotcoord{#1}%
  \edef\!ypivotcoord{#2}%
  \!xpivot=#1\!xunit  \!ypivot=#2\!yunit
  \ignorespaces}
\def\!dcheckforpivot{%
  \!ifnextchar{a}{\!dgetpivot}{\ignorespaces}}
\def\!dgetpivot about #1 #2 {%
  \!xpivot=#1\relax  \!ypivot=#2\relax
  \ignorespaces}

\def\stoprotation{%
  \let\!rotateaboutpivot=\!!!rotateaboutpivot
  \let\!rotateonly=\!!!rotateonly
  \ignorespaces}
 
\def\!!rotateaboutpivot#1#2{%
  \!dimenA=#1\relax  \advance\!dimenA -\!xpivot
  \!dimenB=#2\relax  \advance\!dimenB -\!ypivot
  \!dimenC=\!cosrotationangle\!dimenA
    \advance \!dimenC -\!sinrotationangle\!dimenB
  \!dimenD=\!cosrotationangle\!dimenB
    \advance \!dimenD  \!sinrotationangle\!dimenA
  \advance\!dimenC \!xpivot  \advance\!dimenD \!ypivot
  #1=\!dimenC  #2=\!dimenD
  \ignorespaces}

\def\!!rotateonly#1#2{%
  \!dimenA=#1\relax  \!dimenB=#2\relax 
  \!dimenC=\!cosrotationangle\!dimenA
    \advance \!dimenC -\!rotsign\!sinrotationangle\!dimenB
  \!dimenD=\!cosrotationangle\!dimenB
    \advance \!dimenD  \!rotsign\!sinrotationangle\!dimenA
  #1=\!dimenC  #2=\!dimenD
  \ignorespaces}
\def\!rotsign{}
\def\!!!rotateaboutpivot#1#2{\relax}
\def\!!!rotateonly#1#2{\relax}
\stoprotation

\def\!reverserotateonly#1#2{%
  \def\!rotsign{-}%
  \!rotateonly{#1}{#2}%
  \def\!rotsign{}%
  \ignorespaces}

\def\!getspan span <#1>{%
  \!dshade=#1\relax
  \!ifcoordmode 
    \let\!GRnext=\!GRccheckforAP
  \else
    \let\!GRnext=\!GRdcheckforAP
  \fi
  \!GRnext}
\def\!GRccheckforAP{%
  \!ifnextchar{p}{\!cgetanchor }
    {\!cgetanchor point at {\!xshadesave} {\!yshadesave} }}
\def\!cgetanchor point at #1 #2 {%
  \edef\!xshadesave{#1}\edef\!yshadesave{#2}%
  \!xshade=\!xshadesave\!xunit  \!yshade=\!yshadesave\!yunit
  \ignorespaces}
\def\!GRdcheckforAP{%
  \!ifnextchar{p}{\!dgetanchor}%
    {\ignorespaces}}
\def\!dgetanchor point at #1 #2 {%
  \!xshade=#1\relax  \!yshade=#2\relax
  \ignorespaces}

\def\setshadesymbol{%
  \!ifnextchar<{\!setshadesymbol}{\!setshadesymbol<,,,> }}

\def\!setshadesymbol <#1,#2,#3,#4> (#5#6){%
  \!setputobject{#5}{#6}%
  \setbox\!shadesymbol=\box\!putobject%
  \!shadesymbolxshift=\!xshift \!shadesymbolyshift=\!yshift
%
  \!dimenA=\!xshift \advance\!dimenA \!smidge
  \!override\!dimenA{#1}\!lshrinkage%
  \!dimenA=\!wd \advance \!dimenA -\!xshift
    \advance\!dimenA \!smidge
    \!override\!dimenA{#2}\!rshrinkage
  \!dimenA=\!dp \advance \!dimenA \!yshift
    \advance\!dimenA \!smidge
    \!override\!dimenA{#3}\!bshrinkage
  \!dimenA=\!ht \advance \!dimenA -\!yshift
    \advance\!dimenA \!smidge
    \!override\!dimenA{#4}\!tshrinkage
  \ignorespaces}
\def\!smidge{-.2pt}%

\def\!override#1#2#3{%
  \edef\!!override{#2}%
  \ifx \!!override\empty
    #3=#1\relax
  \else
    \if z\!!override
      #3=\!zpt
    \else
      \ifx \!!override\!blankz
        #3=\!zpt
      \else
        #3=#2\relax
      \fi
    \fi
  \fi
  \ignorespaces}
\def\!blankz{ z}

\setshadesymbol ({\fiverm .})

\def\!startvshade#1(#2,#3,#4){%
  \let\!!xunit=\!xunit%
  \let\!!yunit=\!yunit%
  \let\!!xshade=\!xshade%
  \let\!!yshade=\!yshade%
  \def\!getshrinkages{\!vgetshrinkages}%
  \let\!setshadelocation=\!vsetshadelocation%
  \!xS=\!M{#2}\!!xunit
  \!ybS=\!M{#3}\!!yunit
  \!ytS=\!M{#4}\!!yunit
  \!shadexorigin=\!xorigin  \advance \!shadexorigin \!shadesymbolxshift
  \!shadeyorigin=\!yorigin  \advance \!shadeyorigin \!shadesymbolyshift
  \ignorespaces}
 
\def\!starthshade#1(#2,#3,#4){%
  \let\!!xunit=\!yunit%
  \let\!!yunit=\!xunit%
  \let\!!xshade=\!yshade%
  \let\!!yshade=\!xshade%
  \def\!getshrinkages{\!hgetshrinkages}%
  \let\!setshadelocation=\!hsetshadelocation%
  \!xS=\!M{#2}\!!xunit
  \!ybS=\!M{#3}\!!yunit
  \!ytS=\!M{#4}\!!yunit
  \!shadexorigin=\!xorigin  \advance \!shadexorigin \!shadesymbolxshift
  \!shadeyorigin=\!yorigin  \advance \!shadeyorigin \!shadesymbolyshift
  \ignorespaces}

\def\!lattice#1#2#3#4#5{%
  \!dimenA=#1
  \!dimenB=#2
  \!countB=\!dimenB
%
  \!dimenC=#3
  \advance\!dimenC -\!dimenA
  \!countA=\!dimenC
  \divide\!countA \!countB
  \ifdim\!dimenC>\!zpt
    \!dimenD=\!countA\!dimenB
    \ifdim\!dimenD<\!dimenC
      \advance\!countA 1 
    \fi
  \fi
  \!dimenC=\!countA\!dimenB
    \advance\!dimenC \!dimenA
  #4=\!countA
  #5=\!dimenC
  \ignorespaces}

\def\!qshade#1(#2,#3,#4)#5(#6,#7,#8){%
  \!xM=\!M{#2}\!!xunit
  \!ybM=\!M{#3}\!!yunit
  \!ytM=\!M{#4}\!!yunit
  \!xE=\!M{#6}\!!xunit
  \!ybE=\!M{#7}\!!yunit
  \!ytE=\!M{#8}\!!yunit
  \!getcoeffs\!xS\!ybS\!xM\!ybM\!xE\!ybE\!ybB\!ybC
  \!getcoeffs\!xS\!ytS\!xM\!ytM\!xE\!ytE\!ytB\!ytC
  \def\!getylimits{\!qgetylimits}%
  \!shade{#1}\ignorespaces}
 
\def\!lshade#1(#2,#3,#4){%
  \!xE=\!M{#2}\!!xunit
  \!ybE=\!M{#3}\!!yunit
  \!ytE=\!M{#4}\!!yunit
  \!dimenE=\!xE  \advance \!dimenE -\!xS
  \!dimenC=\!ytE \advance \!dimenC -\!ytS
  \!divide\!dimenC\!dimenE\!ytB
  \!dimenC=\!ybE \advance \!dimenC -\!ybS
  \!divide\!dimenC\!dimenE\!ybB
  \def\!getylimits{\!lgetylimits}%
  \!shade{#1}\ignorespaces}
 
\def\!getcoeffs#1#2#3#4#5#6#7#8{%
  \!dimenC=#4\advance \!dimenC -#2
  \!dimenE=#3\advance \!dimenE -#1
  \!divide\!dimenC\!dimenE\!dimenF
  \!dimenC=#6\advance \!dimenC -#4
  \!dimenH=#5\advance \!dimenH -#3
  \!divide\!dimenC\!dimenH\!dimenG
  \advance\!dimenG -\!dimenF
  \advance \!dimenH \!dimenE
  \!divide\!dimenG\!dimenH#8
  \!removept#8\!t
  #7=-\!t\!dimenE
  \advance #7\!dimenF
  \ignorespaces}

\def\!shade#1{%
  \!getshrinkages#1<,,,>\!nil
  \advance \!dimenE \!xS
  \!lattice\!!xshade\!dshade\!dimenE
    \!parity\!xpos
  \!dimenF=-\!dimenF
    \advance\!dimenF \!xE
  \!loop\!not{\ifdim\!xpos>\!dimenF}
    \!shadecolumn%
    \advance\!xpos \!dshade
    \advance\!parity 1
  \repeat
  \!xS=\!xE
  \!ybS=\!ybE
  \!ytS=\!ytE
  \ignorespaces}

\def\!vgetshrinkages#1<#2,#3,#4,#5>#6\!nil{%
  \!override\!lshrinkage{#2}\!dimenE
  \!override\!rshrinkage{#3}\!dimenF
  \!override\!bshrinkage{#4}\!dimenG
  \!override\!tshrinkage{#5}\!dimenH
  \ignorespaces}
\def\!hgetshrinkages#1<#2,#3,#4,#5>#6\!nil{%
  \!override\!lshrinkage{#2}\!dimenG
  \!override\!rshrinkage{#3}\!dimenH
  \!override\!bshrinkage{#4}\!dimenE
  \!override\!tshrinkage{#5}\!dimenF
  \ignorespaces}

\def\!shadecolumn{%
  \!dxpos=\!xpos
  \advance\!dxpos -\!xS
  \!removept\!dxpos\!dx
  \!getylimits
  \advance\!ytpos -\!dimenH
  \advance\!ybpos \!dimenG
  \!yloc=\!!yshade
  \ifodd\!parity 
     \advance\!yloc \!dshade
  \fi
  \!lattice\!yloc{2\!dshade}\!ybpos%
    \!countA\!ypos
  \!dimenA=-\!shadexorigin \advance \!dimenA \!xpos
  \loop\!not{\ifdim\!ypos>\!ytpos}
    \!setshadelocation
    \!rotateaboutpivot\!xloc\!yloc%
    \!dimenA=-\!shadexorigin \advance \!dimenA \!xloc
    \!dimenB=-\!shadeyorigin \advance \!dimenB \!yloc
    \kern\!dimenA \raise\!dimenB\copy\!shadesymbol \kern-\!dimenA
    \advance\!ypos 2\!dshade
  \repeat
  \ignorespaces}
 
\def\!qgetylimits{%
  \!dimenA=\!dx\!ytC              
  \advance\!dimenA \!ytB
  \!ytpos=\!dx\!dimenA
  \advance\!ytpos \!ytS
  \!dimenA=\!dx\!ybC              
  \advance\!dimenA \!ybB
  \!ybpos=\!dx\!dimenA
  \advance\!ybpos \!ybS}
 
\def\!lgetylimits{%
  \!ytpos=\!dx\!ytB
  \advance\!ytpos \!ytS
  \!ybpos=\!dx\!ybB
  \advance\!ybpos \!ybS}
 
\def\!vsetshadelocation{
  \!xloc=\!xpos
  \!yloc=\!ypos}
\def\!hsetshadelocation{
  \!xloc=\!ypos
  \!yloc=\!xpos}





\def\!axisticks {%
  \def\!nextkeyword##1 {%
    \expandafter\ifx\csname !ticks##1\endcsname \relax
      \def\!next{\!fixkeyword{##1}}%
    \else
      \def\!next{\csname !ticks##1\endcsname}%
    \fi
    \!next}%
  \!axissetup
    \def\!axissetup{\relax}%
  \edef\!ticksinoutsign{\!ticksinoutSign}%
  \!ticklength=\longticklength
  \!tickwidth=\linethickness
  \!gridlinestatus
  \!setticktransform
  \!maketick
  \!tickcase=0
  \def\!LTlist{}%
  \!nextkeyword}

\def\ticksout{%
  \def\!ticksinoutSign{+}}

\ticksout

\def\nogridlines{%
  \def\!gridlinestatus{\!gridlinestoofalse}}
\nogridlines

\def\loggedticks{%
  \def\!setticktransform{\let\!ticktransform=\!logten}}
\def\unloggedticks{%
  \def\!setticktransform{\let\!ticktransform=\!donothing}}
\def\!donothing#1#2{\def#2{#1}}
\unloggedticks

\expandafter\def\csname !ticks/\endcsname{%
  \!not {\ifx \!LTlist\empty}
    \!placetickvalues
  \fi
  \def\!tickvalueslist{}%
  \def\!LTlist{}%
  \expandafter\csname !axis/\endcsname}

\def\!maketick{%
  \setbox\!boxA=\hbox{%
    \beginpicture
      \!setdimenmode
      \setcoordinatesystem point at {\!zpt} {\!zpt}   
      \linethickness=\!tickwidth
      \ifdim\!ticklength>\!zpt
        \putrule from {\!zpt} {\!zpt} to
          {\!ticksinoutsign\!tickxsign\!ticklength}
          {\!ticksinoutsign\!tickysign\!ticklength}
      \fi
      \if!gridlinestoo
        \putrule from {\!zpt} {\!zpt} to
          {-\!tickxsign\!xaxislength} {-\!tickysign\!yaxislength}
      \fi
    \endpicturesave <\!Xsave,\!Ysave>}%
    \wd\!boxA=\!zpt}
  
\def\!ticksin{%
  \def\!ticksinoutsign{-}%
  \!maketick
  \!nextkeyword}

\def\!ticksout{%
  \def\!ticksinoutsign{+}%
  \!maketick
  \!nextkeyword}

\def\!tickslength<#1> {%
  \!ticklength=#1\relax
  \!maketick
  \!nextkeyword}

\def\!tickslong{%
  \!tickslength<\longticklength> }

\def\!ticksshort{%
  \!tickslength<\shortticklength> }

\def\!tickswidth<#1> {%
  \!tickwidth=#1\relax
  \!maketick
  \!nextkeyword}

\def\!ticksandacross{%
  \!gridlinestootrue
  \!maketick
  \!nextkeyword}

\def\!ticksbutnotacross{%
  \!gridlinestoofalse
  \!maketick
  \!nextkeyword}

\def\!tickslogged{%
  \let\!ticktransform=\!logten
  \!nextkeyword}

\def\!ticksunlogged{%
  \let\!ticktransform=\!donothing
  \!nextkeyword}

\def\!ticksunlabeled{%
  \!tickcase=0
  \!nextkeyword}

\def\!ticksnumbered{%
  \!tickcase=1
  \!nextkeyword}

\def\!tickswithvalues#1/ {%
  \edef\!tickvalueslist{#1! /}%
  \!tickcase=2
  \!nextkeyword}

\def\!ticksquantity#1 {%
  \ifnum #1>1
    \!updatetickoffset
    \!countA=#1\relax
    \advance \!countA -1
    \!ticklocationincr=\!axisLength
      \divide \!ticklocationincr \!countA
    \!ticklocation=\!axisstart
    \loop \!not{\ifdim \!ticklocation>\!axisend}
      \!placetick\!ticklocation
      \ifcase\!tickcase
          \relax 
        \or
          \relax 
        \or
          \expandafter\!gettickvaluefrom\!tickvalueslist
          \edef\!tickfield{{\the\!ticklocation}{\!value}}%
          \expandafter\!listaddon\expandafter{\!tickfield}\!LTlist%
      \fi
      \advance \!ticklocation \!ticklocationincr
    \repeat
  \fi
  \!nextkeyword}

\def\!ticksat#1 {%
  \!updatetickoffset
  \edef\!Loc{#1}%
  \if /\!Loc
    \def\next{\!nextkeyword}%
  \else
    \!ticksincommon
    \def\next{\!ticksat}%
  \fi
  \next}    
      
\def\!ticksfrom#1 to #2 by #3 {%
  \!updatetickoffset
  \edef\!arg{#3}%
  \expandafter\!separate\!arg\!nil
  \!scalefactor=1
  \expandafter\!countfigures\!arg/
  \edef\!arg{#1}%
  \!scaleup\!arg by\!scalefactor to\!countE
  \edef\!arg{#2}%
  \!scaleup\!arg by\!scalefactor to\!countF
  \edef\!arg{#3}%
  \!scaleup\!arg by\!scalefactor to\!countG
  \loop \!not{\ifnum\!countE>\!countF}
    \ifnum\!scalefactor=1
      \edef\!Loc{\the\!countE}%
    \else
      \!scaledown\!countE by\!scalefactor to\!Loc
    \fi
    \!ticksincommon
    \advance \!countE \!countG
  \repeat
  \!nextkeyword}

\def\!updatetickoffset{%
  \!dimenA=\!ticksinoutsign\!ticklength
  \ifdim \!dimenA>\!offset
    \!offset=\!dimenA
  \fi}

\def\!placetick#1{%
  \if!xswitch
    \!xpos=#1\relax
    \!ypos=\!axisylevel
  \else
    \!xpos=\!axisxlevel
    \!ypos=#1\relax
  \fi
  \advance\!xpos \!Xsave
  \advance\!ypos \!Ysave
  \kern\!xpos\raise\!ypos\copy\!boxA\kern-\!xpos
  \ignorespaces}

\def\!gettickvaluefrom#1 #2 /{%
  \edef\!value{#1}%
  \edef\!tickvalueslist{#2 /}%
  \ifx \!tickvalueslist\!endtickvaluelist
    \!tickcase=0
  \fi}
\def\!endtickvaluelist{! /}

\def\!ticksincommon{%
  \!ticktransform\!Loc\!t
  \!ticklocation=\!t\!!unit
  \advance\!ticklocation -\!!origin
  \!placetick\!ticklocation
  \ifcase\!tickcase
    \relax 
  \or 
    \ifdim\!ticklocation<-\!!origin
      \edef\!Loc{$\!Loc$}%
    \fi
    \edef\!tickfield{{\the\!ticklocation}{\!Loc}}%
    \expandafter\!listaddon\expandafter{\!tickfield}\!LTlist%
  \or 
    \expandafter\!gettickvaluefrom\!tickvalueslist
    \edef\!tickfield{{\the\!ticklocation}{\!value}}%
    \expandafter\!listaddon\expandafter{\!tickfield}\!LTlist%
  \fi}

\def\!separate#1\!nil{%
  \!ifnextchar{-}{\!!separate}{\!!!separate}#1\!nil}
\def\!!separate-#1\!nil{%
  \def\!sign{-}%
  \!!!!separate#1..\!nil}
\def\!!!separate#1\!nil{%
  \def\!sign{+}%
  \!!!!separate#1..\!nil}
\def\!!!!separate#1.#2.#3\!nil{%
  \def\!arg{#1}%
  \ifx\!arg\!empty
    \!countA=0
  \else
    \!countA=\!arg
  \fi
  \def\!arg{#2}%
  \ifx\!arg\!empty
    \!countB=0
  \else
    \!countB=\!arg
  \fi}
 
\def\!countfigures#1{%
  \if #1/%
    \def\!next{\ignorespaces}%
  \else
    \multiply\!scalefactor 10
    \def\!next{\!countfigures}%
  \fi
  \!next}

\def\!scaleup#1by#2to#3{%
  \expandafter\!separate#1\!nil
  \multiply\!countA #2\relax
  \advance\!countA \!countB
  \if -\!sign
    \!countA=-\!countA
  \fi
  #3=\!countA
  \ignorespaces}

\def\!scaledown#1by#2to#3{%
  \!countA=#1\relax
  \ifnum \!countA<0 
    \def\!sign{-}
    \!countA=-\!countA
  \else
    \def\!sign{}%
  \fi
  \!countB=\!countA
  \divide\!countB #2\relax
  \!countC=\!countB
    \multiply\!countC #2\relax
  \advance \!countA -\!countC
  \edef#3{\!sign\the\!countB.}
  \!countC=\!countA 
  \ifnum\!countC=0 
    \!countC=1
  \fi
  \multiply\!countC 10
  \!loop \ifnum #2>\!countC
    \edef#3{#3\!zero}%
    \multiply\!countC 10
  \repeat
  \edef#3{#3\the\!countA}
  \ignorespaces}

\def\!placetickvalues{%
  \advance\!offset \tickstovaluesleading
  \if!xswitch
    \setbox\!boxA=\hbox{%
      \def\\##1##2{%
        \!dimenput {##2} [B] (##1,\!axisylevel)}%
      \beginpicture 
        \!LTlist
      \endpicturesave <\!Xsave,\!Ysave>}%
    \!dimenA=\!axisylevel
      \advance\!dimenA -\!Ysave
      \advance\!dimenA \!tickysign\!offset
      \if -\!tickysign
        \advance\!dimenA -\ht\!boxA
      \else
        \advance\!dimenA  \dp\!boxA
      \fi
    \advance\!offset \ht\!boxA 
      \advance\!offset \dp\!boxA
    \!dimenput {\box\!boxA} [Bl] <\!Xsave,\!Ysave> (\!zpt,\!dimenA)
  \else
    \setbox\!boxA=\hbox{%
      \def\\##1##2{%
        \!dimenput {##2} [r] (\!axisxlevel,##1)}%
      \beginpicture 
        \!LTlist
      \endpicturesave <\!Xsave,\!Ysave>}%
    \!dimenA=\!axisxlevel
      \advance\!dimenA -\!Xsave
      \advance\!dimenA \!tickxsign\!offset
      \if -\!tickxsign
        \advance\!dimenA -\wd\!boxA
      \fi
    \advance\!offset \wd\!boxA
    \!dimenput {\box\!boxA} [Bl] <\!Xsave,\!Ysave> (\!dimenA,\!zpt)
  \fi}

\normalgraphs
\catcode`!=12 

\hfuzz=2pt
\font\svtnrm=cmr17

\def\Got#1{\hbox{\aa#1}}

\def\tU{\tilde{{U}}}

\def\tX{\tilde{{Y}}}

\def\gsp1{{\Got s}{\Got p}(1)}

\def\tX{\tilde{X}}

\def\secl{1}
\def\seg{2}
\def\ser{3}
\def\sewz{4}
\def\seex{5}
\def\seoe{6}
\def\tU{\tilde{U}}
\def\tX{\tilde{X}}
\def\tx{\tilde{x}}
\def\tP{\tilde{P}}

\def\Se{Sasakian-Einstein }

\centerline{\svtnrm On Sasakian-Einstein Geometry}
\bigskip
\centerline{\sc Charles P. Boyer~~ Krzysztof Galicki~~
\footnote{}{\ninerm During the preparation of this work the authors
were supported by NSF grant DMS-9423752.}}
\bigskip
\bigskip
\baselineskip = 10 truept
\centerline{\bf An Introduction}
\bigskip
In 1960 Sasaki [Sas] introduced a type of metric-contact structure which can
be thought of as the odd-dimensional version of K\"ahler geometry.  This
geometry became known as Sasakian geometry, and although it has been studied
fairly extensively ever since it has never gained quite the reputation of its
older sister -- K\"ahlerian geometry. Nevertheless, it has appeared in an
increasing number of different contexts from quaternionic geometry to
mathematical physics.

In this article we shall focus our attention on a special class of Sasakian
manifolds (and orbifolds) which have the property that the metric $g$ is
Einstein, that is ${\rm Ric}_g=\lambda g$. Such spaces are called
Sasakian-Einstein. Perhaps one reason that the study of \Se manifolds is so
attractive is that they have generic holonomy [BG2]. However, the reason that
they are so tractable is that they are closely related to a reduced holonomy.
In fact, the most geometric definition of a Sasakian structure is: a smooth
manifold $(\cals, g)$ is Sasakian if the metric cone $(C(\cals),
\bar{g}=dr^2+r^2g)$ is K\"ahler, {\it i.e.,} its holonomy group reduces to a
subgroup of $U({m+1\over2})$, where $m={\rm dim}(\cals)=2n+1$. Moreover,
$(\cals,g)$ is Sasakian-Einstein if and only if its metric cone
$(C(\cals),\bar{g})$ is K\"ahler and Ricci-flat, {\it i.e.,} its restricted
holonomy group reduces to a subgroup of $SU({n+1})$ (Calabi-Yau geometry).
In particular, the Sasakian-Einstein geometry properly contains
the so-called 3-Sasakian spaces for which the metric cone is not just
Calabi-Yau but hyperk\"ahler. The 3-Sasakian spaces are intimately related
to quaternionic K\"ahler geometry and from this point of
view they have  been investigated in a series of recent articles [BGM1-6, BGMR,
BG1-2, GS]. The simplest example of a compact simply connected
Sasakian-Einstein manifold is furnished by the unit
round sphere $S^{2n+1}$ whose associated metric cone $\bbc^{n+1}\setminus
\{0\}$ is flat. One of the first examples for which the cone
$C(\cals)$ is not flat was constructed by Tanno [Tan].
Applying certain embedding techniques Tanno showed
that $\cals=S^2\times S^3$ supports a homogeneous
Sasakian-Einstein structure.

Sasakian geometry has associated with it a characteristic vector field [Bl].
This vector field is non-vanishing and thus generates a 1-dimensional
foliation -- the characteristic foliation on the Sasakian manifold $\cals.$ If
we make an additional assumption that the leaves of this foliation are compact
then the space of leaves will be a K\"ahler orbifold. This is both at once a
generalization and a specialization of the well-known Boothby-Wang fibration
[Bl, Hat].  It is a specialization since we are dealing with Sasakian and not
the more general contact geometry. It is a generalization in that we make no
regularity assumption on the foliation, but only assume that the leaves are
compact. We refer to this as {\it quasi-regularity}, and it is this condition
that places us within the category of orbifolds.

In the context of Einstein geometry we are dealing with the orbifold version of
a result of Kobayashi [Kob2] which says that the total space of a principal
circle bundle over a K\"ahler-Einstein manifold of positive scalar curvature
admits an Einstein metric of positive scalar curvature. Recently a
generalization of the Kobayashi's
construction was brought to fruition in a paper by Wang and Ziller [WZ], where
the authors construct Einstein metrics on circle bundles
(and higher dimensional torus bundles) over products of positive
scalar curvature K\"ahler-Einstein manifolds.
We contend that the correct setting for  Wang and Ziller's result in
the case of circle bundles is a commutative  topological
monoid structure $(\calr,\star)$ on the set
$\calr$ of all {\sl regular} Sasakian-Einstein manifolds.
(We refer to ``$\star$" as the ``join".)
In general the Einstein metrics obtained via the Wang-Ziller bundle
construction are not Sasakian-Einstein. However, over each
base $M=M_1\times\cdots\times M_n$ one finds a unique simply connected
circle bundle $\cals$ which is Sasakian-Einstein and as an element of
$\calr$ it is a product of $n$
factors $\cals_1\star\cdots\star\cals_n$.
Furthermore, such $\cals$ has an $n$-dimensional lattice
$L(\cals_1,\ldots,\cals_n)$ of compact Einstein
manifolds canonically associated to it. The points on this lattice
give all the other circle bundles in the Wang-Ziller construction.

It is not hard
to realize that the regularity assumption on $S_i$'s
is much too restrictive even if
one is solely interested in the smooth category of compact Einstein manifolds.
Indeed, our contribution is not merely a novel look at Wang and Ziller's bundle
construction, nor is it the recognition of the central role played by
Sasakian-Einstein geometry; it is, however, mainly the non-trivial
generalization of the monoid $(\calr,\star)$ of regular \Se manifolds to a
monoid structure on the set $\cals\cale$ of compact \Se orbifolds. A crucial
point is that even though the subset of quasi-regular \Se manifolds in
$\cals\cale$ is not closed under the join operation, one can analyze the
conditions necessary for the join to be smooth. A key ingredient in this
construction is
Haefliger's description [Hae] of the classifying space of an orbifold which
allows him to define ``orbifold cohomology''. In turn this allows us to
generalize the notion of the index of a smooth Fano variety to Fano orbifolds;
hence, we can construct a V-bundle with a \Se structure whose total space is
simply connected, and the index is associated to the Sasakian structure. Then
given a pair of quasi-regular \Se manifolds, or more generally orbifolds, we
can define their relative indices by dividing each index by the gcd of their
indices. It also makes sense to define the order of a quasi-regular Sasakian
manifold (orbifold) to be the order as an orbifold of the space of leaves of
the characteristic foliation.  Our main result is:

\noindent{\sc Theorem A}: \tensl Let $\cals_1,\cals_2$ be two simply connected
quasi-regular \Se orbifolds of dimensions $2n_1+1$ and $2n_2+1,$ respectively.
Let $\cals_i$ have orders $m_i$ and relative indices $l_i.$ Then there exists a
multiplication called the {\it join} and denoted by $\star$ such that
$\cals_1\star \cals_2$ is a simply connected quasi-regular \Se orbifold of
dimension $2(n_1+n_2)+1.$ If both $\cals_1$ and $\cals_2$ are smooth \Se
manifolds then $\cals_1\star \cals_2$ is a smooth manifold if and only if
$\gcd(m_1l_2,m_2l_1)=1.$ Furthermore, to each such pair there is a
two-dimensional lattice of Einstein orbifolds each having the same rational
cohomology as $\cals_1\star \cals_2.$ \tenrm

Similarly, one can determine
necessary conditions for a lattice point on $L(\cals_1,\cals_2)$ to be a smooth
Einstein manifold. As a consequence we obtain many new Einstein and
Sasakian-Einstein manifolds.
Using a simple and elegant spectral sequence argument employed by Wang and
Ziller we are able to compute the cohomology rings of many examples of the
joins of \Se manifolds. This will allow us to generalize some of our previous
results [BGMR, BGM2]. For example, we have

\noindent{\sc Corollary} B: \tensl In every odd dimension greater than 5, there
are infinitely many distinct homotopy types of simply connected compact \Se
manifolds having the same rational cohomology groups.  In particular, in each
such dimension, there are infinitely many \Se manifolds with arbitrarily small
injectivity radii.  \tenrm

\noindent{\sc Corollary} C: \tensl In every odd dimension greater than $5$
there are simply connected compact \Se manifolds with any second Betti number.
In particular in each such dimension, there are infinitely many simply
connected \Se manifolds with the property that given any negative real number
$\grk$ there exist no metrics on $\cals$ whose sectional curvatures are all
greater than or equal to $\grk.$ \tenrm

\noindent{\sc Corollary} D: \tensl In each dimension of the form $4n+3$ with
$n\geq 1$ there are compact \Se manifolds that do not admit a 3-Sasakian
structure, and for $n> 2$ there are such manifolds having arbitrary second
Betti number. \tenrm

\noindent{\sc Corollary} E: \tensl In every odd dimension greater than $3$
there are compact simply connected manifolds that admit continuous families of
\Se structures. \tenrm

The second statement in Corollary B follows by a well-known result of Anderson
[An] while the second statement of Corollary C follows from a famous result of
Gromov [Gro].

Concluding this introduction it seems worth mentioning that very recently
Sasakian-Einstein geometry has emerged in the context of dualities of certain
supersymmetric conformal field theories. The whole story begins with an
important conjecture of of Maldacena [Mal] who noticed that the large $N$ limit
of certain conformal field theories in $d$ dimensions can be described in terms
of supergravity (and string theory) on a product of $(d+1)$-dimensional
anti-de-Sitter $AdS_{d+1}$ space with a compact manifold $M$. The idea was
later examined by Witten who proposed a precise correspondence between
conformal field theory observables and those of supergravity [Wit].  It turns
out, and this observation has recently been made by Figueroa [Fi], that $M$
necessarily has real Killing spinors and the number of them determines the
number of supersymmetries preserved. Depending on the dimension and the amount
of supersymmetry, the following geometries are possible: spherical in any
dimension, Sasakian-Einstein in dimension $2k+1$, 3-Sasakian in dimension
$4k+3$, 7-manifolds with weak $G_2$-holonomy, and 6-dimensional nearly K\"ahler
manifolds [AFHB].  The case when ${\rm dim}(M)=5,7$ seems to be of particular
interest.  Several other articles have begun an in-depth investigation of these
ideas.  In particular the paper of Klebanov and Witten [KW] examines this
duality in the case of $M=S^2\times S^3=S^3\star S^3$, Oh and Tatar do the same
for $M=S^3\star S^3\star S^3$ [OhTa].
An article by Morrison and Plesser
formulates an extension of Maldacena's Conjecture in the general case of
Non-Spherical Horizons [MP].  These new developments are certainly a part of
our motivation to embark on a more systematic study of Sasakian-Einstein
geometry and this article takes the first few steps in this direction.

\noindent{\sc Acknowledgments}: We would like to thank Jos\'e
Figueroa-O'Farrill for
telling us about the recent results concerning string branes
at conical singularities, as well as Ben Mann and Jim Milgram for helpful
discussions.  The second named author would like to thank Max-Planck-Institute
f\"ur Mathematik in Bonn for support and hospitality. This project started
during his stay in Bonn in 1997.

\bigskip
\centerline {\bf \S 1. The Classifying Space of an Orbifold and V-bundles}
\medskip

An {\it orbifold} $X$ is an important generalization of a manifold in which
the locally Euclidean charts  $(U_i,\phi_i)$ are replaced by
the {\it local uniformizing systems} $\{\tU_i,\grG_i,\varphi_i\}$, where
$\tU_i$ is an open subset of $\bbr^n$ (or $\bbc^n$) containing the origin,
$\grG_i$ is a finite group of diffeomorphisms ( or biholomorphisms) which
can be taken to be elements of $O(n)$ (or $U(n)$),
$\varphi_i:\tU_i\ra{1.3} U_i$ is a continuous map onto an subset $U_i\subset
X$ such that $\varphi_i\circ \grg=\varphi_i$ for all $\grg\in\grG_i$
and the induced natural map of $\tU_i/\grG_i$ onto $U_i$ is a homeomorphism.
The finite group $\grG_i$ is called a {\it local uniformizing group}.
When it is defined the least common multiple of the orders of the
local uniformizing groups $\grG_i$ is called the {\it order} of the orbifold
$X,$ and is denoted by $\hbox{Ord}(X).$ We shall be particularly
interested in the case when $X=\calz$ a compact {\it complex orbifold}.
Notice that for a compact orbifold the order of the orbifold is always defined.
We refer to our expository paper [BG2]
as well as the original literature for more details.

As with manifolds an alternative definition of orbifold
can be given in terms of groupoids. Following Haefliger [Hae]:
Let $\calg_X$ denote the groupoid of germs of diffeomorphisms
generated by the germs of elements in $\grG_i$ and the germs of the
diffeomorphisms $g_{ji}$ described above. Let $\tX =\sqcup_i \tU_i$ denote the
disjoint union of the $\tU_i.$ Then $x,y\in \tX$ are equivalent if there is a
germ $\grg \in \calg_X$ such that $y=\grg(x).$ The quotient space $X
=\tX/\calg_X$ defines an orbifold (actually an isomorphism class of
orbifolds). In the case that an orbifold $X$ is given as the space of leaves of
a foliation $\calf$ on a smooth manifold, the groupoid $\calg_X$ is just the
transverse holonomy groupoid of $\calf.$

Next we review Haefliger's construction [Hae] of the classifying space of an
orbifold $X$ of dimension $n.$ Let $\{\tU_i,\grG_i,\phi_i\}$ be the local
uniformizing systems of $X,$ and consider again the disjoint union $\tX= \sqcup
\tU_i.$ Let $\calg_X$ denote the groupoid of $X,$ that is the groupoid
generated by germs of diffeomorphisms $g_{ij}:\tU_j\ra{1.3} \tU_i$ and germs of
diffeomorphisms of the action of the finite groups $\grG_i.$ Consider the
bundle of linear frames $\tP$ over $\tX.$ The groupoid $\calg_X$ acts freely on
$\tP$ with quotient space $P$ the frame bundle of $X.$ For each positive
integer $N$ let $V^{(N)}=GL(n+N)/GL(N)$ denote the Stiefel manifold with the
standard action of $GL(n).$ $V^{(N)}$ is $N$-universal as a principal $GL(n)$
bundle. The nested sequence $\cdots \subset V^{(N)}\subset V^{(N+1)}\cdots$
gives rise the direct limit $EGL(n).$ We form the associated bundles
$EX^{(N)}=\tP\times_{GL(n)}V^{(N)}$ whose direct limit we denote by $EX.$ The
groupoid $\calg_X$ acts freely on $EX^{(N)}$ for each $N$ and thus on $EX.$ So
$EX$ is universal for the groupoid $\calg_X.$ We denote the quotient
$EX/\calg_X$ by $BX,$ and the finite pieces $EX^{(N)}/\calg_X$ by $BX^{(N)}.$
We also have that $EX=\tP\times_{GL(n)}EGL(n)$ and
$$BX=\tP\times_{GL(n)}EGL(n)/\calg_X = P\times_{GL(n)}EGL(n). \leqno{\secl.1}$$
Furthermore, there are natural projections
$$\matrix{&&BX&&\cr
          &\swarrow && \searrow &\cr
          X&&&& BGL(n).\cr} \leqno{\secl.2}$$
The right arrow is a fibration, but the left arrow $p:BX\ra{1.3} X$ is not. It
has generically contractible fibers and encodes information from the local
uniformizing groups on the singular strata. It is this map that is of most
interest. Haefliger defines the
orbifold cohomology, homology, and homotopy groups by
$$H^i_{orb}(X,\bbz)=H^i(BX,\bbz), \quad H^{orb}_i(X,\bbz)=H_i(BX,\bbz), \quad
\pi_i^{orb}(X)=\pi_i(BX). \leqno{\secl.3}$$
This definition of $\pi_1^{orb}$ is equivalent to Thurston's better known
definition [Thu] in terms of orbifold deck transformations, and when $X$ is a
smooth manifold these orbifold groups coincide with the usual groups.

In the category of orbifolds the concept of a bundle is replaced by that
of a V-bundle. This consists of bundles $B_{\tU_i}$ over the local uniformizing
neighborhoods $\tU_i$ that patch together in a certain way. In particular,
there are group homomorphisms $h_{\tU_i}:\grG_i\ra{1.3} G,$ where $G$ is the
group of the V-bundle that satisfy the condition that:
\item{\secl.4} If $g_{ji}:\tU_i \ra{1.5}\tU_j$ is a diffeomorphism onto an open
set, then there is a bundle map $g_{ij}^*:B_{U_j}|g_{ji}(\tU_i)\ra{1.3}
B_{\tU_i}$ satisfying the condition that if $\grg\in \grG_i,$ and $\grg'\in
\grG_j$ is the unique element such that $g_{ji}\circ \grg=\grg'\circ g_{ji},$
then $h_{\tU_i}(\grg)\circ g_{ji}^*=g_{ji}^*\circ h_{\tU_j}(\grg').$

\noindent
If the fiber $F$ is a vector space and $G$ acts on $F$ as linear
transformations of $F$, then the V-bundle is called a {\it vector}
V-bundle. Similarly, if $F$ is the Lie group $G$ with its right action, then
the V-bundle is called a {\it principal} V-bundle.
The {\it total space} of a V-bundle over $X$ is an orbifold $E$ with local
uniformizing systems $\{B_{\tU_i},\grG_i^*,\varphi_i^*\}.$ By choosing the
local uniformizing neighborhoods of $X$ small enough, we can always take
$B_{\tU_i}$ to be the product $\tU_i\times F$ which we shall heretofore assume.
There is an action of the local uniformizing group $\grG_i$ on $\tU_i\times F$
given by sending $(\tx_i,b)\in \tU_i\times F$ to
$(\grg^{-1}\tx_i,bh_{\tU_i}(\grg)),$ so the local uniformizing groups
$\grG_i^*$ can be taken to be subgroups of $\grG_i.$ We are particularly
interested in the case of a principal bundle. In the case the fiber is the Lie
group $G,$ so the image $h_{\tU_i}(\grG_i^*)$ acts freely on $F.$ Thus the
total space $P$ of a principal V-bundle will be smooth if and only if
$h_{\tU_i}$ is injective for all $i.$
We shall often denote a V-bundle by the
standard notation $\pi:P\ra{1.3} X$ and think of this as an ``orbifold
fibration''. It must be understood, however, that an orbifold fibration is not
a fibration in the usual sense. Shortly, we shall show that it is a fibration
rationally.

We want to have the analogue of an ``atlas of charts'' on $BX.$ Let
$\{\tU_i,\grG_i,\phi_i\}$ be a cover of uniformizing charts for the orbifold
$X,$ and let $P_i$ denote the linear frame bundle over $U_i.$ The group $GL(n)$
acts locally freely on $P_i$ with isotropy group $\grG_i$ fixing the frames
over the center $a_i\in U_i.$ So we have homeomorphisms $P_i/GL(n)\approx
\tU_i/\grG_i \approx U_i.$ Thus, we can cover $BX$ by neighborhoods of the form
$\tU_i\times_{\grG_i}EGL(n)$ where $\grG_i$ is viewed as a subgroup of $GL(n).$
In fact $\grG_i$ can always be taken as a subgroup of $O(n).$ Now by refining
the cover if necessary we have injection maps $g_{ji}:\tU_i\ra{1.3} \tU_j$ and
these induce the change of ``charts'' maps
$$G_{ji}:\tU_i\times_{\grG_i}EGL(n) \ra{1.3} \tU_j\times_{\grG_j}EGL(n)
\leqno{\secl.5}$$
given by $G_{ji}([\tx_i,e])=[g_{ji}(\tx_i),e].$ This is well-defined since the
unique element $\grg_j$ is identified with $\grg_i$
under the identification of $\grG_i$ as a subgroup of $\grG_j.$ This will allow
us to construct local data on $\tU_i\times_{\grG_i}EGL(n)$ by considering
smooth (or holomorphic) data on the $\tU_i$ and continuous data on $EGL(n)$
which is invariant under the $\grG_i$ action. Since the $g_{ji}$ are
diffeomorphisms (or biholomorphisms) this local data will then patch to give
global data on $BX.$ For example we denote by $\cale$ the sheaf of germs of
complex-valued functions on $BX$ that are smooth in $\tU_i$ and continuous in
$EGL(n).$ We shall call global (local) sections of the sheaf $\cale$ smooth
functions on $BX.$ Similarly we refer to ``smooth'' maps from $BX.$ For
example, if $f:X_1\ra{1.3} X_2$ is a smooth map of orbifolds, then this induces
a smooth map $Bf:BX_1\ra{1.3} BX_2.$ Similarly the map $p:BX\ra{1.3} X$
is smooth, since its local covering maps are smooth.

Next we give a characterization of V-bundles.

\noindent{\sc Theorem} \secl.6: \tensl Let $X$ be an orbifold. There is a
one-to-one correspondence between isomorphism classes of V-bundles on $X$ with
group $G$ and generic fiber $F$ and isomorphism classes of bundles on $BX$ with
group $G$ and fiber $F.$ \tenrm

\noindent{\sc Proof}: A V-bundle on $X$ is a bundle on $\tU_i$ for each local
uniformizing neighborhood $\tU_i$ together with a group homomorphism
$h_{\tU_i}\in Hom(\grG_i,G)$ that satisfy the compatibility condition
\secl.4. This gives an action of $\grG_i$ on $B_{\tU_i}.$ Now cover
$BX$ by neighborhoods of the form $\tU_i\times_{\grG_i}EGL(n)$ where we make
use of the fact that the local uniformizing groups $\grG_i$ can be taken as
subgroups of $O(n)\subset GL(n).$ There is an action of $\grG_i$ on
$\tU_i\times F\times EGL(n)$ given by $(\tx_i,u,e)\mapsto (\grg^{-1}\tx_i,
uh_{U_i}(\grg),e\grg),$ and this gives a $G$-bundle over
$\tU_i\times_{\grG_i}EGL(n)$ with fiber $F$ for each $i.$ Moreover, the
compatibility condition \secl.4 guarantees that these
bundles patch together to give a $G$-bundle on $BX$ with fiber $F.$

Conversely, given a $G$-bundle on $BX$ with fiber $F,$ restricting to
$\tU_i\times_{\grG_i}EGL(n)$ gives a $G$-bundle there. Since for each $i$
$\tU_i\times_{\grG_i}EGL(n)$ is the Eilenberg-MacLane space $K(\grG_i,1),$
there is a one-to-one correspondence between isomorphism classes of $G$-bundles
on $\tU_i\times_{\grG_i}EGL(n)$ and conjugacy classes of
group homomorphisms $Hom(\grG_i,G).$ The fact that these $G$-bundles come from
a global $G$-bundle on $BX$ implies that the compatibility conditions
are satisfied. The correspondence can be seen to be
bijective.  \hfill\za

In what follows we shall not distinguish between V-bundles on $X$ and bundles
on $BX.$ Thus, we have

\noindent{\sc Proposition} \secl.7: \tensl The isomorphism classes of
V-bundles on $X$ with group $G$ are in one-to-one correspondence with elements
of the non-Abelian cohomology set $H^1(BX,\calg)=H^1_{orb}(X,\calg)$ where
$\calg$ is the sheaf of germs of maps to the group $G.$ \tenrm

We now turn to the Abelian case with coefficients in a sheaf. Let $\cale$
denote sheaf of germs of smooth complex-valued functions on $BX$ and
$\cale^*$ the sub-sheaf of no-where vanishing complex valued smooth
functions. The isomorphism classes of complex line bundles on $BX,$ hence
complex line V-bundles on $X,$ are in one-to-one correspondence with the
elements of the cohomology group $H^1_{orb}(X,\cale^*).$ Now $BX$ is an
infinite paracompact CW complex, and $\cale$ is a fine sheaf, so the
exponential sequence gives an isomorphism
$$H^1_{orb}(X,\cale^*)\fract{c_1}{\ra{1.3}} H^2_{orb}(X,\bbz),
\leqno{\secl.8}$$
and we define the {\it Chern class} of a line V-bundle $\call$ to be the image
$c_1(\call)\in H^2_{orb}(X,\bbz).$ More generally one can define the Chern
classes of complex vector V-bundles by using the splitting principle.

We are particularly interested in holomorphic line V-bundles.
In [BG1] the authors introduced the group
$\hbox{Pic}^{orb}(\calz)$ of holomorphic line V-bundles on a complex
orbifold $\calz.$ On $B\calz$ define the sheaf $\cala$ of germs of holomorphic
functions on $B\calz$ to be the sub-sheaf of $\cale$ consisting of functions
that are holomorphic in the $\tU_i.$ Similarly, $\cala^*\subset$
denotes the sub-sheaf of germs of nowhere vanishing functions.
As in Proposition \secl.7 we have
$\hbox{Pic}^{orb}(\calz)\simeq H^1_{orb}(\calz,\cala^*).$
Now let $\calz$ be a compact complex Fano orbifold, that is, the anti-canonical
line V-bundle $K^{-1}_\calz$ is ample. In this case $c_1(K^{-1}_\calz)>0,$ and
as in the smooth case we have

\noindent{\sc Definition} \secl.9: \tensl The {\it index} of a Fano orbifold
$\calz$ is the largest positive integer $m$ such that ${c_1(K^{-1}_\calz) \over
m}$ is an element of $H^2_{orb}(\calz,\bbz).$ The index of $\calz$ is denoted
by $\hbox{Ind}(\calz).$ \tenrm

In view of the fact that the requisite vanishing theorem is lacking in the
singular case, we shall need

\noindent{\sc Lemma} \secl.10: \tensl Let $\calz$ be a Fano orbifold with
$\hbox{Ind}(\calz)=m.$ Then there is a holomorphic line V-bundle $\call\in
\hbox{Pic}^{orb}(\calz)$ such that $\call^m=K^{-1}_\calz.$ \tenrm

\noindent{\sc Proof}: The idea of the proof is simple, but since we are working
on $B\calz$ we write out the details. First on $B\calz$ we define the following
sheaf $\cale^{p,q}$ of ``differential forms'': Let $(z_1,\cdots,z_n)$ be
complex coordinates on $\tU_i.$ Then using the standard multi-index notation,
we construct the sheaf $\cale^{p,q}$ whose stalks are spanned by elements of
the form $f_{IJ}(z,e)dz_I\wedge d\bar{z}_J$ where $I=i_1\cdots i_p$ and $J=
j_1\cdots j_q$ are the usual multi-indices, and $f$ is a smooth function on
$\tU_i\times EGL(2n)$ satisfying $f(\grg^{-1}z,e\grg)=f(z,e)$ for $\grg\in
\grG_i.$ We have, as usual, differential operators $\bar{\partial}$ and the
Dolbeault complex
$$\cdots  \fract{\bar{\partial}}{\ra{1.5}}
\cale^{p,q}\fract{\bar{\partial}}{\ra{1.5}} \cale^{p,q+1}
\fract{\bar{\partial}}{\ra{1.5}}\cdots. \leqno{\secl.14}$$
Notice that $\cale^{0,0}=\cale$ and $\hbox{ker}(\bar{\partial}:\cale^{0,0}
\ra{1.3} \cale^{0,1})=\cala.$
By the isomorphism \secl.8 there is a complex line V-bundle $\call$ such that
$\call^m=K^{-1}.$ The transition function $h_{ij}$ for $\call$ are nowhere
vanishing local sections of $\cale$ over $\tU_i\times_{\grG_i}EGL(2n)\cap
\tU_j\times_{\grG_j}EGL(2n)$ and the transition functions for $\call^m=K^{-1}$
are $h^m_{ij}.$ But since $K^{-1}$ is holomorphic $\bar{\partial}h^m_{ij}=0;$
hence, $\bar{\partial}h_{ij}=0$ implying that $\call$ is holomorphic. \hfill\za

The groups $H^p_{orb}(\calz,\bbz)$ encode information from the local uniformizing
groups as well as the ordinary cohomology groups. In fact by studying the
Leray spectral sequence of the map $p:B\calz\ra{1.3} \calz$ it is
easy to see that these groups coincide with the ordinary groups in
certain cases. For any sheaf
$\calf$ on $B\calz$ we let $R^qp(\calf)$ denote the derived functor sheaves,
that is the sheaves associated to the presheaves $U\mapsto
H^q(p^{-1}(U),\calf).$ Then Leray's theorem says that there is a spectral
sequence $E^{p,q}_r$ with $E_2$ term given by
$E^{p,q}_2= H^p(\calz,R^qp(\calf))$
converging to $H^{p+q}_{orb}(\calz,\calf).$
If $\calf$ is one of the constant sheaves $\bbq,\bbr,\bbc$
one easily sees that spectral sequence
collapses and we recover a result of Haefliger [Hae], namely

\noindent{\sc Corollary} \secl.11: \tensl $H^q_{orb}(\calz,A)\simeq
H^q(\calz,A)$ for $A=\bbq,\bbr,\bbc$ or $\bbz_p,$ where $p$ is
relatively prime to $\hbox{Ord}(\calz).$ \tenrm

\bigskip
\centerline {\bf \S 2. Sasakian-Einstein Geometry}
\medskip

\def\tcals{\tilde{\cals}}
\def\Ric{\hbox{Ric}}

Let us begin with a brief review of properties of Sasakian-Einstein spaces.
There are several equivalent definitions of a Sasakian structure in the
literature (See [Bl, YKon]), but the most geometric definition is also
the most recent [BG2]:

\noindent{\sc Definition \seg.1}: \tensl
Let $(\cals,g)$ be a Riemannian manifold of real dimension $m$.
We say that
$(\cals,g)$ is {\it Sasakian} if the holonomy group of the
metric cone on $\cals$ \break $(C(\cals),{\bar g})=
(\bbr_+\times\cals, \ dr^2+r^2g)$ reduces to a subgroup of $U({m+1\over2})$.
In particular, $m=2n+1, n\geq1$ and $(C(\cals),\bar{g})$ is K\"ahler.
\tenrm

That this definition is equivalent to the other known definitions
of a Sasakian structure given in terms
of the triple $\{g,\xi,\Phi\},$ where $g$ is Riemannian metric, $\xi$ is a
nowhere vanishing Killing vector field, and $\Phi$ is a
tensor field on $\cals,$ was shown in our expository article [BG2].
The Killing vector field $\xi$ is called the
{\it characteristic} or {\it Reeb} vector field. The 1-form $\eta$ defined to
be the 1-form dual to $\xi$ with respect to the metric $g$ is called the {\it
characteristic 1-form} of the Sasakian structure, and
underscores the {\it contact} nature of a Sasakian structure. Indeed, a
Sasakian structure is usually defined as a normal contact metric structure
[Bl, YKon]. The terminology normal means that the almost CR-structure defined
by $\Phi$ on the orthogonal complement to subbundle of $T\cals$ defined by the
characteristic vector field $\xi$ is integrable. It is easy to generalize this
definition to that of a Sasakian {\it orbifold}.  One simply requires that
$\{g,\xi,\Phi\}$ be invariant under the action of the local uniformizing groups
of the orbifold.

We are interested in Sasakian-Einstein geometry. We have

\noindent{\sc Definition-Proposition} \seg.2: \tensl A Sasakian manifold
(orbifold) $(\cals,g,\xi,\Phi)$ is {\it Sasakian-Einstein} if its Riemannian
metric $g$ is Einstein.  The Ricci tensor $\hbox{Ric}$ of any Sasakian manifold
(orbifold) of dimension $2n+1$ satisfies
$\hbox{\Ric}(X,\xi)=2n\eta(X).$
Thus, if the metric $(\cals,g)$ is \Se, then the scalar curvature of $g$ is
positive and equals $2n(2n+1)$.

\noindent Furthermore, a Sasakian manifold $(\cals,g)$ is \Se if and only if
the cone metric $\bar{g}$ is Ricci-flat, {\it i.e.}, $(C(\cals),\bar{g})$ is
Calabi-Yau. In particular, it follows that the restricted holonomy group ${\rm
Hol}^0(\bar{g})\subset SU(n+1).$ \tenrm

Then an immediate corollary of \seg.2 and Myers Theorem is that any
complete Sasakian-Einstein manifold is compact with finite fundamental group.
In fact there is a stronger result due to Hasegawa and Seino:

\noindent{\sc Proposition} \seg.3 [HS]: \tensl Let $\cals$ be a complete
Sasakian manifold such that $\Ric(X,X)\geq \grd >-2$ for all unit vector fields
$X$ on $\cals.$ Then $\cals$ is compact with finite fundamental group. \tenrm

The case when the characteristic vector field $\xi$ of a compact Sasakian
manifold $\cals$ generates a free circle action has been well studied.  In this
case $\cals$ is the total space of a principal $S^1$-bundle whose base space
$\calz$ is a Hodge manifold; hence, a smooth projective algebraic variety. This
is a special case of the well-known Boothby-Wang fibration and is due to
Hatakeyama [Hat]. The contact form $\eta$ is nothing but a connection 1-form on
the bundle $\pi:\cals\ra{1.3} \calz$ and the curvature form $d\eta$ is just the
pullback by $\pi$ of the K\"ahler form on $\calz.$ Moreover, the fibers of
$\pi$ are totally geodesic. Now one can bring to bare O'Neill's well-known
formulae for Riemannian submersions, and these work equally as well in the
quasi-regular case discussed below.

Now generally, since the Killing vector field $\xi$ has unit norm, it defines a
1-dimensional foliation $\calf$ on $\cals,$ and we are interested in the case
that all the leaves of $\calf$ are compact. The assumption that all leaves are
compact is equivalent to the assumption that $\calf$ is {\it quasi-regular},
{\it i.e.},
each point $p\in \cals$ has a cubical neighborhood $U$ such that any leaf
$\call$ of $\calf$ intersects a transversal through $p$ at most a finite number
of times $N(p).$ Furthermore, $\cals$ is called {\it regular} if  $N(p)=1$ for
all $p\in \cals.$ In this case, the foliation $\calf$ is simple, and defines a
global submersion. In the quasi-regular case it is well-known [Tho] that
$\xi$ generates a locally free circle action on $\cals,$ and that the space of
leaves is a compact orbifold (or V-manifold)[Mol]. We shall denote the space of
leaves of the foliation $\calf$ on $\cals$ by $\calz.$ Then the natural
projection $\pi:\cals\ra{1.3} \calz$ is an orbifold submersion and a Siefert
fibration. Actually, much more is true.

\noindent{\sc Theorem \seg.4}: \tensl
Let $(\cals,g)$ be a compact quasi-regular Sasakian manifold of
dimension $2n+1$, and let $\calz$ denote the space of leaves of the
characteristic foliation. Then
\item{(i)}The leaf space $\calz$ is a compact complex orbifold with a
K\"ahler metric $h$ and K\"ahler form $\gro$ which defines an integral class
$[\gro]$ in $H^2_{orb}(\calz,\bbz)$ in such a way that $\pi:(\cals,g) \ra{1.3}
(\calz,h)$ is an orbifold Riemannian submersion. The fibers of $\pi$ are
totally geodesic submanifolds of $\cals$ diffeomorphic to $S^1.$
\item{(ii)} $\calz$ is also a normal projective algebraic variety which is
$\bbq$-factorial.
\item{(iii)} The orbifold $\calz$ is Fano if and only if $\Ric_g >-2.$ In
this case $\calz$ as a topological space is simply connected, and as an
algebraic variety is uniruled with Kodaira dimension $\kappa(\calz)=-\infty.$
\item{(iv)} $(\cals,g)$ is Sasakian-Einstein if and only if $(\calz,h)$ is
K\"ahler-Einstein with scalar curvature $4n(n+1).$
\tenrm

\noindent{\sc Proof}: (i): It is well-known [Mol] that $\calz$ is a compact
orbifold, and the remainder follows exactly as in the regular case except that
now the K\"ahler class $[\gro]$ defines an integral class in
$H^2_{orb}(\calz,\bbz)$ which is generally only a rational class in the
ordinary cohomology $H^2(\calz,\bbq).$ This proves (i). To prove (ii) we notice
that as the usual case $\gro$ is a positive $(1,1)$-form and thus represents a
holomorphic line V-bundle on $\calz.$ (ii) now follows from the Kodaira-Baily
embedding theorem [Ba2].

In the Sasakian case the O'Neill tensors $T$ and $N$ (See [Bes]) vanish,
and the tensor $A$ satisfies $A_XY=-g(\Phi X,Y)\xi$ and $A_X\xi=\Phi X$ for $X$
and $Y$ horizontal. Then one easily sees from [Bes, 9.36] that
$$\Ric_g(X,Y)=\Ric_h(X,Y) -2g(X,Y). \leqno{\seg.5}$$
The first statement of (iii) now follows from this. Moreover, combining \seg.5
with the last statement of \seg.2 implies (iv). In (iii) the simple
connectivity of $\calz$ is basically Kobayashi's argument [Kob1] with the usual
Riemann-Roch replaced by the singular version of Baum, Fulton, and MacPherson
(see [BG1] for details). The uniruledness and Kodaira dimension follow from
Miyaoka and Mori [MiMo]. \hfill\za

Next we give an inversion theorem to Theorem \seg.4. In the regular case this
goes back to a construction of Kobayashi [Kob2] together with Hatakeyama [Hat].
There is also a description in the context of Sasakian-Einstein geometry in
[FrKat2,BFGK] in the regular case. The non-regular case follows by applying
Hatakeyama's results to our more general setup. First we have

\noindent{\sc Definition} \seg.6: \tensl A compact K\"ahler orbifold is called
a {\it Hodge orbifold} if the K\"ahler class $[\gro]$ lies in
$H^2_{orb}(\calz,\bbz).$ \tenrm

With this definition a restatement of the Kodaira-Baily embedding theorem [Ba2]
is:

\noindent{\sc Theorem} [Baily]\seg.7: \tensl A Hodge orbifold is a projective
algebraic variety. \tenrm

Now our inversion theorem is:

\noindent{\sc Theorem} \seg.8: \tensl Let $(\calz,h)$ be a Hodge orbifold.
Let $\pi:\cals\ra{1.3} \calz$ be the $S^1$ V-bundle whose first Chern
class is $[\gro],$ and let $\eta$ be a connection 1-form in $\cals$ whose
curvature is $2\pi^*\gro,$ then $(\cals,\eta)$ with the metric
$\pi^*h+\eta\otimes\eta$ is a Sasakian orbifold.  Furthermore, if all the local
uniformizing groups inject into the group of the bundle $S^1,$ the total space
$\cals$ is a smooth Sasakian manifold. \tenrm

\noindent{\sc Remark} \seg.9: The orbifold structure of $\calz$ is crucial
here. Consider the weighted projective space $\bbc\bbp^2(p_1,p_2,p_3)$ defined
by the usual weighted $\bbc^*$ action on $\bbc^3-\{0\}$ where $p_1,p_2,p_3$ are
pairwise relatively prime integers. As an algebraic variety
$\bbc\bbp^2(p_1,p_2,p_3)$ is equivalent [Kol] to
$\bbc\bbp^2/\bbz_{p_1}\times\bbz_{p_2}\times\bbz_{p_3}.$ But as orbifolds these
are distinct, since the former has $\pi_1^{orb}=0,$ whereas the latter has
$\pi_1^{orb}=\bbz_{p_1}\times\bbz_{p_2}\times\bbz_{p_3}.$ This is to be
contrasted with the smooth case where equivalence as complex manifolds
coincides with equivalence as algebraic varieties. Of course, the metrics are
also different. In the latter case the metric is just the Fubini-Study metric
pushed to the quotient which is K\"ahler-Einstein. The Sasakian structure on
the corresponding V-bundle is just the standard Sasakian-Einstein structure on
the lens space $S^5/\bbz_{p_1}\times\bbz_{p_2}\times\bbz_{p_3}.$ In the former
case the Sasakian structure is one of the non-standard deformed structures on
$S^5$ described in Example 7.1 of [YKon] which is ellipsoidal and not Einstein.
Likewise, the metric $h$ on $\bbc\bbp^2(p_1,p_2,p_3)$ is not K\"ahler-Einstein.

We are particularly interested in constructing simply connected
Sasakian-Einstein manifolds.  Thus,

\noindent{\sc Corollary} \seg.10: \tensl Let $\calz$ be a compact Fano orbifold
with $\pi_1^{orb}(\calz)=0.$ Let $\pi:\cals\ra{1.3} \calz$ be the $S^1$
V-bundle whose first Chern class is ${c_1(\calz)\over \hbox{Ind}(\calz)}.$
Suppose further that the local uniformizing groups of $\calz$ inject into
$S^1.$ Then there is a metric $g$ on the total space $\cals$ such that $\cals$
is a compact simply connected Sasakian manifold with $\Ric_g>-2.$ Furthermore,
if $\calz$ is K\"ahler-Einstein then $(\cals,g)$ is Sasakian-Einstein. \tenrm

\noindent{\sc Proof}: The only part that does not follow immediately from our
results is the simple connectivity. Suppose $\cals$ were not simply connected,
then by compactness and the bound on the Ricci tensor there would be at most a
finite cover $\tilde{\cals}.$ But since $\pi_1^{orb}(\calz)=0$ there is no
nontrivial cover of $\calz.$ So $\tilde{\cals}$ must be the total space of a
V-bundle on $\calz$ covering $\cals.$ But then we must have
$\tilde{\cals}=\cals$ since the first Chern class of the V-bundle
$\pi:\cals\ra{1.3}\calz$ is not divisible in $H^2_{orb}(\calz,\bbz)$ by the
definition of $\hbox{Ind}(\calz).$ \hfill\za

\noindent{\sc Definition \seg.11}: \tensl Let $(\cals,g)$ be a compact
quasi-regular Sasakian orbifold with $\calz$ the space of leaves of the
characteristic foliation. We define the {\it order} of $\cals,$ denoted by
$\hbox{Ord}(\cals)$ to be the order of the orbifold $\calz.$ When $\Ric_g>-2,$
we define the {\it index} of $\cals,$ denoted $\hbox{Ind}(\cals),$ to be the
index of $\calz.$  \tenrm

\noindent{\sc Warning}: The order of a quasi-regular Sasakian orbifold defined
here is {\sc not} the order of $\cals$ as an orbifold. There are many smooth
Sasakian manifolds with large order, whereas $\hbox{Ord}(\cals)=1$ means that
$\cals$ is smooth and regular.

Notice that the index is always defined for any quasi-regular \Se orbifold.
Both the index $\hbox{Ind}(\cals)$ and $\hbox{Ord}(\cals)$ are invariants of
the Sasakian structure on $\cals.$ If $\cals$ is a 3-Sasakian manifold the
space $\calz$ is independent, up to isomorphism, of the choice of
characteristic vector field $\xi$ in the Lie algebra $\gs\gu(2).$ Thus, it
makes sense to talk of both the index and order of a 3-Sasakian manifold.
Moreover, from the existence of the contact line V-bundle on $\calz$ one has
$\hbox{Ind}(\cals)\geq n+1$ where $\hbox{dim}~\cals =4n+3.$ More generally, one
easily sees that for the standard Sasakian structure on $S^{2n+1}$ that
$\hbox{Ind}(S^{2n+1})= n+1.$ Note that $\hbox{Ind}(S^{4n+3})= 2n+2$, and it was
shown in [BG1] using Kawasaki's Riemann-Roch theorem for orbifolds [Kaw] that
this is the only simply connected 3-Sasakian manifold with index $2n+2.$
Furthermore, this is the largest possible index of a 3-Sasakian manifold.
Summarizing we have

\noindent{\sc Proposition} \seg.12: \tensl The index $\hbox{Ind}(\cals)$ of a
3-Sasakian manifold $\cals$ of dimension $4n+3$ is either $n+1$ or $2n+2$ and
equals $2n+2$ if and only if the universal cover of $\cals$ is $S^{4n+3}.$
\tenrm

More detailed information about the index is available in the regular case
which is discussed in the next section. Finally we end this section by
mentioning a related result of Vaisman [Vai,DO]. If $\cals$ is any Sasakian
manifold then $\cals\times S^1$ is locally conformal K\"ahler with parallel
Lee form [Vai]. Such manifolds have been called generalized Hopf manifolds or
Vaisman manifolds [DO]. Moreover, the universal cover of every generalized Hopf
manifold is of the form $\cals\times \bbr$ where $\cals$ is Sasakian. The \Se
manifolds $\cals$ discussed here give a subclass of generalized Hopf manifold
with the property that the local K\"ahler metrics are also Ricci flat. Thus,
one might refer to this subclass of manifolds as locally conformal Calibi-Yau
manifolds. There are obvious translations of our results on \Se manifolds to
the class of locally conformal Calabi-Yau manifolds.

bigskip
\centerline{\bf \S 3. Regular \Se manifolds}
\medskip

The following is a translation to regular Sasakian-Einstein
geometry of known results [Ko, K-O, Wi1, Wi2] about the index of smooth Fano
varieties:

\noindent{\sc Theorem} \ser.1: \tensl Let $\cals$ be a regular
Sasakian-Einstein manifold of dimension $2n+1.$ Then
\item{(i)} $\hbox{Ind}(\cals)\leq n+1.$
\item{(ii)} If $\hbox{Ind}(\cals)=n+1$ then the universal cover $\tilde{\cals}$
of $\cals$ is  $S^{2n+1}$ with its standard Sasakian structure.
\item{(iii)} If $\hbox{Ind}(\cals)=n$ then $\tilde{\cals}$ is a circle bundle
over the complex quadric $Q_n(\bbc).$
\item{(iv)} if $r=\hbox{Ind}(\cals)\geq {n+2\over 2},$ then $b_2(\cals)=0$
unless $\cals$ is a circle bundle over $\bbc\bbp^{r-1}\times \bbc\bbp^{r-1}.$
\item{(v)} If $r=\hbox{Ind}(\cals)= {n+1\over 2},$ then $b_2(\cals)=0$
unless $\cals$ is a circle bundle over $\calz$ where $\calz$ is either
$\bbc\bbp^{r-1}\times Q_r(\bbc),$ $\bbp(T^*\bbc\bbp^r),$ or possibly
$\bbc\bbp^{2r-1}$ blown up along a $\bbc\bbp^{r-2}.$ \tenrm

Items (i)-(iii) follow from results of Kobayashi and Ochiai [KO], while items
(iv) and (v), due to Wi\'sniewski [Wi1, Wi2], use Mori theory. The examples in
[BGM1] and [BGMR] show that (iv) and (v) cannot hold in the non-regular case.

For dimensional reasons we have immediately the following

\noindent{\sc Corollary} \ser.2: \tensl Let $\cals$ be a regular
Sasakian-Einstein manifold of dimension $4n+3.$ Then if $\hbox{Ind}(\cals)>
n+1$, we have $b_2(\cals)=0.$  \tenrm

A related result follows from the boundedness of deformation types of smooth
Fano varieties [Kol]. This implies that in any given dimension the Betti
numbers of a regular \Se manifold are bounded by constants only depending on
the dimension. In particular, the bound on the second Betti number for regular
\Se 5-manifolds is $8,$ whereas for regular \Se 7-manifolds it is $9.$ These
bounds are actually sharp and are discussed in Theorem \ser.5 and Proposition
\sewz.12 below, respectively. The results of [BGMR] show that these bounds no
longer hold in the non-regular case.

A special case of regular Sasakian-Einstein manifolds are the homogeneous ones.
Recall the following well-known terminology. Let $G$ be a complex semi-simple
Lie group. A maximal solvable complex subgroup $B$ is called a Borel subgroup,
and $B$ is unique up to conjugacy.  Any complex subgroup $P$ that contains $B$
is called a parabolic subgroup. Then the homogeneous space $G/P$ is called a
generalized flag manifold. A well-known result of Borel and Remmert [BR] says
that every connected homogeneous K\"ahler manifold is a product of a torus and
a generalized flag manifold.

\noindent{\sc Definition} \ser.3: \tensl A Sasakian manifold $(\cals,g)$
is called a {\it homogeneous Sasakian} manifold if there is a transitive
group $K$ of isometries on $\cals$ that preserve the Sasakian structure, that
is, if $\phi^k\in \hbox{Diff}~\cals$ corresponds to $k\in K,$ then
$\phi^k_*\xi =\xi.$ (This implies that both $\Phi$ and $\eta$ are also
invariant under the action of $K.$) \tenrm

Note that if $\cals$ is compact then $K$ is a compact Lie group.

\noindent{\sc Theorem} \ser.4: \tensl Let $(\cals,g')$ be a complete
homogeneous Sasakian manifold with $\Ric_{g'}\geq \grd >-2.$
Then $(\cals,g')$ is a compact regular homogeneous Sasakian manifold,
and there is a homogeneous Sasakian-Einstein metric $g$ on $\cals$
that is compatible with the underlying normal contact structure.
Moreover, $\cals$ is the total space of an $S^1$-bundle over
a generalized flag manifold $G/P$ whose Chern class is $\displaystyle{m\cdot
{c_1(G/P)\over \hbox{Ind}(G/P)}}$
for some positive integer $m.$ Conversely, given
any generalized flag manifold $G/P$ and any positive
integer $m,$ the total space $\cals_m$
of the circle bundle $\pi: \cals_m \ra{1.3} G/P$ whose
Chern class is given by the above formula  has
a homogeneous Sasakian-Einstein metric $g.$  \tenrm

\noindent{\sc Proof}: By Boothby and Wang [BW] $\cals$ is regular, and
by \seg.4 and \seg.5, $\cals$ is compact and fibers over a smooth
simply connected Fano variety $\calz.$  Moreover, since
there is a transitive group of isometries on
$(\cals,g')$ that commutes with the $S^1$ action generated by $\xi,$ there is
a positive K\"ahler metric $h'$ on $\calz$ together with a transitive group
of isometries $K'$ preserving the K\"ahler structure. It follows from Borel and
Remmert [BR] that $\calz$ is a generalized flag manifold $G/P$ where the
complex semi-simple Lie group $G$ is the complexification of $K'$
and $P\subset G$  is a parabolic subgroup.
Then by a well-known result of Matsushima [Bes] there is a positive
K\"ahler-Einstein metric $h$ on $G/P$ and an automorphism $\phi$ of
the complex structure such that $h=\phi^*h'.$ Furthermore, the group
$K=\phi^{-1}K'\phi$ acts as a transitive group of isometries on
$(G/P,h).$ Thus, by the Kobayashi construction (cf. Corollary \seg.10) there
is a \Se metric $g$ on $\cals$ that is compatible with the
normal contact structure, and the isometry group $K$ lifts to a group
$\check{K}$ on $\cals$ that leaves $\pi^*h$ invariant and commutes with the
circle group $S^1$ generated by $\xi.$ It follows after adjusting the scale
of $h$ that $\pi^*h+\eta\otimes \eta$ is a \Se metric on $\cals.$

The statement about the Chern class follows from Kobayashi's
construction, and the converse follows easily as well.  \hfill\za

The classification of all simply connected regular
Sasakian-Einstein manifolds depends on the classification of all smooth Fano
varieties with a K\"ahler-Einstein metric. This is a deep and important problem
(the well-known Calabi problem for $c_1$ positive) which has recently met with
a great deal of success [TY, Ti1, Ti3, Sui] but whose complete solution is still
at large. The resolution of this problem would also describe all regular \Se
manifolds. Indeed, until recently there was a folklore conjecture that stated
that any (smooth) Fano variety with no holomorphic vector fields admits a
K\"ahler-Einstein metric.  This was first shown to be false in the orbifold
category by Ding and Tian [DT] and more recently in the smooth manifold
category by Tian [Ti3]. (The folklore conjecture is true in the case of smooth
del Pezzo surfaces [TY, Ti2]).

In dimension 5, Friedrich and Kath [FrKat1, BFGK] reduce
the classification of compact regular Sasakian-Einstein 5-manifolds
to the classification of Fano surfaces with a K\"ahler-Einstein
metric [TY, Ti2]. So Tian's later result [Ti2] actually proves a
stronger theorem than appears in [FrKat1, BFGK].

\noindent{\sc Theorem} \ser.5: \tensl Let $(\cals,g)$ be a regular
Sasakian-Einstein 5-manifold. Then $\cals =\tcals/\bbz_m$ where the universal
cover $(\tcals,g)$ is precisely one of the following:
\item{(1)} $S^5$ with its standard (Sasakian) metric.
\item{(2)} The Stiefel manifold $V_2(\bbr^4)$ of 2-frames in $\bbr^4$ with the
unique Sasakian metric that is a Riemannian submersion over $\bbc\bbp^1\times
\bbc\bbp^1$ with its standard complex structure and symmetric metric.
\item{(3)} The total space $S_k$ of the $S^1$ bundles $S_k\rightarrow P_k$ for
$3\leq k\leq 8,$ where $P_k$ is $\bbc\bbp^2\#k\overline{\bbc\bbp^2}$ with any
of its complex structures for which $c_1$ is positive. Moreover, for each such
complex structure there is a unique \Se metric $g$ on $S_k.$ \tenrm

We denote the circle bundle $S_k$ with a fixed \Se structure by $\cals_k.$ A
well-known result of Smale says that any simply connected 5-manifold with $H_2$
torsion free is diffeomorphic to the connected sum $S^5\# k(S^2\times S^3)$
for some nonnegative integer $k,$ which corresponds to the $k$ above.  Indeed,
Theorem \ser.5 implies that there is a regular \Se structure on each of these
manifolds for $k\leq 8$ with the exception of $k=2.$ Furthermore, for $5\leq
k\leq 8$ there are continuous families of \Se structures, which are
inequivalent (i.e. non-trivial deformations) since the Tian
families on $P_k$ are inequivalent. This should be contrasted
with the 3-Sasakian case which is infinitesimally rigid [PP].  Let
$\calr_5(k)$ denote the set of regular \Se structures on $S^5\# k(S^2\times
S^3).$ Then summarizing we have

\noindent{\sc Corollary} \ser.6: \tensl
\item{(1)} $\calr_5(k)$ is empty for $k=2$ and $k\geq 9.$
\item{(2)} $\calr_5(k) = \{\hbox{point}\}$ for $k=0,1,3,4.$
\item{(3)} $\hbox{Dim}~\calr_5(k)\geq k-4$ for $5\leq k\leq 8.$
\tenrm

An appropriate topology on $\calr_5(k)$ will be discussed in section 5, so the
concept of dimension is valid. Another corollary of [Ti2] is

\noindent{\sc Corollary} \ser.7: \tensl Let $\cals$ be a regular Sasakian
5-manifold and let $\calz$ denote the space of leaves of the characteristic
foliation.  Suppose also that $\calz$ is Fano. Then, $\cals$ admits
a compatible \Se metric if and only if the Lie algebra $\ga\gu\gt(\calz)$ of
holomorphic vector fields on $\calz$ is reductive. \tenrm

We finish this section with several examples in higher dimensions. These are
$S^1$ bundles over Fano manifolds known to have K\"ahler-Einstein metrics.

\noindent{\sc Example} \ser.8: {\it Fermat hypersurfaces.} Consider the Fermat
hypersurfaces $F_{d,n+1}$ of degree $d$ in $\bbp^{n+1}$ given in homogeneous
coordinates by
$$z_0^d+\cdots +z_{n+1}^d=0.\leqno{\ser.9}$$
They are Fano for $d\leq n+1$ and Nadel has shown that $F_{d,n+1}$ has a
K\"ahler-Einstein metric for ${n+1\over 2}\leq d\leq n+1.$ (It was shown
previously in [Sui] and [Ti1] that $F_{n,n+1}$ and $F_{n+1,n+1}$ are
K\"ahler-Einstein.) Let $\cals_{d,n+1}$ denote the $S^1$ bundle over
$F_{d,n+1}$ determined by Corollary \seg.10. The index of $F_{d,n+1},$ hence of
$\cals_{d,n+1},$ is $n+2-d.$ Now the topology of smooth hypersurfaces is
well-known [Dim]. The homology is torsion free
and by the famous Lefschetz theorem
$H_*(F_{d,n+1},\bbz)=H_*(\bbp^n,\bbz)$ except in the middle dimension $n$ where
it is determined by its degree [Dim]. Moreover, $\cals_{d,n+1}$ is just the
link of the hypersurface \ser.9 which is well-known to be $(n-1)$-connected
[Dim]. Thus, from the Serre spectral sequence of the circle bundle one finds
$$H^i(\cals_{d,n+1},\bbq)= \cases{\bbq &if $i=0,2n+1$;\cr
       \bbq^{b_n} &if $i=n,n+1$;\cr
       0 &otherwise,\cr}$$
where the $n$th Betti number is given by
$$b_n=b_n(\cals_{d,n+1})=(-1)^n\Bigl(1+{(1-d)^{n+2}-1\over d}\Bigr).$$
For example in dimension 7 ($n=3$) we have three \Se manifolds, namely
$\cals_{4,4},\cals_{3,4},$ and $\cals_{2,4}.$ We find
$$b_3(\cals_{4,4})=60, \quad b_3(\cals_{3,4})=10, \quad b_3(\cals_{2,4})=0.$$
To the best of the authors' knowledge the first two give the first \Se
manifolds with $b_3\neq 0.$ Thus, for topological reasons both $\cals_{4,4}$
and $\cals_{3,4}$ cannot admit a 3-Sasakian structure. More generally the \Se
manifold $\cals_{n+2,2n+2}$ has the same dimension and
index as a 3-Sasakian manifold, but it cannot admit a 3-Sasakian structure
since $b_{2n+1}>0.$ Notice that $F_{2,n+1}$ is just the quadric (which is
homogeneous and admits a K\"ahler-Einstein metric) discussed previously, and
that if $n$ is odd it has the cohomology groups of $\bbc\bbp^n,$ but differs in
the ring structure.  Notice also that generally $\cals_{d,n+1}$ is just the
Brieskorn manifold described by Equation \ser.9 as a submanifold of $S^{2n+1}.$

A similar analysis can be used to discuss other examples of circle bundles over
known K\"ahler-Einstein Fano manifolds. For example Nadel [Na] shows the
existence of K\"ahler-Einstein metrics on certain Fano complete intersections.
Also Koiso and Sakane [Sa, KS1, KS2] proved the existence of K\"ahler-Einstein
metrics on certain almost homogeneous Fano manifolds,
and in particular toric K\"ahler-Einstein Fano manifolds [Mab]. We intend to
study the \Se circle bundles over these manifolds elsewhere.

\bigskip
\centerline{\bf \S 4. The \Se Monoid}
\medskip
In this section we apply a construction due to Wang and Ziller to define a
multiplication on the set of quasi-regular Sasakian-Einstein orbifolds.

\noindent{\sc Definition} \sewz.1: \tensl We denote by $\cals\cale$ the set of
compact quasi-regular \Se orbifolds with $\pi_1^{orb}=0,$ by $\cals\cale^s$ the
subset of $\cals\cale$ that are smooth manifolds, and by $\calr\subset
\cals\cale^s$ the subset of compact, simply connected, regular \Se manifolds.
The set $\cals\cale$ is topologized with the $C^{m,\gra}$ topology, and the
subsets are given the subspace topology. \tenrm

The condition $\pi_1^{orb}=0$ is made to avoid complications, and since every
\Se orbifold is covered by one with $\pi_1^{orb}=0,$ there is no loss of
generality. The set $\cals\cale$ is graded by dimension, that is,
$$\cals\cale =\bigoplus_{n=0}^\infty \cals\cale_{2n+1},$$
and similarly for $\cals\cale^s$ and $\calr.$ In the definition of Sasakian
structure it is implicitly assumed that $n>0.$ So we want to extend the
definition of a \Se structure to the case when $n=0.$ This can easily be done
since a connected one dimensional orbifold is just an interval with possible
boundary, or a circle.  So we can just take $\xi ={\partial\over \partial t},
\eta =dt,\Phi=0,$ with the flat metric $g=dt^2.$ In this case the space of
leaves $\calz$ of the characteristic foliation is just a point. The unit circle
$S^1$ with this structure will play the role of the identity in our monoid. We
define the {\it index} of any element of $\cals\cale_1$ to be $0$ and
understand that $\gcd(0,m)=m$ for any integer $m.$ Then we see that
$\cals\cale^s_1=\calr_1$ is a single point, namely the
unit circle with its flat metric. Also it follows from a theorem of Hamilton
[Bes] that $\cals\cale^s_3=\calr_3=\{pt\},$ namely $S^3$
with its standard \Se structure.  Moreover, we know of no examples of elements
of $\cals\cale^s_5$ that are not in $\calr_5.$

We now define a graded multiplication
$$\cals\cale_{2n_1+1}\times \cals\cale_{2n_2+1}\ra{1.5}
\cals\cale_{2(n_1+n_2)+1}$$
as follows: Let $\cals_1,\cals_2\in \cals\cale$ of dimension $2n_1+1$ and
$2n_2+1$ respectively. Their respective space of leaves $\calz_1$ and $\calz_2$
are K\"ahler-Einstein Fano orbifolds of complex dimension $n_1,n_2,$ and
metrics $h_1,h_2,$ respectively. Moreover, the scalar curvature of
$(\calz_i,h_i)$ is $4n_i(n_i+1).$ Now the product orbifold $(\calz_1\times
\calz_2,h_1+h_2)$ is K\"ahler, but in general not K\"ahler-Einstein. However,
by rescaling we see that the metric
$$h'={(n_1+1)h_1+(n_2+1)h_2\over n_1+n_2+1} \leqno{\sewz.2}$$
is K\"ahler-Einstein and Fano with scalar curvature $4(n_1+n_2)(n_1+n_2+1).$
Furthermore, if $\hbox{Ind}(\calz_i)$ denotes the indices of $\calz_i,$ then
the index of $\calz_1\times \calz_2$ is just
$\gcd(\hbox{Ind}(\calz_1),\hbox{Ind}(\calz_2)).$ Thus, by Corollary \seg.10 the
$S^1$ V-bundle on $\calz_1\times \calz_2$ whose first Chern class is
$${c_1(\calz_1\times \calz_2)\over \hbox{Ind}(\calz_1\times \calz_2)}=
{c_1(\calz_1)+c_1(\calz_2)\over \gcd(\hbox{Ind}(\calz_1),\hbox{Ind}(\calz_2))}
\leqno{\sewz.3}$$
is simply connected in the orbifold sense with a \Se
structure determined by Theorem \seg.8. We denote this multiplication by
$\star,$ and refer to it as the {\it join}. It is clear from the construction
that if $\cals_1,\cals_2\in \cals\cale$ then also $\cals_1\star \cals_2\in
\cals\cale$ with the \Se structure as defined above, and that this procedure
can be iterated in such a way that $\star$ is associative. It is also
commutative up to isomorphism.  Furthermore, $\star$ is continuous in both
factors. We have arrived at:

\noindent{\sc Theorem} \sewz.4: \tensl The operation $\star$ defined above
gives $\cals\cale$ the structure of a commutative associative topological
monoid.  \tenrm

Actually we are interested in the subset $\cals\cale^s$ of smooth \Se
manifolds. However, it is certainly not true that $\cals\cale^s$ is a
submonoid. We want conditions that guarantee that if $\cals_i\in \cals\cale^s$
for $i=1,2$ then $\cals_1\star \cals_2\in \cals\cale^s.$ For $i=1,2$ we define
the {\it relative indices} of the pair $(\cals_1,\cals_2)$ by
$$l_i= {\hbox{Ind}(\cals_i)\over
\gcd(\hbox{Ind}(\cals_1),\hbox{Ind}(\cals_2))}.\leqno{\sewz.5}$$
Then, $\gcd(l_1,l_2)=1,$ and we have,

\noindent{\sc Proposition} \sewz.6: \tensl  Let $\cals_i\in \cals\cale^s$ with
orders $m_i,$ and relative indices $l_i,$ respectively. Then
\item{(i)} $\hbox{Ind}(\cals_1\star
\cals_2)=\gcd(\hbox{Ind}(\cals_1),\hbox{Ind}(\cals_2)).$
\item{(ii)} $\cals_1\star \cals_2\in \cals\cale^s$ if and only if
$\gcd(m_1l_2,m_2l_1)=1.$
\item{(iii)} $\calr$ is a submonoid.
\item{(iv)} If $\cals_1$ is regular (i.e., $m_1=1$) and $\hbox{Ind}(\cals_2)$
divides $\hbox{Ind}(\cals_1),$ then $\cals_1\star \cals_2\in \cals\cale^s$
independently of $m_2.$ \tenrm

\noindent{\sc Proof}: (i) and (iii) are clear from the discussion above, and
(iv) is a special case of (ii), which we now prove. Since $\cals_1\times
\cals_2$ is a 2-torus V-bundle over $\calz_1\times \calz_2,$ the $S^1$ V-bundle
$\cals_1\star \cals_2$ can be realized as a quotient $(\cals_1\times
\cals_2)/S^1$ by choosing a certain homomorphism $S^1\ra{1.3} S^1\times S^1.$
For the V-bundle with Chern class given by \sewz.3 we have the action on
$\cals_1\times \cals_2$ given by $(x,y)\mapsto (\grt^{l_2}x,\grt^{-l_1}y).$ The
condition that $\cals_1\star \cals_2$ is in $\cals\cale^s$ is that there are no
fixed points $(x,y)$ under the above action. If $m_x,m_y$ denote the orders of
the local uniformizing groups (leaf holonomy groups) of $x,y,$ respectively,
then the condition that $(x,y)$ be a fixed point is that
$\gcd(m_xl_2,m_yl_1)=g>1.$ That this condition never hold for all pairs $(x,y)$
is precisely the condition $\gcd(m_1l_2,m_2l_1)=1.$ \hfill\za

This proposition will allow us to construct smooth examples of quasi-regular
\Se manifolds in all odd dimensions greater than or equal to $9.$ Before
discussing some important examples, we introduce some more terminology.

\noindent{\sc Definition} \sewz.7: \tensl We say that $\cals\in \cals\cale$ is
{\it $\cals\cale$-irreducible} if writing $\cals =\cals_1\star \cals_2$ implies
that either $\cals_1$ or $\cals_2$ is $S^1$ with its flat structure. $\cals$ is
{\it $\cals\cale$-reducible} if it is not $\cals\cale$-irreducible. \tenrm

It is clear that $\cals\cale$-irreducibility corresponds to Riemannian
irreducibility on the space of leaves. Reducibility first occurs in dimension
5, and up to diffeomorphism there is precisely one regular
$\cals\cale$-reducible 5-manifold, namely the homogeneous Stiefel manifold
$V_2(\bbr^4)\approx S^2\times S^3\approx S^3\star S^3.$ In terms of the
homogeneous \Se manifolds, reducibility has a group theoretic interpretation. A
well-known theorem of Borel [Bo] states a simply connected homogeneous
K\"ahler manifold corresponding to a complex semi-simple Lie group $G$ is
Riemannian irreducible if and only if $G$ is simple. It is easy to see that in
terms of \Se geometry this implies

\noindent{\sc Proposition} \sewz.8: \tensl A simply connected \Se
homogeneous manifold corresponding to a semi-simple Lie group $G$ is
$\cals\cale$-irreducible if and only if $G$ is simple. \tenrm

Moreover, our construction also implies

\noindent{\sc Proposition} \sewz.9: \tensl Let $\cals_1,\cals_2$ be simply
connected homogeneous \Se manifolds, then $\cals_1\star \cals_2$ is a simply
connected homogeneous \Se manifold. \tenrm

This proposition states that the subset $\calh$ of simply connected homogeneous
\Se manifolds forms a submonoid of $\calr.$

As noted by Theorem \ser.5, the regular \Se 5-manifolds have been classified.
This immediately gives a classification of the regular $\cals\cale$-reducible
7-manifolds.

\noindent{\sc Proposition} \sewz.10: \tensl Any simply connected regular
$\cals\cale$-reducible 7-manifold is one of the following:
$$S^3\star S^3\star S^3, \quad S^3\star S^5, \quad S^3\star \cals_k$$
for $3\leq k\leq 8$ where $\cals_k$ is one of the \Se circle bundles over the
del Pezzo surface discussed in Theorem \ser.5. \tenrm

These examples have already been noted in [BFGK]. Non simply connected examples
are obtained by quotienting by a cyclic subgroup of the circle generated by the
characteristic vector field. Notice also that for $5\leq k\leq 8$ the
7-manifolds $S^3\star S_k$ have continuous families of \Se structures on
them. We shall discuss these shortly, but first we mention that in
the list given in Proposition \sewz.10 only the first two are homogeneous.  In
fact as already noted in [BG2] we have

\noindent{\sc Proposition} \sewz.11: \tensl The simply connected homogeneous
\Se 7-manifolds are precisely one of the following:
\item{(i)} An $\cals\cale$-reducible manifold $S^3\star S^3\star S^3$ or
$S^3\star S^5.$
\item{(ii)} The real Stiefel manifold $V_{5,2}.$
\item{(iii)} A homogeneous 3-Sasakian 7-manifold $S^7$ or $\cals(1,1,1).$
\tenrm

Of course, there are other regular $\cals\cale$-irreducible 7-manifolds as
noted in Examples \ser.8 and \ser.10, but so far a complete classification of
regular \Se 7-manifolds is lacking. Actually the results of Mori and Mukai
[MoMu] on the classification of smooth Fano 3-folds with $b_2\geq 2$ go quite
far toward a classification of regular \Se 7-manifolds, but this is beyond the
scope of the present paper, and is currently under investigation.
Nevertheless, there are some immediate consequences of their work worth
mentioning here, namely

\noindent{\sc Proposition} \sewz.12: \tensl Let $\cals$ be a compact regular
\Se 7-manifold. Then
\item{(i)} $b_2(\cals)\leq 9.$
\item{(ii)} If $b_2(\cals)\geq 5$ then $\cals$ is reducible. Explicitly $\cals
\approx S^3\star \cals_k$, where $4\leq k\leq 8.$ \tenrm

Notice that the upper bound on the second Betti number for regular \Se
7-manifolds is realized by $S^3\star S_8,$ and there is a continuous
4-dimensional family of \Se structures on these manifolds. It is interesting
to contemplate whether it is generally true that the bound on $b_2$ is realized
by a $\cals\cale$-reducible element in $\calr.$ If this were true then the
regular \Se 9-manifold with the largest second Betti number would be $S_8\star
S_8$ with $b_2=17.$ The answer probably lies in Mori theory. For general odd
dimension greater than 3 our join construction gives

\noindent{\sc Proposition} \sewz.13: \tensl For each positive integer $n\geq 2$
there exists smooth manifolds of dimension $2n+1$ that admit continuous
families $\calf$ of regular \Se structures. Examples of such the manifolds are
$$S_k,\quad \overbrace{S^3\star\cdots\star S^3}^{n-2\rm\;times}\star S_k$$
for $5\leq k\leq 8,$ as well joins with other regular \Se manifolds. Moreover,
the dimension of the family $\calf$ of \Se structures on $S_k$ is $\geq k-4.$
\tenrm

This result should be contrasted with the 3-Sasakian case which is
infinitesimally rigid [PP]. Manifolds that admit continuous families of
non-regular \Se structures also exist as will be discussed below. However, they
begin in dimension 11.

Next we turn to the more lucrative non-regular case. Propositions \seg.12 and
\sewz.6 imply the following:

\noindent{\sc Corollary} \sewz.14: \tensl Let $\cals$ be any compact simply
connected 3-Sasakian manifold of dimension $4n+3$ that is not a sphere. Then
for any positive integer $r$ the join $S^{2(n+1)r-1}\star \cals$ is a smooth
\Se manifold of dimension $4n+2nr+2r+1.$ In particular, if $n=1$, so that
$\hbox{dim}~\cals =7,$ then $S^3\star \cals$ is a smooth \Se
9-manifold. \tenrm

\noindent{\sc Corollary} \sewz.15: \tensl Let $\cals$ be any compact simply
connected \Se 7-manifold of index $2,$ for example, a 3-Sasakian 7-manifold
that is not $S^7.$ Then $S^{2m+1}\star \cals$ is smooth if $m$ is odd or if $m$
is even and the order of $\cals$ is odd.  \tenrm

The last case in Corollary \sewz.14 will prove to be of much interest to us.
Notice that $\hbox{Ind}(S^3\star
\cals)=\gcd(\hbox{Ind}(S^3),\hbox{Ind}(\cals))=2,$ so that the procedure
iterates arriving at:

\noindent{\sc Corollary} \sewz.16: \tensl Let $\cals$ be any compact simply
connected 3-Sasakian manifold of dimension $7$ that is not a 7-sphere. Then for
any positive integer $r$ the r-fold join $S^3\star\cdots \star S^3\star\cals$
is a smooth \Se $2r+7$-manifold of index $2.$ \tenrm

The examples of a 3-Sasakian 7-manifold that we have in mind are the toric
3-Sasakian manifolds $\cals(\grO_k)$ of [BGMR] where $\grO_k$ is a $k$ by $k+2$
matrix of integers which satisfy certain gcd conditions (see [BGMR] for
details), and $b_2(\cals(\grO_k))=k.$  In the case of the join of two
non-regular \Se manifolds and/or when the relative indices are different, the
gcd conditions in Proposition \sewz.6 are generally fairly restrictive.  For
example, let $\cals_1,\cals_2\in \cals\cale^s_{4n+3}$ be two 3-Sasakian
manifolds neither of which are spheres.  Then $l_1=l_2=1,$ and the smoothness
conditions become $\gcd(m_1,m_2)=1.$ Even this is restrictive since orders tend
to be large and have many divisors.  However, consider the 3-Sasakian manifolds
$\cals(1,\cdots,1,2p_i+1)$ discussed in [BGM2]. In these cases $m_i=p_i+1,$ so
if we choose $\gcd(p_1+1,p_2+1)=1,$ which is easy to satisfy, we get a smooth
join. Notice, however, that $\cals\star \cals$ is never smooth for $\cals$
non-regular.

Finally we consider some non-regular examples of manifolds that admit
continuous families of \Se structures. Then Proposition \sewz.6 implies

\noindent{\sc Proposition} \sewz.17: \tensl Let $\cals$ be a simply connected
\Se manifold, and $S_k$ the circle bundle over the del Pezzo surface
defined in Theorem \ser.5. Then
\item{(i)} $S_k\star \cals$ admits a continuous family of dimension $\geq k-4$
of inequivalent non-regular \Se structures when $5\leq k\leq 8.$
\item{(ii)} $S_k\star \cals$ is a smooth manifold if and only if
$\gcd(\hbox{Ind}(\cals),\hbox{Ord}(\cals))=1.$ In particular, $S_k\star
\cals(\grO_k)$ is smooth if and only if $\hbox{Ord}(\cals(\grO_k))$ is odd.
\tenrm

We have only computed the orders of the 3-Sasakian manifolds $\cals(\grO_k)$ in
the case $k=1.$ These are the 3-Sasakian 7-manifolds $\cals(p_1,p_2,p_3)$ of
[BGM2] where the $p_i$'s are pairwise relatively prime.  Furthermore, the
order is the least common multiple of
$({p_1+p_2\over 2})({p_+p_3\over
2})({p_2+p_3\over 2})$ if all $p_i$'s are odd, and the least common multiple of
the product
$(p_1+p_2)(p_1+p_3)(p_2+p_3)$ if one of the $p_i$'s is even.  In the last case
one of the sums is even so this can be eliminated. In the first case we find a
solution if $p_i=4r_i+1$ for some natural numbers $r_i.$ Summarizing we have

\noindent{\sc Corollary} \sewz.18: \tensl The 11 dimensional orbifolds
$S_k\star \cals(p_1,p_2,p_3)$ are smooth manifolds if and only if $p_i=4r_i+1$
for some natural numbers $r_i$ which satisfy $\gcd(4r_i+1,4r_j+1)=1.$
Furthermore, if $5\leq k\leq 8$ the manifolds $S_k\star \cals(p_1,p_2,p_3)$
admit a continuous family of dimension $\geq k-4$ of inequivalent
non-regular \Se structures. In particular, for each natural
number $r$ and for $5\leq k\leq 8$ the 11-manifolds $S_k\star
\cals(1,1,4r+1)$ admit such continuous families of non-regular \Se structures.
\tenrm

The gcd condition in Proposition \sewz.17 seemingly becomes less restrictive
when the dimension of $\cals$ is such that $n+1$ is a large prime. However, it
is in dimension 7 that we have the most interesting examples of 3-Sasakian
manifolds, those having arbitrary second Betti number constructed in [BGMR].
However, as previously mentioned it is only in the case $b_2=k=1$ that we
have computed the order.  It would be interesting to see whether there
exist continuous families of \Se structures on manifolds with any second Betti
number.

Notice that once we find solutions $\cals$ to the conditions in Proposition
\sewz.16, we can easily construct solutions in higher dimension, for example by
joining with $S^3.$ Since $\hbox{Ind}(\cals_k\star \cals)=1$ the smoothness
conditions for $S^3\star \cals_k\star \cals$ are automatically satisfied. Thus,

\noindent{\sc Corollary} \sewz.19: \tensl Let $\cals$ be a simply connected
\Se manifold of dimension $2n+1$ satisfying
$\gcd(\hbox{Ind}(\cals),\hbox{Ord}(\cals))=1,$ then
$\overbrace{S^3\star\cdots\star S^3}^{m\rm\;times} \star S_k\star \cals$ is a
smooth $2(m+n)+5$-manifold that admits a continuous family of dimension $\geq
k-4$ of \Se structures. In particular, the $(2m+11)$-manifolds
$\overbrace{S^3\star\cdots\star S^3}^{m\rm\;times} \star S_l\star
\cals(\grO_k)$  with $\hbox{Ord}(\cals(\grO_k))$ odd has second Betti number
$b_2=k$ and continuous families of inequivalent \Se structures. \tenrm

We should mention that we do not know whether the condition
$\hbox{Ord}(\cals(\grO_k))$ be odd can be satisfied for arbitrary $k.$

\bigskip
\centerline {\bf \S 5. The Cohomology of Some Joins}
\medskip

We first obtain some general information about the low Betti numbers of the
join $\cals_1\star \cals_2$ of two \Se manifolds (orbifolds). The following
lemma follows from harmonic theory, and was given in [BG1]:

\noindent{\sc Lemma} \seex.1: \tensl Let $\cals$ be a \Se orbifold of dimension
$2n+1,$ then for $0\leq r\leq n$ we have
$b_r(\cals)=b_r(\calz)-b_{r-2}(\calz).$ \tenrm

This lemma can be used to show:

\noindent{\sc Lemma} \seex.2: \tensl Let $\cals_i\in \cals\cale_{2n_i+1},$ then
\item{(1)} $b_2(\cals_1\star \cals_2)=b_2(\cals_1)+b_2(\cals_2)+1$ if $n_i\geq
1,$
\item{(2)} $b_3(\cals_1\star \cals_2)=b_3(\cals_1)+b_3(\cals_2)$ if $n_i\geq
3,$
\item{(3)} $b_4(\cals_1\star \cals_2)=b_4(\cals_1)+b_4(\cals_2)+b_2(\cals_1)
b_2(\cals_2)+b_2(\cals_1)+b_2(\cals_2)+1$ if $n_i\geq 4.$ \tenrm

When $n_i$ is outside the indicated range, the formula is slightly different,
but is easily worked out. For general $b_r$ the formulas are increasingly more
complicated and are different depending on whether $r$ is
even or odd, or whether the range conditions are satisfied or not. In order to
determine the cohomology of $\cals_1\star \cals_2$ in specific cases,
we shall employ a more elegant technique using spectral sequences used by Wang
and Ziller [WZ].

Let $\cals_i\in \cals\cale^s,$ and consider the commutative diagram of
fibrations
$$\matrix{\cals_1\times \cals_2&\ra{2.6} &B(\cals_1\star \cals_2)&\ra{2.6}
&BS^1\cr
\decdnar{=}&&\decdnar{}&&\decdnar{\psi}\cr
\cals_1\times \cals_2&\ra{2.6} &B\calz_1\times B\calz_2&\ra{2.6}
&BS^1\times BS^1.\cr}\leqno{\seex.3}$$
The maps are all the obvious ones. In particular, $\psi$ is determined by the
$S^1$ action of the previous section, namely, $\psi(\grt)=
(\grt^{l_2},\grt^{-l_1}).$ The point is that the differentials in the Serre
spectral sequence of the top fibration are determined through naturality by the
differentials in the Serre spectral sequence of the bottom fibration. Wang and
Ziller apply this method to computing the integral cohomology ring of more
general circle bundles (torus bundles as well) over products of projective
spaces. For us, this corresponds to the case $S^{2m+1}\star S^{2n+1},$ which is
homogeneous. We refer to [WZ] for the cohomology ring in this case. At first we
shall apply this method to a more general situation where only rational
information can be obtained; however, there are several cases of interest
to us where we have enough information about the differentials to compute the
integral cohomology groups. Recall that from its definition the join
$\cals_1\star \cals_2$ of two simply connected \Se manifolds is necessarily
simply connected (in the orbifold sense if $\cals_1\star \cals_2$ is not
smooth). Nevertheless, we can easily obtain non-simply connected \Se manifolds
with cyclic fundamental group by dividing by a cyclic subgroup of the circle
generated by the characteristic vector field $\xi.$ Hereafter, in this section
all joins are simply connected.

\noindent{\sc Theorem} \seex.4: \tensl Let $\cals$ be any simply connected
3-Sasakian 7-manifold which is not the 7-sphere. Then
$$H^q(S^{2m+1}\star \cals,\bbq)\approx
\cases{H^q(S^{2m+1}\times \calz,\bbq) & if $m>2$;\cr H^q(S^2\times
\cals,\bbq)& if $m=1$;\cr \cases{\bbq &if $q=0,11$;\cr
       \bbq^{k+1} &if $q=2,4,7,9$;\cr
       0 &if $q=1,10,3,8$;} &if $m=2,$\cr}$$
where $k=b_2(\cals).$ Moreover, for the case when
$m=2$ if $k>1$ then $H^5(S^5\star
\cals,\bbq) \approx H^6(S^5\star \cals,\bbq)\approx 0;$ whereas if $k=1$ there
are two possibilities; either $H^5(S^5\star \cals,\bbq) \approx H^6(S^5\star
\cals,\bbq)\approx 0,$ or $H^5(S^5\star \cals,\bbq) \approx H^6(S^5\star
\cals,\bbq)\approx \bbq.$ \tenrm

Notice in this case that the rational cohomology of $S^{2m+1}\star \cals$
depends only on $b_2(\cals).$ Moreover, when $m=2$ and $k>1$ it does not have
the rational cohomology of a product. Even in other cases $S^{2m+1}\star \cals$
cannot be a product. For example, if $m$ is odd and greater than $1,$ then by
Corollary \sewz.8 $S^{2m+1}\star \cals$ is smooth, but if $\cals$ is
non-regular, $S^{2m+1}\times \calz$ is not. However, in the second possibility
for $m=2$ and $k=1$ we see that $H^q(S^5\star \cals,\bbq)\approx
H^q(\bbc\bbp^2\times \cals,\bbq).$

\noindent{\sc Proof of Theorem \seex.4}: We consider the Serre spectral
sequence of the fibrations in diagram \seex.3. We also make use of the fact
(Corollary \secl.17) that rationally $H^q$ and $H^q_{orb}$ coincide. If $m>3$
then the orientation class $u$ of $S^{2m+1}$ in the fiber $S^{2m+1}\times
\cals$ occurs after the orientation class of $\cals,$ so by commutativity of
diagram \seex.3, $E^{p,q}_2$ of the Serre spectral sequence for the top
fibration in \seex.3 coincides rationally with the $E_2$ term for the spectral
sequence for the fibration $\cals\ra{1.3} B\calz\ra{1.3} BS^1$ for $q<2m+1.$ It
follows that the $E_2$ term of the former converges to the cohomology of the
product $S^{2m+1}\times \calz.$ For $m=3$ there are two 7-classes in the
cohomology of the fiber, say $u$ and $v.$ By naturality of the spectral
sequences, we have that $d_8(u)=l_2^4s^4$ and $d_8(v)=l_1^4s^4$ where $l_i$ are
given by \sewz.5. The same argument goes through as before, but now it is the
7-class $l_1^4u-l_2^4v$ that survives to the limit. So again rationally the
cohomology is the cohomology of the product $S^7\times \calz.$

Next consider the case $S^3\star \cals$ for which $(l_1,l_2)=(1,1).$ Here we
discuss this case integrally with the added assumption that
$H^3(\cals,\bbz)=0,$ since it will be treated in Theorem \seex.6 below.  By
the Leray-Serre Theorem  the $E_2$ term of the spectral sequence for the fibration
$S^3\times\cals\ra{1.3} S^3\star\cals\ra{1.3} BS^1$ is given by

$$
\beginpicture
\setcoordinatesystem units <2pt,1.5pt>
\setlinear
\put {\hbox to .2pt{\hfill}} [t] at 20 65

\plot 0 96  0 0  116 0 /

\multiput {$\bullet$} at 0 0  0 18  0 36  0 54  0 72  0 90 /
\multiput {$\bullet$} at 18 0  18 18 18 36  18 54  18 72 18 90 /
\multiput {$\bullet$} at 36 0  36 18 36 36  36 54  36 72 36 90 /
\multiput {$\bullet$} at 54 0  54 18 54 36  54 54  54 72 54 90 /
\multiput {$\bullet$} at 72 0  72 18 72 36  72 54  72 72  /
\multiput {$\bullet$} at 90 0  90 18 90 36  90 54  /
\multiput {$\bullet$} at 108 0  108 18 108 36  /

\multiput {$\bbz^{2k}$} at 6 94 40 94 /
\multiput {$\bbz^k$} at  4 40 40 40 76 40 111 40 /
\multiput {$\bbz$} at  3 3 40 3 76 3 112 3 3 57 40 57 76 57 /
\multiput {$T$} at  4 76 40 76 76 76 /
\multiput {$0$} at 4 20 20 20  38 20  56 20 74 20 92 20  20 3 20 38 20 56 20 74
20 92 56 3 56 20 56 38 56 56 56 74 56 92 92 3 92 38 92 56 111 20 /

\put {$p$} at 124 0
\put {$q$} at  0 110
\put {\hbox{Diagram \seex.5}: $E^{p,q}_2$ for $S^3\times\cals\ra{1.3} S^3\star
\cals\ra{1.3} BS^1$} [t] at  45 -15
\endpicture $$
Again let $u$ denote the orientation class of $S^3$ and $s$ the 2-class of
$BS^1,$ then by naturality of the diagram \seex.3, we have $d_4(u)=s^2.$ Now
$S^3\star \cals$ is a 9-manifold, so by Poincar\'e duality it suffices to
consider the diagram up to and including dimension 4. But it is easily seen
from Diagram \seex.5 that no other differentials can occur in this range.
Furthermore, the torsion groups $T$ in row 4 do not occur rationally, and the
result follows.

Finally we consider the case $m=2.$ Over $\bbq$ the $E_2$
term of the spectral sequence in this case is obtained from
Diagram \seex.5 by tensoring with $\bbq,$ putting $0$'s
in the fourth row ({\it i.e.,} $E^{p,3}_2=0$), and changing the $\bbq^{2k}$'s
in row six to $\bbq^{k+1}.$ Now
there are no differentials below level 5, so by Poincar\'e duality we need only
concern ourselves with the $E^{0,5}_2$ term. Here there are $k+1$ classes,
namely the orientation class $u$ of $S^5$ and the 5-classes $\grb_i\in
H^5(\cals,\bbq)$ that are the Poincar\'e duals of the 2-classes in
$H^2(\cals,\bbq).$ By looking at the rational spectral sequence of the fibration
$\cals\ra{1.3} B\calz\ra{1.3} BS^1,$ one sees that there are two possibilities
(details will be given elsewhere).  The first possibility is that
$d_4(\grb_i)=\sum_ja_{ij} s^2\otimes \gra_j$ where the rank of the matrix
$(a_{ij})$ is $k.$ In this case we also have $d_6(u)=s^3.$ This implies that
$H^5(S^5\star \cals,\bbq) \approx H^6(S^5\star\cals,\bbq)\approx 0.$ The
second possibility only occurs if $k=1,$ and in this case we get $H^5(S^5\star
\cals,\bbq) \approx H^6(S^5\star \cals,\bbq)\approx \bbq.$ \hfill\za

Generally, this only determines the rational cohomology; however, for
$S^3\star \cals,$ since $(l_1,l_2)=(1,1)$ the only torsion class in dimension 4
occurs in the $\cals$ factor on the fiber. This class will survive to
$E_\infty.$ Thus, it remains to solve an extension problem to determine
the integral cohomology in this case. Below dimension 4 integral information
follows easily from Diagram \seex.5; however, since determination of the full
integral cohomology requires some more detailed knowledge we shall specialize
to the case of the 3-Sasakian manifolds $\cals =\cals(p_1,p_2,p_3)$ described
in [BGM2] for the complete picture.   Actually we do not directly solve the
extension problem, but rather take a different tact.

\noindent{\sc Theorem} \seex.6: \tensl Let $\cals$ be a simply connected \Se
7-manifold of index $2$ with $H^3(\cals,\bbz)\approx 0$ and second Betti number
$b_2(\cals)=k.$ Then we have
$H^2(S^3\star \cals,\bbz)\approx \bbz^{k+1}$ and $H^3(S^3\star
\cals,\bbz)\approx 0.$
More specifically, if $\cals(p_1,p_2,p_3)$ is one of
the simply connected 3-Sasakian 7-manifolds described in [BGM2] with the
$p_i$'s pairwise relatively prime.  Then there is an isomorphism
$$H^q(S^3\star \cals(p_1,p_2,p_3),\bbz)\approx \cases{\bbz & if $q=0,9,5$;\cr
                    \bbz^2 & if $q=2,7$;\cr
                    \bbz\oplus \bbz_{\grs_2} & if $q=4;$\cr
                    \bbz_{\grs_2} & if $q=6;$\cr
                    0 & otherwise,\cr}$$
where $\grs_2 =p_1p_2+p_1p_3+p_2p_3.$ In fact there is a ring isomorphism
$H^*(S^3\star \cals(p_1,p_2,p_3),\bbz)\approx H^*(S^2\times
\cals(p_1,p_2,p_3),\bbz).$ \tenrm

\noindent{\sc Proof}: As mentioned above the first statement follows easily
from Diagram \seex.5. However, the proof of the remainder is much more
involved.  By Theorem \seex.4 and Poincar\'e duality it suffices to show
$H^4(S^3\star \cals(p_1,p_2,p_3),\bbz)\approx \bbz\oplus \bbz_{\grs_2}.$ Our
approach is to consider $S^3\times SU(3)$ as a 2-torus bundle over $S^3\star
\cals(p_1,p_2,p_3).$ Consider the action of $S^1\times S^1$ on $S^3\times
SU(3)$ given by
$$(q,\bba)\longmapsto (q\grr,\hbox{diag}(\grt^{p_1},\grt^{p_2},\grt^{p_3})\bba~
\hbox{diag}(\grr,\grr^{-1},\grt^{-p_1-p_2-p_3})).\leqno{\seex.7}$$
It is well-known that the cohomology ring of $SU(3)$ is isomorphic to
the exterior algebra over 2 generators $e_3,e_5$ in cohomological dimension $3$
and $5,$ respectively. So the $E_2$ term of the Leray-Serre spectral sequence
for the fibration
$$S^3\times SU(3)\ra{1.8} S^3\star \cals(p_1,p_2,p_3)\ra{1.8} BS^1\times BS^1
\leqno{\seex.8}$$
is given by
$$
\beginpicture
\setcoordinatesystem units <2pt,1.5pt>
\setlinear
\put {\hbox to .2pt{\hfill}} [t] at 20 65

\plot 0 96  0 0  116 0 /

\multiput {$\bullet$} at 0 0  0 18  0 36  0 54  0 72  0 90 /
\multiput {$\bullet$} at 18 0  18 18 18 36  18 54  18 72 18 90 /
\multiput {$\bullet$} at 36 0  36 18 36 36  36 54  36 72 36 90 /
\multiput {$\bullet$} at 54 0  54 18 54 36  54 54  54 72 54 90 /
\multiput {$\bullet$} at 72 0  72 18 72 36  72 54  72 72  /
\multiput {$\bullet$} at 90 0  90 18 90 36  90 54  /
\multiput {$\bullet$} at 108 0  108 18 108 36  /

\multiput {$\bbz$} at 3 3 3 94 38 94 /
\multiput {$0$} at  3 38 38 38 74 38 111 38 3 74 38 74 74 74 /
\multiput {$\bbz^3$} at  76 3 /
\multiput {$\bbz^4$} at  40 57 112 3 /
\multiput {$\bbz^6$} at  76 57 /
\multiput {$\bbz^2$} at  3 57 40 3  /
\multiput {$0$} at 3 20 20 20  38 20  56 20 74 20 92 20  20 3 20 38 20 56 20 74
20 92 56 3 56 20 56 38 56 56 56 74 56 92 92 3 92 38 92 56 111 20 /

\put {$p$} at 124 0
\put {$q$} at  0 110
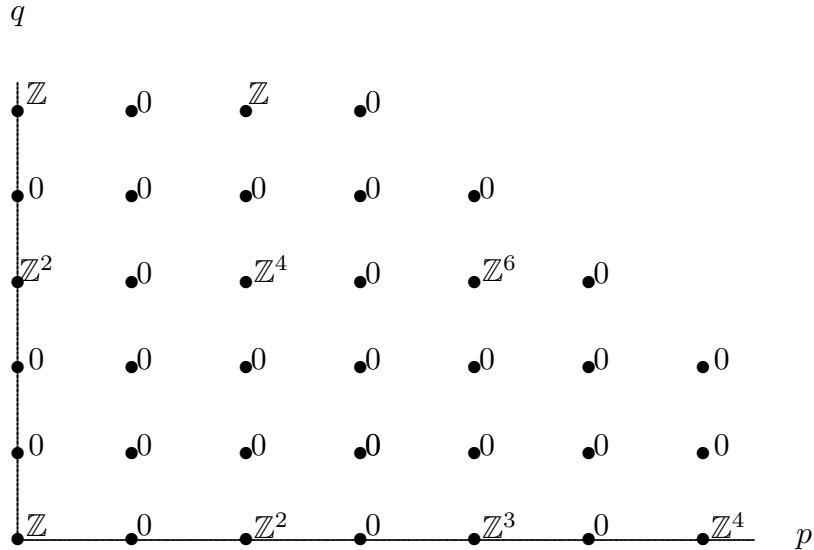
\put {\hbox{Diagram \seex.9}: $E^{p,q}_2$ for $S^3\times SU(3)\ra{1.3}
S^3\star \cals(p_1,p_2,p_3)\ra{1.3} BS^1\times BS^1$} [t] at 45 -15 \endpicture
$$
Let $u,e_3$ denote the 3-classes in $E^{0,3}_2.$ We need to determine the
differentials $d_4(u)$ and $d_4(e_3).$ The theorem follows easily from

\noindent{\sc Lemma} \seex.10: \tensl The differentials satisfy
\item{(i)} $d_4(u)=t^2,$
\item{(ii)} $d_4(e_3)=\grs_2s^2-t^2,$
where $t$ and $s$ are the corresponding 2-classes in $BS^1\times BS^1.$
\tenrm

\noindent{\sc Proof}: First we notice from the form of the action \seex.8 that
the fibration \seex.8 factors through the fibration
$$S^3\times SU(3)\ra{1.8} S^3\times \cals(p_1,p_2,p_3)\ra{1.8}
BS^1,\leqno{\seex.11}$$
where the $S^1$ action is defined by the subgroup obtained by putting
$\grr=\hbox{id}$ in \seex.7. This $S^1$ acts only on the $SU(3)$ factor, and
the $d_4(e_3)$ was determined in [BGM2] to be $d_4(e_3)=\pm\grs_2s^2.$ Thus, by
naturality for the fibration \seex.8 we must have the form
$d_4(e_3)=\grs_2s^2+at^2+bst$ where $a$ and $b$ are to be determined.
Actually, since the fibration \seex.8 factors through both the fibration
\seex.11 and the fibration
$$S^3\times SU(3)\ra{1.8} (S^3\times SU(3))/S^1(\grr)
\ra{1.8} BS^1,\leqno{\seex.12}$$
where the action $S^1(\grr)$ is obtained from \seex.8 by putting
$\grt=\hbox{id},$ it follows that $b=0.$ Now the action $S^1(\grr)$ is a
subgroup of a 2-torus action which acts on the $S^3$ and $SU(3)$ factors
separately. So we have a commutative diagram of fibrations
$$\matrix{S^3\times SU(3)&\ra{2.6} &(S^3\times SU(3))/S^1(\grr)&\ra{2.6}
&BS^1\cr
\decdnar{=}&&\decdnar{}&&\decdnar{}\cr
S^3\times SU(3)&\ra{2.6} &S^2\times M_{1,-1}&\ra{2.6}
&BS^1\times BS^1,\cr}\leqno{\seex.13}$$
where $M_{1,-1}$ is a certain Aloff-Wallach manifold described in Eschenburg
[Esch1-2].
Now for the bottom row we have the Hopf fibration on the first factor,
so $d_4(u)=t_1^2,$ and Eschenburg shows that $d_4(e_3)=t_2^2$ for the second
factor. By naturality these pull back to $d_4(u)=t^2$ and $d_4(e_3)=t^2$ for
the fibration on the top row of \seex.13. Now combining \seex.13 with \seex.9
we have the commutative diagram of fibrations
$$\matrix{S^3\times SU(3)&\ra{2.6} &S^3\star \cals(p_1,p_2,p_3)&\ra{2.6}
&BS^1\times BS^1\cr
\decdnar{=}&&\decdnar{}&&\decdnar{\psi}\cr
S^3\times SU(3)&\ra{2.6} &S^2\times \calz(p_1,p_2,p_3)&\ra{2.6}
&BS^1\times BS^1\times BS^1.\cr}\leqno{\seex.14}$$
Here the map $\psi$ is induced by $(\grr,\grt)\mapsto (\grr,\grr^{-1},\grt).$
Pulling back by this map proves the lemma, and hence, Theorem \seex.6.
\hfill\za

Now let $\cals$ be a simply connected \Se 7-manifold of index $2$ with
vanishing $H^3(\cals,\bbz),$  and consider the $r$-fold iterates $S^3\star
\cdots\star S^3\star \cals$ which by Proposition \sewz.6 are all have smooth of
index $2.$ For general $\cals$ the argument is with rational coefficients, but
the argument will work integrally for any such $\cals$ that satisfies
$H^*(S^3\star \cals,\bbz)\approx H^*(S^2\times \cals,\bbz).$ So assume by
induction that in the fibration
$$S^3\times (S^3\star\cdots \star S^3\star
\cals) \ra{1.6} S^3\star\cdots \star S^3\star \cals \ra{1.6} BS^1$$
the fiber is homologically $S^3\times S^2\times\cdots S^2\times \cals.$ Then in
the spectral sequence there are no differentials off of any of the $S^2$
factors, so all the differentials come from the $S^3\times \cals$ term as in
Theorems \seex.4 and \seex.7. It follows that the total space is homologically
$S^2\times\cdots S^2\times \cals,$ and we have arrived at

\noindent{\sc Theorem} \seex.15: \tensl Let $\cals$ be a simply connected
\Se 7-manifold of index $2$ with $H^3(\cals,\bbz)=0.$ Then there are group
isomorphisms
$$H^2(S^3\star\cdots \star S^3\star \cals,\bbz)\approx H^2(S^2\times \cdots
\times S^2\times \cals,\bbz),\quad H^3(S^3\star\cdots \star S^3\star
\cals,\bbq)\approx 0,$$
and ring isomorphisms
$$H^*(S^3\star\cdots \star S^3\star \cals,\bbq)\approx H^*(S^2\times \cdots
\times S^2\times \cals,\bbq),$$
and
$$H^*(S^3\star\cdots \star S^3\star \cals(p_1,p_2,p_3),\bbz)\approx
H^*(S^2\times \cdots \times S^2\times \cals(p_1,p_2,p_3),\bbz).$$
\tenrm

It is interesting to speculate as to whether $S^3\star\cdots \star S^3\star
\cals$ is actually homeomorphic (or even diffeomorphic) to
$S^2\times\cdots\times S^2\times \cals$ in this theorem.   It is not true for
the homogeneous space $S^3\star S^7,$ as it was shown by Wang and Ziller [WZ],
where in their notation $S^3\star S^7$ is $M^{1,3}_{1,2},$  that it is a
nontrivial $\bbr\bbp^7$ bundle over $S^2.$ But this has more to do
with the fact that the relative indices for $S^3\star S^7$ are $(1,2)$ instead
of $(1,1).$ Indeed, Wang and Ziller show that for the circle action
corresponding to $(1,1)$ their $M^{1,3}_{1,1}$ is diffeomorphic to $S^2\times
S^7.$ However, their proof does not seem to generalize to our more general
(replacing $S^7$ by $\cals$) situation. Even if $S^3\star\cdots \star S^3\star
\cals$ were diffeomorphic to $S^2\times\cdots\times S^2\times \cals$ the \Se
metric could not be the product metric.

In [BGMR] the authors constructed ``toric'' 3-Sasakian 7-manifolds
$\cals(\grO_k)$ with arbitrary second Betti number. As mentioned previously
these manifolds depend on a certain $k$ by $k+2$ matrix $\grO_k.$
Moreover, it was shown in [BGM6] that these manifolds are simply connected and
that $H^3(\cals(\grO_k),\bbz)=0.$ Thus the results in [BGMR] and [BGM6]
together with Theorem \seex.7 imply

\noindent{\sc Corollary} \seex.16: \tensl There exist compact simply connected
\Se manifolds with arbitrary second Betti number $b_2$ and $b_3=0$ in every odd
dimension greater than $5,$ namely the manifolds $S^3\star\cdots \star S^3\star
\cals(\grO_k).$ \tenrm

It is still an open question whether this corollary holds in dimension 5 as
well. For the 3-Sasakian 7-manifolds $\cals(p_1,p_2,p_3)$ of [BGM2], Theorem
\seex.7 implies

\noindent{\sc Corollary} \seex.17: \tensl Let $p_1,p_2,p_3$ be pairwise
relatively prime positive integers. Then among the \Se manifolds
$S^3\star\cdots \star S^3\star \cals(p_1,p_2,p_3)$ there are infinitely many
that are homotopically distinct in every odd dimension greater than $5.$
Furthermore, all these manifolds have $b_2=1$ and $b_3=0.$ \tenrm

Next we consider \Se manifolds $S^3\star \cals_{3,4}$ and $\cals_{3,4}\star
\cals$ where again $\cals$ is a simply connected \Se 7-manifold of index 2 with
$H^3(\cals,\bbz)=0,$ and $\cals_{3,4}$ is the \Se circle bundle over the cubic
Fermat surface $F_{3,4}.$ As mentioned in Example \ser.8, $\cals_{3,4}$ has
index 2, so again we get $(l_1,l_2)=(1,1).$ From the spectral sequence for the
fibration $\cals_{3,4}\ra{1.3} F_{3,4}\ra{1.3} BS^1$ we see that the
differentials must satisfy $d_2(u_a)=d_4(u_a)=0$ and $d_2(v_a)=t\otimes u_a$
where $u_a$ are the 3-classes in $\cals_{3,4},$ and $v_a$ are their Poincare
duals. We can then analyze the spectral sequences as before using naturality of
the diagram \seex.3. We find

\noindent{\sc Proposition} \seex.19: \tensl Let $\cals$ be a simply connected
\Se 7-manifold of index $2$ with $H^3(\cals,\bbz)=0$ and $b_2(\cals)=k.$ Let
$\cals_{3,4}$ denote the \Se circle bundle over the cubic Fermat surface
$F_{3,4}.$  Then the \Se 9-manifolds $S^3 \star S_{3,4}$ satisfy
$$H^q(S^3\star \cals_{3,4},\bbz)\approx \cases{\bbz & if $q=0,9,2,7$;\cr
                    \bbz^{10} & if $q=3,6$;\cr
                     0 & otherwise,\cr}$$
whereas for the \Se 13-manifolds $S_{3,4}\star \cals$ satisfy
$$H^2(S_{3,4}\star \cals,\bbz)\approx \bbz^{k+1}, \qquad H^3(S_{3,4}\star
\cals,\bbz)\approx \bbz^{10}.$$ \tenrm

We can now consider the iterated join of $S^3$ with these manifolds, beginning
with $S^3\star\cdots\star S^3\star \cals_{3,4}.$ Notice that $S^3\star
\cals_{3,4}$ is not cohomologically $S^2\times \cals_{3,4}.$ In fact, they are
the same through dimension 3, but differ at 4. However, looking at the spectral
sequence for the fibration
$$S^3\times S^3\star \cals_{3,4}\ra{1.8} S^3\star S^3\star \cals_{3,4}\ra{1.8}
BS^1,$$
we see that indeed $S^3\star S^3\star \cals_{3,4}$ is cohomologically
$S^2\times S^3\star \cals_{3,4}.$ It is now easy to see that the induction
procedure now works in this case. In the case of $S^3\star\cdots\star S^3\star
\cals_{3,4}\star \cals$ where as usual $\cals$ is any simply connected \Se
7-manifold of index 2 with vanishing $H^3(\cals,\bbz),$ we notice that the
iteration argument works through cohomology dimension 3, and we arrive at:

\noindent{\sc Theorem} \seex.20: \tensl For the \Se $2n+7$-manifolds
$\overbrace{S^3\star\cdots\star S^3}^{n\rm\;times}\star \cals_{3,4}$ there is a
ring isomorphism
$$H^*(\overbrace{S^3\star\cdots\star
S^3}^{n\rm\;times}\star\cals_{3,4},\bbz)\approx
H^*(\overbrace{S^2\times\cdots\times S^2}^{n-1\rm\;times} \times
S^3\star\cals_{3,4},\bbz).$$
Let $\cals$ be a simply connected
\Se 7-manifold of index $2$ with $H^3(\cals,\bbz)=0.$ Then for $q=1,2,3$ there
are group isomorphisms
$$H^q(\overbrace{S^3\star\cdots\star S^3}^{n\rm\;times}\star\cals_{3,4}\star
\cals,\bbz)\approx H^q(\overbrace{S^2\times\cdots\times
S^2}^{n-1\rm\;times}\times S^3\star\cals_{3,4}\star \cals,\bbz).$$ \tenrm

Proposition \seex.19 and Theorem \seex.20 have the following corollary:

\noindent{\sc Corollary} \seex.21: \tensl There exist compact simply connected
\Se manifolds with arbitrary second Betti number and non-vanishing third Betti
number in every odd dimension greater than $11,$ namely, the manifolds
$S^3\star\cdots S^3\star\cals_{3,4}\star \cals(\grO_k).$ In particular, for
every integer $n>2$ there exist compact simply connected \Se manifolds with
arbitrary second Betti number of dimension $4n+3$ which cannot admit a
3-Sasakian structure.  \tenrm

It is interesting to ask whether in the case of $k=b_2(\cals)=1$ we can find
infinitely many homotopically distinct manifolds with $b_3\neq 0.$ One can
certainly imitate the proof of Theorem \seex.7 with $S^3$ replaced by the
7-manifold $S_{3,4}.$ However, in the present case there is only one 3-class in
$\cals_{3,4}\times SU(3)$ that does not survive to $E_\infty,$ namely $e_3$ in
$SU(3),$ and this transgresses to a generator in $E^{0,4}_4.$ Thus,
there is no torsion produced at this stage as in Theorem \seex.7.

Finally we consider the case of certain manifolds which admit continuous
families of \Se structures. In general it is difficult to determine the
necessary differentials in the relevant spectral sequences, but for the case of
the regular 7-manifolds $S^3\star \cals_k$ it is quite tractable. Of course as
discussed previously it is only for the range $3\leq k\leq 8$ that we have a
\Se structure, and only in the range $5\leq k\leq 8$, where there are continuous
families.

\noindent{\sc Theorem} \seex.22: \tensl The integral cohomology ring of the
7-manifolds $S^3\star \cals_k$ is given by
$$H^q(S^3\star \cals_k,\bbz)\approx \cases{\bbz & if $q=0,7$;\cr
                    \bbz^{k+1} & if $q=2,5$;\cr
                    \bbz^k_2 & if $q=4;$\cr
                     0 & otherwise,\cr}$$
with the ring relations determined by $\gra_i\cup \gra_j=0, s^2=0, 2\gra_i\cup
s=0,$ where $\gra_i,s$ are the $k+1$ two classes with $i=1,\cdots k.$ \tenrm

\noindent{\sc Proof}: A result of Smale [Sm, BFGK] says that $\cals_k$ is
diffeomorphic to $\#_k(S^2\times S^3).$ So the $E_2$ term of the spectral
sequence for the fibration $S^3\times \cals_k\ra{1.3} S^3\star \cals_k\ra{1.3}
BS^1$ is that of Diagram \seex.5 with the $\bbz$ in row four replaced
by $\bbz^{k+1}$ and the $\bbz^{2k}$ in row six replaced by $\bbz^{k+1}.$
By the Kobayashi-Ochiai Theorem [KO] we have $\hbox{Ind}(\cals_k)=1.$ Thus, we
have $(l_1,l_2)=(2,1).$ As usual we can compute the relevant differentials by
naturality using Diagram \seex.3. Letting $u$ denote orientation class of
$S^3,$ and $\grb_i$ the 3-classes in $\cals_k,$ we find $d_4(u)=s^2,$ and
$d_2(\grb_i)=-2s\otimes \gra_i.$ This gives rise to torsion classes
$E_3^{2,2}\approx \bbz_2^k$ which survive to $E_\infty.$ The remainder now
follows by naturality and Poincar\'e duality. \hfill\za

It would be interesting to construct infinite sequences of homotopically
distinct manifolds of a given dimension that admit continuous families of \Se
structure. For example, the 11-manifolds $\cals_k\star \cals(1,1,4r+1)$ are
good candidates. However, at this time we do not appear to have enough
knowledge of the differentials in the spectral sequence to determine whether
the torsion classes $H^4(\cals(1,1,4r+1),\bbz)\approx \bbz_{8r+3}$ give arise
to enough torsion in the total space to homotopically distinguish these
manifolds.

Of course, we can construct higher dimensional examples with continuous
families of \Se structures, for example by joining with more copies of $S^3.$
However, since the relative indices in this case are $(l_1,l_2)=(2,1)$ more
2-torsion is produced with each iteration making the spectral sequence more
difficult to analyze.  Particular cases where this does not happen are with
the 9-manifolds $S_k\star S_{k'}$ and the 11-manifolds $S_k\star \cals_{4,4}$
and further iterates, where we recall the \Se circle bundle $\cals_{4,4}$ over
the quartic Fermat surface.  These all have relative indices $(1,1)$ and will
have no torsion. For example our methods give

\noindent{\sc Theorem} \seex.23: \tensl The integral cohomology ring of the
11-manifolds $\cals_k\star \cals_{4,4}$ is given by
$$H^q(S_k\star \cals_{4,4},\bbz)\approx \cases{\bbz & if $q=0,11,4,7$;\cr
                     \bbz^{k+1} & if $q=2,9$;\cr
                     \bbz^{60} & if $q=3,8;$\cr
                     \bbz^{60k} & if $q=5,6;$\cr
                     0 & otherwise,\cr}$$
with the ring relations determined by $\gra_i\cup \gra_j=0, s^3=0, \gra_i\cup
s=0,u_a\cup u_b=0, s\cup u_a=0,$ where $\gra_i,s$ are the $k+1$ 2 classes
with $i=1,\cdots k,$ and $u_a$ are the 3 classes with $a=1,\cdots 60.$ \tenrm

\noindent and

\noindent{\sc Theorem} \seex.24: \tensl The integral cohomology ring of the
9-manifolds $\cals_k\star \cals_{k'}$ is given by
$$H^q(S_k\star S_{k'},\bbz)\approx \cases{\bbz & if $q=0,9$;\cr
                     \bbz^{k+k'+1} & if $q=2,7$;\cr
                     \bbz^{kk'+1} & if $q=4,5;$\cr
                     0 & otherwise,\cr}$$
with the ring relations determined by $\gra_i\cup \gra_j=0, s^3=0, \gra_i\cup
s=0,\gra'_i\cup s=0,\gra'_i\cup \gra'_j=0$ where $\gra_i,\gra'_j,s$ are the
$k+k'+1$ 2 classes with $i=1,\cdots k;j=1\cdots k'.$ \tenrm

It is straightforward to work out further iterations all of which have
continuous families of \Se structures. In particular, the $(2m+11)$-orbifolds
$\overbrace{S^3\star\cdots\star S^3}^{m\rm\;times} \star S_l\star
\cals(\grO_k)$ discussed in Corollary \sewz.18 satisfy $b_2=l+k+m+1$ and
$b_3=0.$

\def\lcm{{\rm lcm}}
\bigskip
\centerline{\bf \S 6. A Lattice of Einstein Orbifolds}
\medskip

In their paper [WZ] Wang and Ziller constructed Einstein metrics on the
total space of many torus bundles over K\"ahler-Einstein manifolds. The purpose
of the this section is to discuss the various Einstein manifolds that arise as
a result of our extension of their method to the orbifold category.

Given any two \Se manifolds $\cals_1,\cals_2$ one can define an orbifold with
an Einstein metric by generalizing the Wang-Ziller procedure to treat
circle V-bundles over K\"ahler-Einstein orbifolds. We simply relax the
condition in Proposition \sewz.6 that the $l_i$ be the respective relative
indices of $\cals_i.$ For any pair $(l,k)$ of positive integers we let
$S^1(l,k)$ denote the circle action on $\cals_1\times \cals_2$ given by
$(x,y)\mapsto (\grt^{k}x,\grt^{-l}y).$ Then we define orbifolds as the quotient
$$M(\cals_1,\cals_2;l,k)={\cals_1\times \cals_2\over S^1(l,k)}.
\leqno{\seoe.1}$$
It is clear that when $(l,k)$ are the relative indices $(l_1,l_2)$ we recover
our join operation, that is $M(\cals_1,\cals_2;l_1,l_2)=\cals_1\star \cals_2.$
The following result follows directly from our results together with those of
Wang and Ziller [WZ].

\noindent{\sc Theorem} \seoe.2: \tensl Let $\cals_i\in \cals\cale^s$ with
orders $m_i,$ respectively. Then
\item{(i)} $M(\cals_1,\cals_2;l,k)$ is a compact orbifold.
\item{(ii)} $\pi_1^{orb}(M(\cals_1,\cals_2;l,k))=0$ if and only if
$\gcd(l,k)=1.$
\item{(iii)} If $\gcd(m_1k,m_2l)=1$ then $M(\cals_1,\cals_2;l,k)$ is a compact
simply connected manifold.
\item{(iv)} $H^q(M(\cals_1,\cals_2;l,k),\bbq)\approx H^q(\cals_1\star
\cals_2,\bbq)$ for all $q$ and for all $l,k.$
\item{(v)} $M(\cals_1,\cals_2;l,k)$ admits a homothety class of
Einstein metrics of positive scalar curvature.
\item{(vi)} $M(\cals_1,\cals_2;l,k)$ has a Sasakian structure that is \Se if
and only if $(l,k)=m(l_1,l_2)$ for $m\in \bbz^+.$ \tenrm

\noindent
{\sc Proof}: (i) through (iv) and (vi) are clear from
the discussion in previous sections.
The homothety class of Einstein metrics on $M(\cals_1,\cals_2;l,k)$
in (v) is obtained as follows: One considers the product
$\cals_1\times\cals_2$ with a Riemannian metric given by
$g(\lambda_1,\lambda_2)=
\lambda_1g_1+\lambda_2g_2$, where $g_1,g_2$ are the Sasakian-Einstein
metrics on each factor. Then a calculation identical to that of [WZ]
shows that there exists a ratio $\lambda_1/\lambda_2$ for which
the quotient metric
$\hat{g}(\lambda_1,\lambda_2)$ on $M(\cals_1,\cals_2;l,k)$
obtained by Riemannian submersion
via the quotient map
$$\pi(\lambda_1,\lambda_2): \cals_1\times \cals_2
\longrightarrow M(\cals_1,\cals_2;l,k)$$
is Einstein. This ratio depends
on $(n_1,n_2)$,  $(l_1,l_2)$ and $(l,k)$.  In the Sasakian-Einstein case
of (vi) the scaling factors $(\lambda_1,\lambda_2)$ are in addition
completely determined by the condition that submerssed Einstein
metric on the associated leaf space
of the characteristic foliation be of correct Einstein constant
(see \sewz.4).
\hfill\za

The set $\{M(\cals_1,\cals_2;l,k)|(l,k)\in \bbz^+\times \bbz^+\}$ can be given
the structure of a lattice as follows:
\item{(1)} $M(\cals_1,\cals_2;l,k)\leq M(\cals_1,\cals_2;l',k')$ if and only if
$l\leq l'$ and $k\leq k';$
\item{(2)} $M(\cals_1,\cals_2;l,k)\vee M(\cals_1,\cals_2;l',k')
=M(\cals_1,\cals_2;\lcm(l,l'),\lcm(k,k'));$
\item{(3)} $M(\cals_1,\cals_2;l,k)\wedge M(\cals_1,\cals_2;l',k')
=M(\cals_1,\cals_2;\gcd(l,l'),\gcd(k,k')).$

\noindent This lattice is isomorphic to the lattice $\bbz^+\times \bbz^+$ with
the same operations.

Hence, associated to every pair $(\cals_1,\cals_2)$ of \Se manifolds there is a
lattice $L(\cals_1,\cals_2)$ of compact Einstein orbifolds, all of positive
scalar curvature of course, and all having the same rational cohomology. Notice
that if $\gcd(m_1,m_2)>1$ then no member of the lattice is a smooth manifold.
If however, $\gcd(m_1,m_2)=1$ then there will be infinitely many smooth
manifolds in the lattice.  There is a half-line of \Se orbifolds in the lattice
$L(\cals_1,\cals_2)$ given by $\{M(\cals_1,\cals_2;ml_1,ml_2)|m=1,2,\cdots \},$
precisely one of which is simply connected in the orbifold sense, namely
$M(\cals_1,\cals_2;l_1,l_2).$ We want to extend the lattice
$L(\cals_1,\cals_2)$ in two ways.  First we allow either $l$ or $k$ to take the
values $0,$ but not both simultaneously.  Then we see that
$M(\cals_1,\cals_2;l,0)= \cals_1\times \calz_2$ and $M(\cals_1,\cals_2;0,k)=
\calz_1\times \cals_2.$ Notice that these added elements do not necessarily
have the same rational cohomology as the elements of $L(\cals_1,\cals_2).$
Second and more importantly we want to include limit points of sequences.  To
do this we consider the Cheeger $\rho^*$-topology as discussed for example in
Wang and Ziller [WZ] by noticing that this easily extends to the orbifold
category.  This topology defines a distance between two Riemannian orbifolds
that basically depends on the curvature and its covariant derivative.  Thus,
two Riemannian orbifolds can be close even though their topology may be quite
different. Using Proposition 4.3 of [WZ] we find

\noindent{\sc Proposition} \seoe.3: \tensl Any element
$M(\cals_1,\cals_2;l,k)\in L(\cals_1,\cals_2)$ can be written as the limit
$$M(\cals_1,\cals_2;l,k)= \lim_{t\rightarrow
\infty}M(\cals_1,\cals_2;lt+a,kt+b)$$
for arbitrary $(a,b)\in \bbz^+\times \bbz^+.$ \tenrm

\noindent This proposition says that the completion
$\check{L}(\cals_1,\cals_2)$ of $L(\cals_1,\cals_2)$ in the Cheeger
$\grr^*$-topology behaves like the projective scheme over the lattice
$L(\cals_1,\cals_2).$

We are mainly interested in those points of the lattice that are represented by
smooth manifolds.  Necessary condition for the smoothness of the join
$\cals_1\star \cals_2$ are that $\gcd(m_i,l_i)=1$ and $\gcd(m_1,m_2)=1.$ Thus,
if $\gcd(m_1,m_2)=1$ there will be infinitely many pairs $(l,k)$ such that
$M(\cals_1,\cals_2;l,k)$ is a smooth Einstein manifold even though
$\cals_1\star \cals_2$ may not be smooth.  An example is given by taking
$\cals_1=S^{4n+1}$ and $\cals_2$ any 3-Sasakian 7-manifold whose order is even.
Then $S^{4n+1}\star \cals$ is not smooth, but $M(S^{4n+1},\cals;l,k)$ will be
smooth for infinitely many $(l,k),$ namely those that satisfy $\gcd(k,m_2l)=1.$

\noindent{\sc Example} \seoe.4: As a special case let us take $\cals_1$ to be
regular and $k=1.$ Then the gcd condition in (iii) above is automatically
satisfied, so for each $l,$ $M(\cals_1,\cals_2;l,1)$ is a smooth Einstein
manifold of dimension $\hbox{dim}~\cals_1 +\hbox{dim}~\cals_2 -1.$ Now let us
further specialize by putting $\cals_1=S^3$ and taking $\cals_2$ to be any
simply connected 3-Sasakian 7-manifold $\cals$ that is not $S^7.$ In this case
$M(S^3,\cals;1,1)=S^3\star \cals$ and we get a sequence
$\{M(S^3,\cals;l,1)\}_{l=1}^\infty$ of compact simply connected Einstein
manifolds of dimension $9.$ Furthermore, with respect to the Cheeger
$\grr^*$-topology  we see this sequence converges to $S^2\times \cals$ with the
product Einstein metric. One can also see that this sequence of Einstein
manifolds contains a subsequence of homotopically distinct Einstein manifolds.

\noindent{\sc Example} \seoe.5: In the case of $\cals_1 =\cals_l$ a circle
bundle over a del Pezzo surface, and $\cals_2 =\cals(\grO_k),$ we know that
$\cals_l\star \cals(\grO_k)$ will be a smooth manifold if and only if
$\hbox{Ord}(\cals(\grO_k))$ is odd. However, even if $\cals_l\star
\cals(\grO_k)$ is not a smooth manifold, $M(S_l,\cals(\grO_k);1,1)$ is always a
smooth Einstein manifold. Of course, it is not \Se. It is interesting to ponder
the question whether the continuous family of \Se structures on $S_l$ when
$5\leq l\leq 8$ induces an effective continuous family of Einstein metrics on
$M(S_l,\cals(\grO_k);1,1).$

An interesting question is whether any elements in a given lattice
$L(\cals_1,\cals_2)$ are homeomorphic or diffeomorphic. This is a difficult
question in general; however, Wang and Ziller give some examples in the case of
spheres. For example, it follows directly from Smale's classification of simply
connected 5-manifolds with spin [Sm] that $M(S^3,S^3;l,k)$ are all
diffeomorphic to $S^2\times S^3$ for all relatively prime $(l,k).$

\noindent{\sc Remark} \seoe.6: Clearly
Theorem \seoe.2 can be further extended to the
case of $n$ compact \Se orbifolds $(\cals_i,g_i), i=1,\ldots,n$.
Let ${\bf x}=(x_1\ldots,x_n)$
and consider the product
$\cals_1\times\cdots\times\cals_n$ with a metric
$g({\bf x})=x_1g_1+\cdots+x_ng_n$. Now, there is
a canonical $n$-torus action on $\cals_1\times\cdots\times\cals_n$.
A choice of an arbitrary circle subgroup in this torus is obtained via a
homomorphism $h({\bf p}):S^1\rightarrow T^n$,
where ${\bf p}=(p_1,\ldots, p_n)$, and
$p_i\in\bbz$ are the ``winding numbers" on each factor.
Now, the quotient of
$T^n$ by this circle gives a $T^{n-1}({\bf p})$-torus
action on the product $\cals_1\times\cdots\times\cals_n$
with the quotient space $M(\cals_1,\ldots,
\cals_n;{\bf p};{\bf x})$ being a compact orbifold. The combinatorial
conditions for smoothness are considerably more complicated here, but in
principle one certainly expects to find an analogue of (iii) of Theorem
\seoe.2.  In any case one can produce many examples of smooth quotients. As
in the $n=2$ case we do get an $n$-dimensional lattice
$L(\cals_1,\ldots,\cals_n)$ of compact orbifolds which all admit Einstein
metrics of positive scalar curvature by the argument identical to that of [WZ].
Also, there is a unique simply-connected Sasakian-Einstein space
$\cals_1\star\cdots\star\cals_n$ at some lattice point of $L$ determined by the
relative indices of the factors.  Hence, we can think of
$L(\cals_1,\ldots,\cals_n)$ as completely determined by
$\cals_1\star\cdots\star\cals_n$.

\noindent{\sc Remark} \seoe.7: Finally, even Wang and Ziller's torus bundle
construction can be reformulated and generalized in the language of the monoid
$(\cals\cale,\star)$. In the discussion
of the previous remark one can consider a homomorphism
$h(\Omega):T^s\rightarrow T^n$, where $\Omega\in\calm(s\times n;\bbz)$
is a ``weight matrix". This yields a $T^{n-s}(\Omega)$-torus
action on the product $\cals_1\times\cdots\times\cals_n$
with the quotient space $M(\cals_1,\ldots,
\cals_n;\Omega;{\bf x})$ being a compact orbifold and a $T^s$-bundle
over the product $\calz_1\times\cdots\times\calz_n$. Taking all weight
matrices $\Omega$
defines an $ns$-dimensional lattice $L(\cals_1,\ldots,\cals_n; s)$ whose
points are all compact Einstein orbifolds.

\medskip
\centerline{\bf Bibliography}
\medskip
\font\ninesl=cmsl9
\font\bsc=cmcsc10 at 10truept
\parskip=1.5truept
\baselineskip=11truept
\ninerm
\item{[AFHS]} {\bsc B. S. Acharya, J. M. Figueroa-O'Farrill, C. M. Hull, and
B. Spence}, {\ninesl Branes at Conical Singularities and Holography},
preprint, August 1998; hep-th/9808014.
\item{[An]} {\bsc M. Anderson}, {\ninesl Convergence and rigidity of manifolds
under Ricci curvature bounds}, Invent. Math. 102 (1990), 429-445.
\item{[Ba1]} {\bsc W. L. Baily}, {\ninesl The decomposition theorem for
V-manifolds}, Amer. J. Math. 78 (1956), 862-888.
\item{[Ba2]} {\bsc W. L. Baily}, {\ninesl On the imbedding of V-manifolds in
projective space}, Amer. J. Math. 79 (1957), 403-430.
\item{[Bes]} {\bsc A. L. Besse}, {\ninesl Einstein Manifolds},
Springer-Verlag, New York (1987).
\item{[BFGK]} {\bsc H. Baum, T. Friedrich, R. Grunewald, and I. Kath},
{\ninesl Twistors and Killing Spinors on Riemannian Manifolds},
Teubner-Texte f\"ur Mathematik, vol. 124, Teubner, Stuttgart, Leipzig, 1991.
\item{[BG1]} {\bsc C. P. Boyer and  K. Galicki}, {\ninesl
The twistor space of a 3-Sasakian manifolds},
Int. J. Math. 8 (1997), 31-60.
\item{[BG2]} {\bsc C. P. Boyer and  K. Galicki}, {\ninesl
3-Sasakian Manifolds}, preprint, January 1998, MPI-98/19;
(http://www.mpim-bonn.mpg.de/html/preprints/preprints.html);
to appear in {\it Essays on Einstein Manifolds}, International Press;
C. LeBrun and M. Wang, Eds.
\item{[BGM1]} {\bsc C. P. Boyer, K. Galicki, and B. M. Mann}, {\ninesl
Quaternionic reduction and Einstein manifolds}, Comm.
Anal. Geom. 1 (1993), 1-51.
\item{[BGM2]} {\bsc C. P. Boyer, K. Galicki, and B. M. Mann}, {\ninesl
The geometry and topology of 3-Sasakian manifolds},
J. reine angew. Math. 455 (1994), 183-220.
\item{[BGM3]} {\bsc C. P. Boyer, K. Galicki, and B. M. Mann}, {\ninesl
On strongly inhomogeneous Einstein manifolds}, Bull. London Math. Soc. 28
(1996), 401-408.
\item{[BGM4]} {\bsc C. P. Boyer, K. Galicki, and B. M. Mann}, {\ninesl
Hypercomplex structures on Stiefel manifolds}, Ann. Global Anal. Geom. 14
(1996), 81-105.
\item{[BGM5]} {\bsc C. P. Boyer, K. Galicki, and B. M. Mann}, {\ninesl
A note on smooth toral reductions of spheres},
Manuscripta Math. 95 (1998), 149-158.
\item{[BGM6]} {\bsc C. P. Boyer, K. Galicki, and B. M. Mann}, {\ninesl
Hypercomplex structures from 3-Sasakian structures}, J. reine angew. Math. 501
(1998), 115-141.
\item{[BGMR]} {\bsc C. P. Boyer, K. Galicki, B. M. Mann, and E. Rees},
{\ninesl Compact 3-Sasakian 7-manifolds with arbitrary second Betti number},
Invent. Math. 131 (1998), 321-344.
\item{[Bl]} {\bsc D. E. Blair}, {\ninesl Contact Manifolds in
Riemannian Geometry}, Lecture Notes in Mathematics 509, Springer-Verlag,
New Yrok 1976.
\item{[Bo]} {\bsc A. Borel}, {\ninesl K\"ahlerian coset spaces of semi-simple
Lie groups}, Proc. Nat. Acad. Sci. USA 40 (1954), 1147-1151. Reprinted in A.
Borel, Oeuvres, Collected Papers, Vol. 1, Springer-Verlag, Berlin Heidelberg,
1983.
\item{[BR]} {\bsc A. Borel and R. Remmert}, {\ninesl \"Uber kompakte homogene
K\"ahlersche Mannigfaltigkeiten}, Math. Ann. 145 (1962), 429-439.
\item{[BW]} {\bsc W.M. Boothby and H.C. Wang}, {\ninesl On Contact Manifolds},
Ann. of Math. 68 (1958), 721-734.
\item{[Dim]} {\bsc A. Dimca}, {\ninesl Singularities and Topology of
Hypersurfaces}, Springer-Verlag, New York, 1992.
\item{[DO]} {\bsc S. Dragomir and L. Ornea}, {\ninesl Locally Conformal
K\"ahler Geometry}, Progress in Mathematics 155, Birkh\"auser, Boston, 1998.
\item{[DT]} {\bsc W. Ding and G. Tian}, {\ninesl K\"ahler-Einstein metrics and
the generalized Futaki invariants}, Invent. Math. 110 (1992), 315-335.
\item{[Esch1]} {\bsc J. H. Eschenburg}, {\ninesl New examples of manifolds
with strictly positive curvature}, Invent. Math. 66 (1982), 469-480.
\item{[Esch2]} {\bsc J. H. Eschenburg}, {\ninesl Cohomology of biquotients},
Manuscripta Math. 75 (1992), 151-166.
\item{[Fi]} {\bsc J. M. Figueroa-O'Farrill}, {\ninesl Near-Horizon
Geometries of Supersymmetric Branes}, preprint, July 1998; hep-th/9807149.
\item{[FrKat1]} {\bsc T. Friedrich and I. Kath}, {\ninesl Einstein
manifolds of dimension five with small first eigenvalue of the Dirac operator},
J. Diff. Geom. 29 (1989), 263-279.
\item{[FrKat2]} {\bsc T. Friedrich and I. Kath}, {\ninesl Compact
seven-dimensional
manifolds with Killing spinors}, Comm. Math. Phys. 133 (1990), 543-561.
\item{[Gro]} {\bsc M. Gromov}, {\ninesl Curvature, diameter and Betti numbers},
Comment. Math. Helvetici 56 (1981) 179-195.
\item{[GS]} {\bsc K. Galicki and S. Salamon}, {\ninesl On
Betti numbers of 3-Sasakian manifolds},  Geom. Ded. 63 (1996), 45-68.
\item{[Hae]} {\bsc A. Haefliger}, {\ninesl Groupoides d'holonomie et
classifiants}, Ast\'erisque 116 (1984), 70-97.
\item{[Hat]} {\bsc Y. Hatakeyama}, {\ninesl Some notes on differentiable
manifolds with almost contact structures}, T\^ohuku Math. J. 15 (1963),
176-181.
\item{[HS]} {\bsc I. Hasegawa and M. Seino}, {\ninesl Some remarks on Sasakian
geometry--applications of Myers' theorem and the canonical affine connection},
J. Hokkaido Univ. Education 32 (1981), 1-7.
\item{[Kaw]} {\bsc T. Kawasaki}, {\ninesl The Riemann-Roch theorem for
complex
V-manifolds}, Osaka J. Math. 16 (1979), 151-159.
\item{[KO]} {\bsc S. Kobayashi and T. Ochiai}, {\ninesl Characterizations of
complex projective spaces and hyperquadrics}, J. Math.  Kyoto. Univ. 13 (1973),
31-47.
\item{[Kob1]} {\bsc S.  Kobayashi}, {\ninesl On compact K\"ahler manifolds with
positive Ricci tensor}, Ann. Math. 74 (1961), 381-385.
\item{[Kob2]} {\bsc S.  Kobayashi}, {\ninesl Topology of positively pinched
Kaehler manifolds}, T\^ohoku Math. J. 15 (1963), 121-139.
\item{[Kol]} {\bsc J. Koll\'ar}, {\ninesl Rational Curves on Algebraic
Varieties}, Springer-Verlag, New York, 1996.
\item{[KS1]} {\bsc N. Koiso and Y. Sakane}, {\ninesl Non-homogeneous
K\"ahler-Einstein metrics on compact complex manifolds} in Curvature and
Topology of Riemannian manifolds, Lecture Notes in Mathematics 1201,
Springer-Verlag, 1986, 165-179.
\item{[KS2]} {\bsc N. Koiso and Y. Sakane}, {\ninesl Non-homogeneous
K\"ahler-Einstein metrics on compact complex manifolds II}, Osaka J. Math. 25
(1988), 933-959.
\item{[KW]} {\bsc I. R. Klebanov and E. Witten}, {\ninesl
Superconformal Field Theory on Threebranes at a Calabi-Yau Singularity},
preprint, July 1998; hep-th/9807080.
\item{[Ma]} {\bsc J. Maldacena}, {\ninesl The large $N$ limit of
superconformal field theories and supergravity}, preprint, 1997;
hep-th/9711200.
\item{[Mab]} {\bsc T. Mabuchi}, {\ninesl Einstein-K\"ahler forms, Futaki
invariants and convex geometry on toric Fano varieties}, Osaka J. Math. 24
(1987), 705-737.
\item{[Mol]} {\bsc P. Molino}, {\ninesl Riemannian Foliations}, Progress in
Mathematics 73, Birkh\"auser, Boston, 1988.
\item{[MiMo]} {\bsc Y. Miyaoka and S. Mori}, {\ninesl A numerical criterion
for uniruledness}, Ann.  Math. 124 (1986), 65-69.
\item{[MoMu]} {\bsc S. Mori and S. Mukai}, {\ninesl Classification of Fano
3-folds with $\scriptstyle{B_2\geq 2.}$}, Manuscripta Math. 36 (1981), 147-162.
\item{[MP]} {\bsc D. R. Morrison and M. R. Plesser},
{\ninesl Non-Spherical Horizons, I}, preprint, October 1998;
hep-th/9810201.
\item{[Na]} {\bsc A.M. Nadel} {\ninesl Multiplier ideal sheaves and existence
of K\"ahler-Einstein metrics of positive scalar curvature}, Ann. Math. 138
(1990), 549-596.
\item{[OhTa]} {\bsc K. Oh and R. Tatar}, {\ninesl Three Dimensional
SCFT from M2 Branes at Conical Singularities}, preprint, October 1998;
hep-th/9810244.
\item{[PP]} {\bsc H. Pedersen and Y. S. Poon}, {\ninesl
A note on rigidity of
3-Sasakian manifolds}, Proc. Am. Math. Soc., to appear (1998).
\item{[Sa]} {\bsc Y. Sakane}, {\ninesl Examples of compact Einstein K\"ahler
manifolds with positive Ricci tensor}, Osaka J. Math. 23 (1986), 585-617.
\item{[Sas]} {\bsc  S. Sasaki}, {\ninesl On differentiable manifolds with
certain structures which are closely related to almost contact structure},
T\^ohoku Math. J. 2 (1960), 459-476.
\item{[Siu]} {\bsc Y.-T. Siu}, {\ninesl The existence of K\"ahler-Einstein
metrics on manifolds with positive anticanonical line bundle and a suitable
finite symmetry group}, Ann. Math. 127 (1988), 585-627.
\item{[Sm]} {\bsc S. Smale}, {\ninesl On the structure of 5-manifolds},
Ann. Math. 75 (1962), 38-46.
\item{[Tan]} {\bsc S. Tanno}, {\ninesl Geodesic flows
on $C_L$-manifolds and Einstein metrics on $S^3\times S^2$}, in
{\it Minimal submanifolds and geodesics (Proc.
Japan-United States Sem., Tokyo, 1977)}, pp. 283-292, North Holland,
Amsterdam-New York, 1979.
\item{[Thu]} {\bsc W. Thurston}, {\ninesl The Geometry and Topology of
3-Manifolds}, Mimeographed Notes, Princeton Univ. Chapt. 13 (1979).
\item{[Tho]} {\bsc C.B. Thomas}, {\ninesl Almost regular contact manifolds},
J. Differential Geom. 11 (1976), 521-533.
\item{[Ti1]} {\bsc G. Tian}, {\ninesl On K\"ahler-Einstein metrics on certain
K\"ahler manifolds with $C_1(M)>0$}, Invent. Math. 89 (1987), 225-246.
\item{[Ti2]} {\bsc G. Tian}, {\ninesl On Calabi's Conjecture for complex
surfaces with positive first Chern class}, Invent. Math. 101 (1990), 101-172.
\item{[Ti3]} {\bsc G. Tian}, {\ninesl K\"ahler-Einstein metrics with positve
scalar curvature}, Invent. Math. 137 (1997), 1-37.
\item{[TY]} {\bsc G. Tian and S.-T. Yau}, {\ninesl K\"ahler-Einstein
metrics on complex surfaces with $c_1>0$}, Comm. Math. Phys. 112 (1987),
175-203.
\item{[Vai]} {\bsc I. Vaisman}, {\ninesl Locally conformal K\"ahler manifolds
with parallel Lee form}, Rend. Matem. Roma 12 (1979), 263-284.
\item{[WZ]} {\bsc M. Wang and W. Ziller}, {\ninesl Einstein metrics on
principal torus bundles}, J. Diff. Geom. 31 (1990), 215-248.
\item{[Wi1]} {\bsc J.A.  Wi\'sniewski}, {\ninesl On a conjecture of Mukai},
Manus. Math. 68 (1990), 135-141.
\item{[Wi2]} {\bsc J.A.  Wi\'sniewski}, {\ninesl On Fano manifolds of large
index}, Manus. Math. 70 (1991), 145-152.
\item{[Wit]} {\bsc E. Witten}, {\ninesl Anti-de Sitter space and holography},
preprint, 1998; hep-th/9802150.
\item{[YKon]} {\bsc K. Yano and M. Kon}, {\ninesl
Structures on manifolds}, Series in Pure Mathematics 3,
World Scientific Pub. Co., Singapore, 1984.
\medskip
\bigskip \line{ Department of Mathematics and Statistics
\hfil October 1999} \line{ University of New Mexico \hfil} \line{ Albuquerque,
NM 87131 \hfil } \line{ email: cboyer@math.unm.edu, galicki@math.unm.edu \hfil}
\bye